\documentclass[reqno]{amsart}
\usepackage{amscd}
\usepackage[dvips]{graphics}
\usepackage{color}
\usepackage{graphicx,subfigure,float,caption2,epstopdf}

\def\beq{\begin{equation}}
\def\eeq{\end{equation}}
\def\ba{\begin{array}}
\def\ea{\end{array}}

\def\cal{\mathcal}

\newenvironment{abs}{\textbf{Abstract}\mbox{  }}{ }
\newenvironment{key words}{\textbf{Keywords}\mbox{  }}{ }
\newtheorem{theorem}{Theorem}[section]
\newtheorem*{theorem*}{Theorem A }
\newtheorem{lemma}[theorem]{Lemma}
\newtheorem{prop}[theorem]{Proposition}
\newtheorem{crl}[theorem]{Corollary}

\theoremstyle{definition}
\renewenvironment{proof}{\noindent{\textbf{Proof.}}}{\hfill$\Box$}
\newtheorem{rem}[theorem]{Remark}

\theoremstyle{remark}
 
\numberwithin{equation}{section}

\allowdisplaybreaks

\begin{document}

\title[Stein-Weiss type inequalities on the upper half space]{\textbf{Stein-Weiss type inequalities with partial variable weight on the upper half space and related  weighted inequalities}}
\author{Jingbo Dou and Jingjing Ma}

\address{Jingbo Dou, School of Mathematics and Statistics,
Shaanxi Normal University, Xi'an, 710119, P. R. China}
\email{jbdou@snnu.edu.cn}

\medskip

\address{Jingjing Ma, School of Mathematics and Statistics,
Shaanxi Normal University,
Xi'an, 710119, P. R. China}
\email{mjjwnm@126.com}

\maketitle

% ----------------------------------------------------------------
\noindent
\begin{abs}
In this paper, we establish a class of Stein-Weiss type inequality with partial variable weight functions on the upper half space using a weighted Hardy type inequality. Overcoming the impact of  weighted functions, the existence of extremal functions is proved via the concentration compactness principle, whereas Riesz rearrangement inequality is not available. Moreover, the cylindrical symmetry with respect to $t$-axis and the explicit forms on the boundary of all nonnegative extremal functions are discussed via the method of moving planes and  method of moving spheres, as well as, regularity results are obtained by the regularity lift lemma and bootstrap technique. As applications, we obtain some weighted Sobolev inequalities with partial variable weight function for Laplacian and fractional Laplacian.
\end{abs}\\
\begin{key words} Hardy-Littlewood-Sobolev inequality, Stein-Weiss type inequality, weighted Sobolev inequality, fractional integral operator, extremal function, cylindrical symmetry, regularity, moving sphere method.
 \end{key words}
\indent

%---------------------------------------------------------------------------------
\section{Introduction\label{Section 1}}
The classical sharp Hardy-Littlewood-Sobolev (HLS) inequality states as
\begin{equation}\label{class-HLS}
\big|\int_{\mathbb{R}^n}\int_{\mathbb{R}^n}f(x)|x-y|^{-\lambda}g(y)dxdy\big|\le N(n, p,\lambda)\|f\|_{L^p({\mathbb{R}^n})}\|g\|_{L^t({\mathbb{R}^n})}
\end{equation}for all $f\in L^p(\mathbb{R}^n), \,  g\in L^t(\mathbb{R}^n), \, 1<p, \, t<\infty, \,$ $ 0<\lambda<n$ and $1/p+1/t+\lambda/n=2$ (see e.g., \cite{HL1928, HL1930, So1963, Stein}).
Lieb \cite{Lieb1983} proved the existence of extremal functions to inequality \eqref{class-HLS} by rearrangement technique and computed the best constants for conformal case $t=p,$  and $p=2$ or $t=2$. Moreover, Chen, Ou and Li \cite{CLO2006} and Li \cite{L2004} classified the all positive solutions to the extremal functions for $t=p$. Frank and Lieb \cite{FL2012} gave a new proof of the best constant to inequality \eqref{class-HLS} by a rearrangement free technique.

Recall the following basic estimate
\begin{equation}\label{Tr-0}
|u(x)|\le\frac{1}{\omega_{n-1}}\int_{\mathbb{R}^{n}}\frac{|\nabla u(y)|}{|x-y|^{n-1}}dy
\end{equation}holds for  any $u\in C^\infty_0(\mathbb{R}^{n})$, where $\omega_{n-1}=\frac{2\pi^{\frac{n}2}}{\Gamma(\frac{n}2)}$ is the surface area of the unit sphere $\mathbb{S}^{n-1}.$ Employing \eqref{class-HLS} and \eqref{Tr-0}, Sobolev \cite{So1963} established classical Sobolev inequality (also see Aubin \cite{A1998}) as follows
\begin{equation}\label{Sobolev-inq}
\big(\int_{\mathbb{R}^{n}}|u(x)|^{\frac{np}{n-p}}dx\big)^\frac{n-p}{n}\le S_{n,p}\int_{\mathbb{R}^{n}} |\nabla u(x)|^pdx
\end{equation}for $u\in W^{1,p}(\mathbb{R}^{n}),$ where $1<p<n$, $ S_{n,p}$ is the best constant. The inequality \eqref{Sobolev-inq} plays an important role in PDE and functional analysis, etc.

The sharp HLS inequality is also closely related to Moser-Trudinger-Onofri and Beckner inequalities \cite{B1993}, as well as Gross's logarithmic Sobolev inequality \cite{Gr}. All these inequalities play a significant role in solving global geometric problems, such as Yamabe problem, Ricci flow problem, etc.. In recent years, from global analysis, the integral (curvature) equation involving the HLS inequality \eqref{class-HLS} was studied by Zhu \cite{Z2016}, Dou and Zhu \cite{DZ2019} and Dou, Guo and Zhu \cite{DGZ2020}.

Stein and Weiss \cite{SW1958} established the following weighted HLS inequality (also called Stein-Weiss inequality)
\begin{equation}\label{class-WHLS}
\big|\int_{\mathbb{R}^n}\int_{\mathbb{R}^n}\frac{f(x) g(y)}{|x|^\beta|x-y|^{\lambda}|y|^\alpha}dx dy\big|\le N(n,p,\alpha,\beta,\lambda)\|f\|_{L^p({\mathbb{R}^n})}\|g\|_{L^t{(\mathbb{R}^n)}}
\end{equation}for all $f\in L^p(\mathbb{R}^n),\, g\in L^t(\mathbb{R}^n), 1<p, \, t<\infty, 0<\lambda<n, \alpha<n/{p'}, \beta<n/{t'}, \alpha+\beta\ge0$ and
$1/p+1/t+(\alpha+\beta+\lambda)/n=2.$ Throughout this paper and in here, we always let $p'$ denote the conjugate index of $p.$

Under the restriction $\alpha, \, \beta\ge0,$ Lieb \cite{Lieb1983} discussed the existence of extremal functions to inequality \eqref{class-WHLS} by rearrangement technique. Later, Chen, Lu and Tao \cite{CLT2019} employed the concentration compactness principle to show the existence of extremal  functions to HLS inequality \eqref{class-WHLS} and HLS inequality on the Heisenberg group respectively, which dropped out the restriction $\alpha, \, \beta\ge0.$ However, it is an open question to compute the best constant and explicit extremal functions to inequality \eqref{class-WHLS} in the general case due to the impact of the weighted functions.

Chen and Li \cite{CL2008} classified the extremal functions to inequality \eqref{class-WHLS} and showed the best constant for a special case $\lambda=2, \,  \alpha=\beta=\frac{n}{p+1}+\frac{n-2}2$ and $p=t.$ Beckner \cite{B2008a,B2008b} discussed the best constant of the inequality \eqref{class-WHLS} with $\alpha+\beta+\lambda=n$ and $p=t'.$ Wu, Shi and Yan \cite{WSY2014} showed the best constant of the inequality \eqref{class-WHLS} for three cases: $(i)$ $p=q=1$; $(ii)$ $p=1$ and $1<t\le\infty$ or $t=1$ and $1<p\le\infty$; $(iii)$ $1<p,t<\infty$ and $\frac1p+\frac1t=1$. Jin and Li \cite{JL2006} discussed the symmetry of the positive extremal functions to inequality \eqref{class-WHLS} by using the method of moving planes in integral forms. Moreover, the asymptotic behavior of extremal functions to inequality \eqref{class-WHLS}  was also analyzed around the origin and near infinity with $\alpha, \, \beta\ge0$ in \cite{BLL2011, LLM2011, LM2011, LL2007, O2012}.

In the last two decades, various extensions of HLS inequality have been investigated. For examples, one has sharp HLS inequality  on the upper half space, on Heisenberg group, on compact Riemannian manifolds, on reversed forms, and on weighted forms, see e.g., \cite{CL2018, D2016, DGZ2017, DZ2015a, DZ2015b, G2020, HLZ2012, HZ2016}. In particular, we notice that  Beckner \cite{B2021} established  the following  Stein-Weiss type inequality 
\begin{equation}\label{WHLSD-Beckner}
\big|\int_{\mathbb{R}^{n+1}_+}\int_{\mathbb{R}^{n+1}_+}\frac{f(x,t) g(y,z)}{t^\alpha|(x,t)-(y,z)|^{\lambda} z^\alpha} dxdtdydz\big|
\le C_{\alpha,\lambda,p} \|f\|_{L^p(\mathbb{R}_+^{n+1})}\|g\|_{L^p(\mathbb{R}_+^{n+1})}
\end{equation}
for  $f, g\in L^p(\mathbb{R}_+^{n+1}), $ where $0<\lambda<n+1,0<\alpha,p=1-\frac{\lambda+2\alpha}{2n},$ $C_{\alpha,\lambda,p}$ is the positive constant.
And he showed \eqref{WHLSD-Beckner} is equivalent to a class of HLS inequality on hyperbolic space $\mathbb{H}^n$.  One of our motivations, we extend the inequality   \eqref{WHLSD-Beckner} to  the general case  and discuss the existence of extremal functions.

On the other hand, recently, first author, Sun, Wang and Zhu \cite{DSWZ2020} established a weighted Sobolev inequality on the upper half space as follows. Let $\mathbb{R}^{n+1}_+=\{(x,t)\in \mathbb{R}^n\times \mathbb{R}\,:\, t>0\}$ be the upper half space, for $n\ge 1, \, \beta>-1, \, \alpha>0,$ $\frac{n-1}{n+1}\beta\le\alpha\le\beta+2, \, p^*=\frac{2n+2\beta+2}{n+\alpha-1}.$ There is a positive constant $C_{n+1,\alpha,\beta}>0$ such that  for all $u\in{\cal D}_{\alpha}^{1,2}({\mathbb{R}^{n+1}_+})$
\begin{equation}\label{GGN-2}
\big(\int_{\mathbb{R}^{n+1}_+ }t^\beta |u|^{p^*}dydt\big)^{\frac2{p^*}} \le  C_{n+1,\alpha,\beta} \int_{\mathbb{R}^{n+1}_+ }t^{\alpha}|\nabla u|^2 dydt.
\end{equation}
They proved the best constant $C_{n+1, \alpha,\beta}$ can be achieved by a nonnegative extremal function $u\in{\cal D}_{\alpha}^{1,2}({\mathbb{R}^{n+1}_+})$ for $\frac{n-1}{n+1}\beta<\alpha<\beta+2,$ and the extremal functions were classified. In particular, two cases ($\alpha=\beta$ and $\beta=\alpha-1$) of the solutions were explicitly written out. In fact, inequality \eqref{GGN-2} still holds for $\alpha\le0$, see e.g., \cite{CL2007, WZ2022} and \cite{Maz2011} on $\mathbb{R}^{n+1}$. However,  we could not show  that the inequality \eqref{GGN-2} holds for case $\alpha<0 $ due to some reasons of technique.

 It is obvious that  the weighted function in inequality \eqref{WHLSD-Beckner} is just a partial variable weight $t$. Moreover, HLS inequality is closely relate to Sobolev inequality.
A nature question: Could we establish the inequality \eqref{GGN-2} with $\alpha<0$ employing the Stein-Weiss type inequality with partial variable weight $t$ ?

Motivated by above problems, in this paper, we investigates the extension of Stein-Weiss type inequality to a class of HLS inequality with partial variable weight $t$ on $\mathbb{R}^{n+1}_+$, and gives an affirmative answer to the above question.

  Our main results state as follows.
\begin{theorem}\label{WHLSB}
Let $\alpha,\beta,\lambda,p,r$ satisfy
\begin{equation}\label{WH-exp-0}
\begin{cases}
&n\ge1,0<\lambda<n+1, \, 1< r, \, p<\infty ,\\
&\alpha<\frac{1}{p'}, \beta<\frac {1}{r'},~\alpha+\beta\ge0, \\
&\frac1p +\frac1r+\frac{\alpha+\beta+\lambda}{n+1}=2.
\end{cases}
\end{equation}
Then there is a positive constant $N_{\alpha,\beta,\lambda,p}:=C(n,\alpha,\beta,\lambda,p)$ depending on $n,\alpha,\beta,\lambda$ and $p$, such that for all $f\in L^p(\mathbb{R}_+^{n+1}), \, g\in L^r( \mathbb{R}_+^{n+1}),$
\begin{equation}\label{WHLSD-O}
\big|\int_{\mathbb{R}^{n+1}_+}\int_{\mathbb{R}^{n+1}_+}\frac{f(x,t) g(y,z)}{t^\alpha|(x,t)-(y,z)|^{\lambda} z^\beta} dxdtdydz\big|
\le N_{\alpha,\beta,\lambda,p} \|f\|_{L^p(\mathbb{R}_+^{n+1})}\|g\|_{L^r(\mathbb{R}_+^{n+1})}
\end{equation}holds. Moreover, the upper and lower bounds of constant $N_{\alpha,\beta,\lambda,p}$ is
\begin{eqnarray*}
\max\big\{D_1,D_2,D_3\big\} 
&\le& N_{\alpha,\beta,\lambda,p}
\le\min\big\{\big(\frac{p}{p-1} \big)^{\frac{p-1}{p}}{p}^{\frac1p}D_1,\big(\frac{p}{p-1} \big)^{\frac{p-1}{p}}{p}^{\frac1p}D_2\big\},
\end{eqnarray*}
where
\begin{eqnarray*}
%D_0&=&{\color{blue}\big(\frac{w_n}n\big)^{\frac{r-1}r},\quad N=N(n,p,\lambda),} \\
D_1&=&\big[\frac{(r-1)\cdot\pi^{\frac n2}}{r(n+1-\beta-\lambda)-(n+1)}\cdot\frac{\Gamma\big(\frac{r(1-\beta )-1}{2(r-1)}\big)}{\Gamma\big(\frac{r(n+1-\beta)-(n+1)}{2(r-1 )}\big)}\big]^{\frac{r-1}{r}}\\
\quad&&\times\big[\frac{(p-1)\cdot\pi^{\frac n2}}{p(n+1-\alpha) -(n+1)}\cdot\frac{\Gamma\big(\frac{p(1-\alpha)-1}{2(p-1)}\big)}{\Gamma\big(\frac{p(n+1-\alpha )-(n+1)}{2(p-1)}\big)}\big]^{\frac{p-1}{p}},\\
D_2&=&\big[\frac{(r-1)\cdot\pi^{\frac n2}}{r(n+1-\beta)-(n+1)}\cdot\frac{\Gamma\big(\frac{r(1-\beta) -1}{2(r-1)}\big)}
{\Gamma\big(\frac{r(n+1-\beta)-(n+1)}{2(r-1)}\big)}\big]^{\frac{r-1}{r}}\\
\quad&&\times\big[\frac{(p-1)\cdot\pi^{\frac n2}}{p(n+1-\alpha-\lambda)-(n+1)}\cdot\frac{\Gamma\big(\frac{p(1-\alpha)-1}{2(p-1)}\big)}{\Gamma\big(\frac{p(n+1-\alpha )-(n+1)}{2(p-1)}\big)}\big]^{\frac{p-1}{p}},\\
D_3&=&\max \big\{\big[(2^n-1)\frac{w_{n-1}}n\big]^{\frac{r-1}r}, N(n+1,p,\lambda)\big[(2^n-1)\frac{w_{n-1}}n\big]^{\frac{r-1}r},\\
&&2^{\lambda}\big[(2^n-\frac1{2^{n+2}})\frac{w_n}{n+1}\big]^{\frac{p-1}{p}}
\big[(2^n-1)\frac{w_{n-1}}n\big]^{\frac{r-1}r}\big\},
\end{eqnarray*}
and $N(n+1,p,\lambda)$ is the best constant of  inequality \eqref{class-HLS} in $\mathbb{R}^{n+1}$  with the upper bound (see \cite{EM2001})
 \[
N(n+1,p,\lambda)\le\frac{n+1}{pr(n+1-\lambda)}\big(\frac{w_n}{n+1}\big)^{\frac \lambda {n+1}}\big[\big(\frac{p\lambda }{(n+1)(p-1)}\big)^{\frac \lambda{n+1}}+\big(\frac{r \lambda}{(n+1)(r-1)}\big)^{\frac \lambda {n+1}}\big]
\]
with $w_n=\frac{2\pi^{\frac{n+1}2}}{\Gamma(\frac{n+1}2)}$ is the surface area of unit ball in $\mathbb{R}^{n+1}$.
\end{theorem}
Define the singular integral operator on $\mathbb{R}^{n+1}_+$ as
\[
 E_\lambda f:=E_\lambda f(y,z)=\int_{\mathbb{R}^{n+1}_+}\frac{f(x,t)}{|(x,t)-(y,z)|^\lambda}dxdt,   ~~~~~~~~~~\forall\,(y,z)\in\mathbb{R}^{n+1}_+,~ \lambda \in(0,n+1).
\]
Let $q=r'$, we can write the assumptions of exponents $\alpha,\beta,\lambda,p,r$ in \eqref{WH-exp-0} as
\begin{equation}\label{WH-exp-1}
\begin{cases}
&n\ge1, 0<\lambda<n+1, \, 1< p\le q<\infty,\\
&\alpha<\frac{1}{p'}, \, \beta<\frac1q,~\alpha+\beta\ge0,  \\
&\frac1q=\frac1p-\frac{n+1-(\alpha+\beta+\lambda)}{n+1},
\end{cases}
\end{equation}
and write
\[\|t^\alpha f\|_{L^p(\mathbb{R}^{n+1}_+)}:=\big(\int_{\mathbb{R}^{n+1}_+}t^{\alpha p}|f|^pdxdt\big)^{\frac1p}.\]
The inequality \eqref{WHLSD-O} is equivalent to the following dual form.
 %to the following form with $q=r'$.%
\begin{theorem}\label{WHLSD-theo} Let $\alpha,\beta,\lambda,p,q$ satisfy the conditions \eqref{WH-exp-1}. Then there is a positive constant $N_{\alpha,\beta,\lambda,p}$  such that
\begin{equation}\label{WHLSD-1}
\|z^{-\beta}E_\lambda f\|_{L^q(\mathbb{R}^{n+1}_+)}\le N_{\alpha,\beta,\lambda,p}\|t^\alpha f\|_{L^p(
\mathbb{R}^{n+1}_+)}
\end{equation}holds for any $f\in L^p({\mathbb{R}^{n+1}_+}).$
\end{theorem}
\begin{rem}\label{Extremal remark} $(i)$ Obviously, \eqref{WHLSD-O} is also equivalent to
\begin{equation}\label{WHLSD-O-1}
\big|\int_{\mathbb{R}^{n+1}_+}\int_{\mathbb{R}^{n+1}_+}\frac{f(x,t) g(y,z)}{|(x,t)-(y,z)|^{\lambda}}dxdtdydz\big|
\le N_{\alpha,\beta,\lambda,p}\|t^\alpha f\|_{L^p(\mathbb{R}^{n+1}_+)}\|z^\beta g\|_{L^r(\mathbb{R}^{n+1}_+)}.
\end{equation}

$(ii)$ Define a weighted singular integral operator as
\[K_{\alpha,\beta}f(y,z):=\int_{\mathbb{R}^{n+1}_+}\frac{f(x,t)}{t^\alpha |(x,t)-(y,z)|^{\lambda}z^\beta}dxdt.
\]
Then inequality \eqref{WHLSD-1} is equivalent to the following form
\begin{equation}\label{WHLSD-2}
\| K_{\alpha,\beta} f\|_{L^q( \mathbb{R}^{n+1}_+)}\le N_{\alpha,\beta,\lambda,p}\|f\|_{L^p(
\mathbb{R}^{n+1}_+)}
\end{equation} with the same assumptions of Theorem  \ref{WHLSD-theo}.

$(iii)$ For the conformal case $p= p_\alpha=\frac{2(n+1)}{2(n+1)-\lambda-2\alpha}, \, r=r_\beta=\frac{2(n+1)}{2(n+1)-\lambda-2\beta}$, by Kelvin transformation, the inequality \eqref{WHLSD-O} is equivalent  to  the  following  weighted inequality on unit ball $B^{n+1}=B_1(x^1)\subset\mathbb{R}^{n+1}$ (see Corollary \ref{cor-ball})
\begin{eqnarray*}
& &\big|\int_{B^{n+1}}\int_{B^{n+1}}\big(\frac 2{1-|\xi-x^1|^2}\big)^{\alpha}\big(\frac 2{1-|\eta-x^1|^2}\big)^{\beta}
\frac{F(\xi)G(\eta)}{|\xi-\eta|^\lambda}d\xi d\eta\big|\nonumber\\
&\le& N_{\alpha,\beta,\lambda,p_\alpha}\|F\|_{L^{p_\alpha}(B^{n+1})}\|G\|_{L^{r_\beta}(B^{n+1})}
\end{eqnarray*} for any $F\in{L^{p_\alpha}(B^{n+1})}$ and $G\in{L^{r_\beta}(B^{n+1})},$ where $x^1=(0,\cdots,0,-1)\in\mathbb{R}^{n+1}.$
\end{rem}
The best constant in \eqref{WHLSD-O} or \eqref{WHLSD-1}  can be classified by
\begin{eqnarray*}
N_{\alpha,\beta,\lambda,p}
&=&\sup\big\{\int_{\mathbb{R}^{n+1}_+}\int_{\mathbb{R}^{n+1}_+}\frac{f(x,t) g(y,z)}{t^\alpha |(x,t)-(y,z)|^{\lambda}z^\beta}dxdtdydz\,:\,\\
& &\|f\|_{L^p(\mathbb{R}_+^{n+1})}=\|g\|_{L^r(\mathbb{R}_+^{n+1})}=1\big\} \\
&=&\sup\big\{\|z^{-\beta}E_\lambda f\|_{L^q(\mathbb{R}^{n+1}_+)}\,:\, f\in
L^p{(\mathbb{R}^{n+1}_+}),\,\|t^\alpha f\|_{L^p(\mathbb{R}^{n+1}_+)}=1\big\}.
%&=&\sup\big\{\|S f\|_{L^q(\mathbb{R}^{n+1}_+)}\,:\, \|f\|_{L^p(\mathbb{R}^{n+1}_+)}=1\big\}.~~~~~~~~~~~~~~~~~~~~~~~~~~~~~~~~~~~~~~~~~~~~~~~~~~~~~~~~~~~~~~~~
\end{eqnarray*}
We have
\begin{theorem}\label{C-attained} Let $\alpha,\beta,\lambda,p,q$ satisfy \eqref{WH-exp-1} with $p<q$ excepted $\alpha=\beta=0$.  Then, $N_{\alpha,\beta,\lambda,p}$ can be attained. Namely, there exists some nonnegative function $f$ such that $\|t^\alpha f\|_{L^p(\mathbb{R}^{n+1}_+)}=1$ and $N_{\alpha,\beta,\lambda,p}=\|z^{-\beta}E_\lambda f\|_{L^q(\mathbb{R}^{n+1}_+)}.$
\end{theorem}
\begin{rem}\label{Extremal noexistence} $(i)$ When $p=q,$ the extremal functions of inequality \eqref{WHLSD-1} can not be expected to attain, see e.g., \cite{H1977} for the case $p=q=2, \, \alpha=0, \,  \beta=1/2, \lambda=n-1.$

$(ii)$ For $\alpha=\beta=0,$ Dou and Zhu \cite{DZ2019} proved that the best constant $N_{0,0,\lambda,p}$ did not attain (see Proposition 2.1 in \cite{DZ2019}).
\end{rem}
%\begin{rem}\label{Extr-rem}
 Lieb \cite{Lieb1983} firstly employed the Riesz rearrangement inequalities to prove the existence of extremal functions for inequality \eqref{class-HLS} and \eqref{class-WHLS} (with $\alpha,\beta\ge0$), respectively. Moreover, Dou and Zhu \cite{DZ2015a} also employed the Riesz rearrangement inequalities to prove the existence of  extremal functions  to HLS type inequality on the upper half space. However, the rearrangement inequalities are no longer valid to prove the existence of extremal functions for inequality \eqref{WHLSD-O}, since vertical weights $t$ and $z$ may
cause the nonexistence of minimizing sequence. Hence, we will show the existence of extremal functions  to inequality \eqref{WHLSD-O} by the concentration compactness principle, and drop out the restriction $\alpha,\beta\ge0.$ This technique is also adopted in \cite{CLT2019, HWY2008}.
 %We will modify some details which is different with the procedure in the paper \cite{CLT2019}.
%\end{rem}

The extremal functions to inequality \eqref{WHLSD-O}, up to a positive constant multiplier, satisfy the following Euler-Lagrange system
\begin{equation}\label{Euler-syst-2}
\begin{cases}
f^{p-1}(x,t)=\frac{1}{t^\alpha}\int_{\mathbb{R}^{n+1}_+} \frac{g(y,z)}{z^\beta|(x,t)-(y,z)|^{\lambda}}dydz, &\quad(x,t)\in\mathbb{R}^{n+1}_+,\\
g^{r-1}(y,z)=\frac{1}{z^\beta}\int_{\mathbb{R}^{n+1}_+} \frac{f(x,t)}{t^\alpha|(x,t)-(y,z)|^{\lambda}}dxdt, &\quad(y,z)\in\mathbb{R}^{n+1}_+.
\end{cases}
\end{equation} % with the conditions \eqref{WH-exp-1}.

Using the method of moving planes and the method of moving spheres, we show that the $f(x,t), \, g(x,t)$ are cylindrical symmetric and monotone decreasing with respect to $t$-axis (that is, $u$ and $v$ are radially symmetric and monotone decreasing with respect to some $x_0\in\mathbb{R}^n$) and the explicit form on the boundary of all positive solutions to the system \eqref{Euler-syst-2}.
\begin{theorem}\label{classfy-1} Let $\alpha,\beta,\lambda,p,r$ satisfy \eqref{WH-exp-0}
 and $p_\alpha=\frac{2(n+1)}{2(n+1)-\lambda-2\alpha}, r_\beta=\frac{2(n+1)}{2(n+1)-\lambda-2\beta}.$
Assume that $(f,g)$  is a pair of positive solutions of \eqref{Euler-syst-2} with $f\in L^{p+1}(\mathbb{R}^{n+1}_+),$  then the following results hold.

$(i)$ $f(x,t),g(x,t)$ are radially symmetric and monotone decreasing  with respect to $x$ about some $x_0\in\mathbb{R}^n$.

$(ii)$ For $p=p_\alpha, r=r_\beta,$  $f$ and $g$ have the following forms on the boundary
\[f(y,0)=c_1\big(\frac{1}{|y-y_0|^2+d^2}\big)^{\frac{2n+2-\lambda-2\alpha}{2}},\quad
g(y,0)=c_2\big(\frac{1}{|y-y_0|^2+d^2}\big)^{\frac{2n+2-\lambda-2\beta}{2}},
\]where some constants $c_1,c_2,d>0,$ and $y,y_0\in\mathbb{R}^{n}.$
\end{theorem}

In addition, for $p=p_\alpha, r=r_\beta,$  we need $f,g \in C^0(\overline{\mathbb{R}^{n+1}_+})$ in Theorem \ref{classfy-1} to carry out the method of moving spheres. Hence, under the assumption $f\in L^{p+1}(\mathbb{R}^{n+1}_+)$, we will prove $f,g \in C^0(\overline{\mathbb{R}^{n+1}_+})$. Moreover, we show the  following regularity results  by the regularity lift lemma and bootstrap technique. That is,
\begin{theorem}\label{Regularity-fg}
 Let $\alpha,\beta,\lambda,p,r$ satisfy \eqref{WH-exp-0}. If $(f,g)$  is a pair of positive solutions of \eqref{Euler-syst-2} with $f\in L^{p+1}(\mathbb{R}^{n+1}_+),$ then

$(i)$ $t^{\frac{\alpha}{p-1}}f,t^{\frac{\beta}{r-1}}g\in L^{\infty}(\mathbb{R}^{n+1}_+).$

$(ii)$  $t^{\frac{\alpha}{p-1}}f, t^{\frac{\beta}{r-1}}g\in C^\infty(\mathbb{R}^{n+1}_+\setminus\{t=0\})\cup C^{0,\gamma}_{loc}(\overline{\mathbb{R}^{n+1}_+})$ for any $0<\gamma<1$. In particular, $f, g\in C^\infty(\mathbb{R}^{n+1}_+\setminus\{t=0\})\cup C^{0,\gamma}_{loc}(\overline{\mathbb{R}^{n+1}_+})$.
\end{theorem}

On the same way, we obtain a more general Stein-Weiss type inequality  with  partial variable $\hat{x}\in\mathbb{R}^{m}$ weight function on $\mathbb{R}^{n+m}.$
\begin{theorem}\label{WHS-m}
Let $\alpha,\beta,\lambda,p,r$ satisfy
\begin{equation}\label{EWH-exp-0}
\begin{cases}
&n,m\ge0, 0<\lambda<n+m, \, 1< r, \, p<\infty ,\\
&\alpha<\frac{m}{p'}, \beta<\frac {m}{r'},~\alpha+\beta\ge0, \\
&\frac1p +\frac1r+\frac{\alpha+\beta+\lambda}{n+m}=2.
\end{cases}
\end{equation}
Then  there exists a positive constant $ C_{\alpha,\beta,\lambda,p}:=C(n,  m,\alpha,\beta,\lambda,p)$ depending on $n, m, \alpha,\beta,\lambda$ and $p$, such that for all $f(x,\hat{x})\in L^p(\mathbb{R}^{n+m}), \, g(y,\hat{y})\in L^r(\mathbb{R}^{n+m}),$
\begin{equation}\label{EWHLS}
\big|\int_{\mathbb{R}^{n+m}}\int_{\mathbb{R}^{n+m}}\frac{g(y,\hat{y}) f(x,\hat{x})}{|\hat{x}|^\alpha|(x,\hat{x})-(y,\hat{y})|^{\lambda}|\hat{y}|^\beta}
dxd\hat{x}dyd\hat{y}\big|
\le C_{\alpha,\beta,\lambda,p} \|f\|_{L^p(\mathbb{R}^{n+m})}\|g\|_{L^r(\mathbb{R}^{n+m})}
\end{equation}holds. Moreover, the upper and lower bounds of constant $C_{p,\alpha,\beta,\lambda}=C({n,m,p,\alpha,\beta,\lambda})$ satisfy
\begin{eqnarray*}
\max\big\{Q_1,Q_2,Q_3\big\}
&\le&C_{\alpha,\beta,\lambda, p}
\le\min\big\{\big(\frac{p}{p-1}\big)^{\frac{p-1}{p}}{p}^{\frac1p}Q_1,\big(\frac{p}{p-1}\big) ^{\frac{p-1}{p}}{p}^{\frac1p}Q_2\big\},
\end{eqnarray*}
where
\begin{eqnarray*}
%Q_0&=&{\color{blue}\big(\frac{\pi^{\frac {n+m-1}2}}{\Gamma(\frac{n+m}2)}\big)^{\frac{r-1}r}, \quad N^*=N(n+m,p, \lambda)},\\
%B^*&=&{\color{blue}B^*(n+m,p,\lambda)=w_{n+m-1}(2^{n+m+\lambda p'}-2^{-(n+m-\lambda p')}),}\\
Q_1&=&[\frac{(r-1)\cdot2\pi^{\frac {n+m-1}2}}{r(n+m-\beta-\lambda)-(n+m)}\cdot\frac{\Gamma\big(\frac{r(m-\beta )-m}{2(r-1)}\big)}{\Gamma\big(\frac{r(n+m-\beta)-(n+m)}{2(r-1 )}\big)\Gamma(\frac m2)}]^{\frac{r-1}{r}}\\
\quad&&\times[\frac{(p-1)\cdot2\pi^{\frac {n+m-1}2}}{p(n+m-\alpha)-(n+m)}\cdot\frac{\Gamma\big(\frac{p(m-\alpha)-m}{2(p-1)}\big)}{\Gamma\big(\frac{p(n+m-\alpha )-(n+m)}{2(p-1)}\big)\Gamma(\frac m2)}]^{\frac{p-1}{p}},\\
Q_2&=&[\frac{(r-1)\cdot2\pi^{\frac {n+m-1}2}}{r(n+m-\beta) -(n+m)}\cdot\frac{\Gamma\big(\frac{r(m-\beta) -m}{2(r-1)}\big)}{\Gamma\big(\frac{r(n+m-\beta)-(n+m)}{2(r-1)}\big)\Gamma(\frac m2)}]^{\frac{r-1}{r}}\\
\quad&&\times[\frac{(p-1)\cdot2\pi^{\frac {n+m-1}2}}{p(n+m-\alpha-\lambda)-(n+m)}\cdot\frac{\Gamma\big(\frac{p(m-\alpha)-m}{2(p-1)}\big)}{\Gamma\big(\frac{p(n+m-\alpha )-(n+m)}{2(p-1)}\big)\Gamma(\frac m2)}]^{\frac{p-1}{p}},\\
Q_3&=&\max\big\{\big[\frac{(2^{n+m}-1)w_{n+m-1}}{n+m}\big]^{\frac{r-1}r},
N(n+m,p,\lambda)\big[\frac{(2^{n+m}-1)w_{n+m-1}}{n+m}\big]^{\frac{r-1}r},\\
& &2^\lambda\big[(2^{n+m}-\frac1{2^{n+m}})\frac{w_{n+m-1}}{n+m}\big]^{\frac {p-1}{p}}\big[\frac{(2^{n+m}-1)w_{n+m-1}}{n+m}\big]^{\frac{r-1}r}\big\},
\end{eqnarray*}
and $N(n+m,p,\lambda)$ is the best constant of  inequality \eqref{class-HLS} in $\mathbb{R}^{n+m}$.
 %\[
%N(n+m,p,\lambda)\le{\color{blue}\frac{n+m}{(n+m-\lambda)}\big(\frac{|\mathbb{S}^{n+m-1}|}{n+m}\big)^{\frac \lambda {n+m}}\frac 1{pr}\big(\big(\frac{\lambda/(n+m)}{1-1/p}\big)^{\frac \lambda{n+m}}+\big(\frac{\lambda/(n+m)}{1-1/ r}\big)^{\frac \lambda {n+m}}\big)},
%\]
%{\color{blue}here $|\mathbb{S}^{n+m-1}|=\frac{2\pi^{\frac{n+m}2}}{\Gamma(\frac{n+m}2)}$ is the area of unit sphere $\mathbb{S}^{n+m-1}\subset\mathbb{R}^{n+m}.$}
\end{theorem}

Define the fractional
integral operator (or Riesz potential operator)
\[I_\lambda f(y,\hat{y})=\int_{\mathbb{R}^{n+m}}\frac{f(x,\hat{x})}{|(x,\hat{x})-(y,\hat{y})|^\lambda}dxd\hat{x},  ~~~~\forall (y,\hat{y})\in\mathbb{R}^{n+m}, ~\lambda \in (0,n+m)\]
for Lebesgue-measurable  function $f(x,\hat{x})$ on $ \mathbb{R}^{n+m}$.
Let $q=r'$, the assumptions of exponents  \eqref{EWH-exp-0}  can be written as
\begin{equation}\label{EWH-exp-1}
\begin{cases}
&n, m\ge0,0<\lambda<n+m, \,  1< p\le q<\infty,\\
& \alpha<\frac{m}{p'}, \, \beta<\frac mq,~\alpha+\beta\ge0,\\
&\frac1q=\frac1p-\frac{n+m-(\alpha+\beta+\lambda)}{n+m}.
\end{cases}
\end{equation}

The inequality \eqref{EWHLS} is equivalent to the following dual form on $\mathbb{R}^{n+m}.$
 %with $q=r'$.
\begin{theorem}\label{EWHLSD-theo} Assume that $\alpha,\beta,\lambda,p,q$ satisfy \eqref{EWH-exp-1}. Then there is a positive constant $C_{p,\alpha,\beta,\lambda}$  depending on $n, m, \alpha, \beta,\lambda$ and $p$, such that
\begin{equation}\label{EWHLSD-1}
\||\hat{y}|^{-\beta}I_\lambda f\|_{L^q( \mathbb{R}^{n+m})}\le C_{p,\alpha,\beta,\lambda}\||\hat{x}|^\alpha f\|_{L^p(\mathbb{R}^{n+m})}
\end{equation}holds for any $f\in L^p(\mathbb{R}^{n+m}).$
\end{theorem}
\begin{rem}\label{Extremal remark0} Similar to the results of Remark \ref{Extremal remark}, we have the following equivalent forms about inequality \eqref{EWHLS}.

$(i)$ The inequality \eqref{EWHLS} is also equivalent to
\begin{equation}\label{EWHLS-1}
\big|\int_{\mathbb{R}^{n+m}}\int_{\mathbb{R}^{n+m}}\frac{f(x,\hat{x}) g(y,\hat{y})}{|(x,\hat{x})-(y,\hat{y})|^{\lambda}}dxd\hat{x}dyd\hat{y}\big|
\le C_{p,\alpha,\beta,\lambda}\||\hat{x}|^\alpha f\|_{L^p(\mathbb{R}^{n+m})}\||{\hat{y}}|^\beta g\|_{L^r(\mathbb{R}^{n+m})}
\end{equation} with the same assumptions of  Theorem \ref{WHS-m}.

$(ii)$ Define a weighted singular integral operator
\[H_{\alpha,\beta}f(y,\hat{y})=\int_{\mathbb{R}^{n+m}}\frac{f(x,t)}{|\hat{x}|^\alpha |(x,\hat{x})-(y,\hat{y})|^{\lambda}|\hat{y}|^\beta} dxd\hat{x},
\]
and $\alpha,\beta,\lambda,p,q$ satisfy \eqref{EWH-exp-1}.
Then inequality \eqref{EWHLS}  is equivalent to the following form
\begin{equation}\label{EWHLS-2}
\| H_{\alpha,\beta} f\|_{L^q( \mathbb{R}^{n+m})}\le C_{p,\alpha,\beta,\lambda}\|f\|_{L^p(
\mathbb{R}^{n+m})}.
\end{equation}

$(iii)$ For the conformal case $p= \tilde{p}_\alpha:=\frac{2(n+m)}{2(n+m)-\lambda-2\alpha}, \, r=\tilde{r}_\beta:=\frac{2(n+m)}{2(n+m)-\lambda-2\beta}$, by  the stereographic
projection, the inequality \eqref{EWHLS} is equivalent to the following weighted inequality on unit sphere $\mathbb{S}^{n+m}$
\begin{eqnarray*}
\big|\int_{\mathbb{S}^{n+m}}\int_{\mathbb{S}^{n+m}}
\frac{F(\xi)G(\eta)}{|\xi^m|^\alpha|\xi-\eta|^\lambda|\eta^m|^\beta}d\xi d\eta\big|
\le C_{\tilde{p}_\alpha,\alpha,\beta,\lambda}\|F\|_{L^{\tilde{p}_\alpha}(\mathbb{S}^{n+m})}
\|G\|_{L^{\tilde{r}_\beta}(\mathbb{S}^{n+m})}
\end{eqnarray*} for $F\in{L^{\tilde{p}_\alpha}(\mathbb{S}^{n+m})}$ and $G\in{L^{\tilde{r}_\beta}(\mathbb{S}^{n+m})}$, where $\xi=(\xi^n,\xi^m,\xi_{n+m+1})=(\xi_1,\xi_2,\cdots,$ $\xi_n,\xi_{n+1},\cdots,\xi_{n+m},\xi_{n+m+1})\in\mathbb{S}^{n+m}, $ $ \eta=(\eta^n,\eta^m,\eta_{n+m+1})
=(\eta_1,\eta_2,\cdots,\eta_n,\eta_{n+1},$  $\cdots,\eta_{n+m},\eta_{n+m+1})\in\mathbb{S}^{n+m}$.

$(iv)$ If $n=0, m\ge1$, inequality \eqref{EWHLS} is just equivalent to the inequality \eqref{class-WHLS}. In addition, Chen, Lu and Tao \cite{CLT2019}  also established inequality \eqref{EWHLS} by  Sawyer and Wheeden condition on weighted inequalities for integral operators. They also showed the symmetry results the method of moving plane in integral forms.  In addition Beckner \cite{B2021} also established inequality \eqref{EWHLS}.
\end{rem}

Let $\mathbb{H}^n=\{X=(x,t)\in\mathbb{R}^{n+1,1}:\, <X,X>=|x|^2-t^2=-1, t>0\}$ be the hyperbolic space,  and  $d(w,w')=\frac{|w-w'|}{\sqrt{tz}}$ ($t,z>0$) denotes the Poincar\'{e} metric on  $\mathbb{H}^n$ with $w, w'\in \mathbb{H}^n$.  We  can extend the results  of Beckner in \cite{B2021} as follows.

\begin{theorem} For $F\in L^p(\mathbb{H}^n), G\in L^r(\mathbb{H}^n),\, 1<p,r<\infty, 0<\lambda<n+1, \alpha+\beta\ge0$ and $ \lambda=\frac{n+1}{p'}+\frac{n+1}{r'}-\alpha-\beta, $ the Stein-Weiss type inequality  \eqref{WHLSD-O} ($f,g\ge0$) is equivalent to the following inequality on the hyperbolic space
$$\int_{\mathbb{H}^n}\int_{\mathbb{H}^n}\frac{F(w)G(w')}{d(w,w')^{\lambda}}dwdw'\le C_{p,\alpha,\beta,\lambda }\|F\|_{L^p(\mathbb{H}^n)}\|G\|_{L^r(\mathbb{H}^n)}.$$

Moreover, inequality \eqref{WHLSD-O} is controlled by inequality \eqref{EWHLS} for $f,g\ge0$, $f(x,\hat{x})=|\hat{x}|^{-\alpha} F(w)$, and $ g(y,\hat{y}) =|\hat{y}|^{-\beta} G(w')$, where $w=(x,\hat{x})\in\mathbb{R}^{n+m}, w'=(y,\hat{y})\in\mathbb{R}^{n+m}.$
\end{theorem}

 This theorem can be immediately verified  by  choosing $F(w)=f(x,t)t^{\frac{n+1}p}, G(w')=g(y,z)z^{\frac{n+1}r}$, where $w=(x,t)\in\mathbb{R}^{n+1}_+, w'=(y,z)\in\mathbb{R}^{n+1}_+.$ So we only state the result and omit the proof process.

The best constant in \eqref{EWHLS} can be classified as
\begin{eqnarray*}
C_{p,\alpha,\beta,\lambda}
&=&\sup\big\{\int_{\mathbb{R}^{n+m}}\int_{\mathbb{R}^{n+m}}\frac{f(x,\hat{x}) g(y,\hat{y})}{|\hat{x}|^\alpha |(x,\hat{x})-(y,\hat{y})|^{\lambda}{|\hat{y}|^\beta}}dxd\hat{x}dyd\hat{y}\,:\,\\
& &\|f\|_{L^p(\mathbb{R}^{n+m})}=\|g\|_{L^r(\mathbb{R}^{n+m})}=1\big\} \\
&=&\sup\big\{\|{|\hat{y}|}^{-\beta}I_\lambda f\|_{L^q(\mathbb{R}^{n+m})}\,:\, f\in L^p{({\mathbb{R}^{n+m}})},\,\|{|\hat{x}|}^\alpha f\|_{L^p(\mathbb{R}^{n+m})}=1\big\}.
%&=&\sup\big\{\|S f\|_{L^q(\mathbb{R}^{n+m})}\,:\, \|f\|_{L^p(\mathbb{R}^{n+m})}=1\big\}.~~~~~~~~~~~~~~~~~~~~~~~~~~~~~~~~~~~~~~~~~~~~~~~~~~~~~~~~~~~~~~~~
\end{eqnarray*}

Similar to Theorem \ref{C-attained}, by the concentration compactness principle,  we can show the existence of extremal functions of \eqref{EWHLSD-1}. That is

\begin{theorem}\label{C-mn-attained} Let $\alpha,\beta,\lambda,p,q$ satisfy \eqref{EWH-exp-1} and $q>p.$ Then, $C_{p,\alpha,\beta,\lambda}$ can be attained.
\end{theorem}

Corresponding to Euler-Lagrange equation of inequality \eqref{EWHLS}, up to a positive constant multiplier, we can write it as
\begin{equation}\label{Euler-syst-3}
\begin{cases}
f^{p-1}(x,\hat{x})=\frac{1}{|\hat{x}|^\alpha}\int_{\mathbb{R}^{n+m} } \frac{g(y,\hat{y})}{|\hat{y}|^\beta|(x,\hat{x})-(y,\hat{y})|^{\lambda}}dyd\hat{y},&\quad (x,\hat{x})\in\mathbb{R}^{n+m},\\
g^{r-1}(y,\hat{y})=\frac{1}{|\hat{y}|^\beta}\int_{\mathbb{R}^{n+m}} \frac{f(x,\hat{x})}{|\hat{x}|^\alpha|(x,\hat{x})-(y,\hat{y})|^{\lambda}}dxd\hat{x},&\quad (y,\hat{y})\in\mathbb{R}^{n+m}
\end{cases}
\end{equation}with the assumptions of {\eqref{EWH-exp-0}.
Arguing as in Theorem \ref{classfy-1}, we have
\begin{theorem}\label{classfy-2}
Let $\alpha,\beta,\lambda,p,r$ satisfy \eqref{EWH-exp-0} and $\tilde{p}_\alpha=\frac{2(n+m)}{2(n+m)-\lambda-2\alpha}$ and $ \tilde{r}_\beta=\frac{2(n+m)}{2(n+m)-\lambda-2\beta}$. Assume that $(f, g)$ is  a pair of positive solutions to system \eqref{Euler-syst-3} with $f\in L^{p+1}({\mathbb{R}^{n+m}})$, then the following results hold.

 $(i)$  $f(x,\hat{x}),g(x,\hat{x})$ are radially symmetric and monotone decreasing with respect to $x$ about some $x_0$ in $\mathbb{R}^n$; If $\alpha,\beta\ge0,$ then $f(x,\hat{x}),g(x,\hat{x})$ are radially symmetric and monotone decreasing with respect to $\hat{x}$ about the origin in $\mathbb{R}^m$.

 $(ii)$ For $p=\tilde{p}_\alpha, r=\tilde{r}_\beta,$  $f$ and $g$ have the following forms  on $\mathbb{R}^n\times\{\hat{y}=0\}$
\[
f(y,0)=c_3\big(\frac{1}{|y-y_0|^2+a^2}\big)^{\frac{2(n+m)-\lambda-2\alpha}2},\quad
~g(y,0)=c_4\big(\frac{1}{|y-y_0|^2+a^2}\big)^{\frac{2(n+m)-\lambda-2\beta}2},
\]where some constants $c_3,c_4,a>0,$ and $y,y_0\in\mathbb{R}^{n}$.
\end{theorem}
As applications, we can deduce some weighted Sobolev inequalities on $\mathbb{R}^{n+1}_+$ and $\mathbb{R}^{n+m}$, respectively. We will discuss the extremal problems of these inequalities in the future work.

Define $W_\alpha^{1,p}(\mathbb{R}^{n+1}_+)$ be the completion of  $C^\infty_0(\mathbb{R}^{n+1}_+)$ with respect to the norm
\[
\|u\|_{W_\alpha^{1,p}}=\int_{\mathbb{R}^{n+1}_+}t^{\alpha}|\nabla u|^pdxdt,
\]
and  $W_\alpha^{\gamma,p}(\mathbb{R}^{n+1}_+)$ be the completion of  $C^\infty_0(\mathbb{R}^{n+1}_+)$ with respect to the norm
\[
\|u\|_{W_\alpha^{\gamma,p}}=\int_{\mathbb{R}^{n+1}_+}t^{\alpha}|(-\Delta)^\frac\gamma2 u|^pdxdt
\]
with $0<\gamma<n+1$.

Theorem \ref{WHLSD-theo} implies the following weighted Sobolev inequality with partial variable weight $t$  on $\mathbb{R}^{n+1}_+$.
\begin{theorem}\label{Weighted-Sobolev-1} Let $n\ge1,\, \,1<p<p^*=p^*(\alpha_1, \beta_1):=\frac{p(n+1+\beta_1)}{n+1+\alpha_1-p}$, and
\[
\frac{pn}{p*}-(n+1-p)<\alpha_1<p-1,-1<\beta_1\le\frac{\alpha_1(n+1)}{n+1-p}.\]
Then  there exists some positive $S(n+1,p,\alpha_1,\beta_1)$ such that the following weighted inequality holds
\begin{equation}\label{WS-1}
\big(\int_{\mathbb{R}^{n+1}_+}t^{\beta_1}|u|^{p^*}dxdt\big)^{\frac p{p^*}}\le S(n+1,p,\alpha_1,\beta_1)\int_{\mathbb{R}^{n+1}_+}t^{\alpha_1}|\nabla u|^pdxdt
\end{equation}
for any $u\in W_{\alpha_1}^{1,p}(\mathbb{R}^{n+1}_+).$

Moreover, assume that $m$ is a positive integer with $mp< n$, and $\frac{pn}{p^*_m}-(n+1-mp)<\alpha_m<p-1,  -1< \beta_m\le\frac{\alpha_m(n+1)}{n+1-mp},  1<p<p_m^*=p^*(\alpha_m,   \beta_m):=\frac{p(n+1+\beta_m)}{n+1+\alpha_m-m p},$
the following weighted higher order Sobolev inequalities hold:

$(i)$ If $m=odd$,  then for any $u\in W_{\alpha_m}^{m,p}(\mathbb{R}^{n+1}_+)$ satisfying
$$\lim_{t\to0}\frac{\partial u(x,t)}{\partial{t}}=\lim_{t\to0}\frac {\partial\Delta u(x,t)}{\partial{t}}=\cdots=\lim_{t\to0}\frac {\partial\Delta^{\frac{m-1}2 -1} u(x,t)}{\partial{t}}=0,(m\neq1)$$
 there exists some constant $C(n+1,p,\alpha_m,\beta_m,m)$ such that
\begin{equation}\label{WS-mn-odd}
\big(\int_{\mathbb{R}^{n+1}_+}t^{\beta_m}|u|^{p_m^*}dxdt\big)^{\frac p{p_m^*}}\le C(n+1,p,\alpha_m,\beta_m,m)\int_{\mathbb{R}^{n+1}_+}t^{\alpha_m}\big|\nabla\Delta^\frac{m-1}2 u\big|^pdxdt.
\end{equation}

$(ii)$ If $m=even$, then for any $u\in W_{\alpha_m}^{m,p}(\mathbb{R}^{n+1}_+) $ satisfying
$$\lim_{t\to0}\frac{\partial u(x,t)}{\partial{t}}=\lim_{t\to0}\frac {\partial\Delta u(x,t)}{\partial{t}}=\cdots=\lim_{t\to0}\frac{\partial \Delta^{\frac{m}2-1} u(x,t)}{\partial{t}}=0,$$   there exists some constant $C(n+1,p,\alpha_m,\beta_m,m)$ such that
\begin{equation}\label{WS-mn-even}
\big(\int_{\mathbb{R}^{n+1}_+}t^{\beta_m}|u|^{p_m^*}dxdt\big)^{\frac p{p_m^*}}\le C(n+1,p,\alpha_m,\beta_m,m)\int_{\mathbb{R}^{n+1}_+}t^{\alpha_m}\big|\Delta^\frac{m}2 u\big|^pdxdt.
\end{equation}
\end{theorem}

According to Green function of  fractional Laplacian $(-\Delta)^\frac\gamma2$ with $0<\gamma<n+1$ on the upper half space $\mathbb{R}^{n+1}_+$ with zero boundary, we follow from \eqref{WHLSD-1} that the following weighted Sobolev inequality with partial variable $t$ weight  function holds for the fractional Laplacian $(-\Delta)^\frac\gamma2$.
\begin{theorem}\label{Fractional Sobolev-1} Assume that $0<\gamma<n+1, \, 1<p\le q<\infty, \, \alpha+\beta\ge0, \, \alpha<\frac {1}{p'}, \, \beta<\frac {1}q,$ and
\[\frac1q=\frac1p+\frac{\alpha+\beta-\gamma}{n+1}.\]
Then  there exists some positive constant $S(n+1,p,\alpha,\beta,\gamma)$ such that
\begin{equation*}
\|t^{-\beta}u\|_{L^q({\mathbb{R}^{n+1}_+})}\le S(n+1,p,\alpha,\beta, \gamma)\|t^\alpha(-\Delta)^{\frac\gamma2}u\|_{L^p(\mathbb{R}^{n+1}_+)}
\end{equation*} holds for $u\in W_{\alpha}^{\gamma,p}(\mathbb{R}^{n+1}_+)\cap L_\gamma$, where $ L_\gamma:=\big\{u\in L^1_{loc} \,|\, \int_{\mathbb{R}^{n+1}}\frac{|u(x,t)|}{1+|(x,t)|^{n+1+\gamma}}dxdt\big\}$.
\end{theorem}
Similarly, using Green function of fractional Laplacian $(-\Delta)^\frac\gamma2$ with $0<\gamma<n+m$ on whole space $\mathbb{R}^{n+m}$ and invoking \eqref{EWHLSD-1}, we also have the following weighted Sobolev inequality with partial variable t weight function   for fractional Laplacian $(-\Delta)^\frac\gamma2.$
\begin{theorem}\label{Fractional Sobolev-nm} Assume that $n,m\ge0, \, 0<\gamma<n+m, \, 1<p\le q<\infty, \, \alpha+\beta\ge0, \, \alpha<\frac {m}{p'}, \, \beta<\frac{m}q,$ and
\[
\frac1q=\frac1p+\frac{\alpha+\beta-\gamma}{n+m}.
\]
Then there exists some positive constant $S(n, m,p,\alpha,\beta,\gamma)$ such that
\begin{equation*}
\||\hat{x}|^{-\beta}u\|_{L^q({\mathbb{R}^{n+m}})}\le S(n,m,p,\alpha,\beta,\gamma)\||\hat{x}|^\alpha(-\Delta)^\frac\gamma2u\|_{L^p(\mathbb{R}^{n+m})}
\end{equation*}holds for any $u\in W_{\alpha}^{\gamma ,p}(\mathbb{R}^{n+m})$.

In particular, for $\frac1q=\frac1p+\frac{\alpha+\beta-1}{n+m},$ it holds
\begin{equation}\label{Sobolev-inq-mn}
 \||\hat{x}|^{-\beta}u\|_{L^{q}({\mathbb{R}^{n+m}})}\le S(n, m, p ,\alpha,\beta)\||\hat{x}|^\alpha\nabla u\|_{L^p(\mathbb{R}^{n+m})}
\end{equation}for $u\in W_{\alpha}^{1,p}(\mathbb{R}^{n+m})$.
%where $q=\frac{p(n+m+\beta)}{n+m+\alpha-p}, \, 1<p< p^*_1.$
\end{theorem}

\begin{rem}\label{Sobolev-thm-rem}
$(i)$ We point out that $\alpha$ and $p$ have the restriction condition $1+\frac nq-\frac{n+m}p<\alpha<\frac m{p'}=m-\frac mp$ and $p>1$  under the assumption of Theorem  \ref{Fractional Sobolev-nm}. From the restriction of $\alpha$, it is easy to see that $m>1$ if $p=q$. Particularly, if $m=1$, it must be $p<q$. This is also the reason why we assume $p<p^*$ in Theorem \ref{Weighted-Sobolev-1}. But Maz'ya \cite{Maz2011} established inequality \eqref{Sobolev-inq-mn} with $\alpha>1+\frac nq-\frac{n+m}p$ and $p\ge1$. Because we employed Stein-Weiss type inequality \eqref{EWHLSD-1}, the restriction $\alpha<\frac m{p'}$  and $p>1$ are necessary in here.

$(ii)$ For $ \alpha=\beta,p=2$,  Dong and Phan \cite{DP2023} also established a weighted Sobolev inequality \eqref{WS-1} on the upper half cylinder $(t_0-r^2,t_0)\times B^+_R(x_0)\subset \mathbb{R}^{n+1}_+$ by Stein-Weiss type inequality \eqref{class-WHLS} on $\mathbb{R}^1$. Indeed,  inequality \eqref{WHLSD-O} can provide an immediate proof.
\end{rem}
The paper is organized as follows. In section $2,$ we shall establish inequality \eqref{WHLSD-1} based on the weighted Hardy inequality and the layer cake representation technique. In Section $3,$ we show the existence of extremal functions via the concentration compactness principle. In Section $4,$ by the method of moving planes,  we show the cylindrical symmetry of nonnegative solutions to \eqref{Euler-syst-2}.  In section $5$, we show the form on the boundary of  nonnegative solutions to \eqref{Euler-syst-2} via the method of moving spheres. In section $6$, we give regularity results. Finally, as applications, we establish some weighted Sobolev inequalities in Section $7$.

In the whole paper, we always denote $c, C$ as a general positive constant whose value may be different from line to line.

%~~~~~~~~~~~~~~~~~~~~~~~~~~~~~~~~~~~~~~~~~~~~~~~~~~~~~~~~~~~~~~~~~~~~~~~~~~~~
\section{\textbf{The rough Stein-Weiss type inequality on the upper half space} \label{Section 2}}
In this section, we shall devote to establishing the rough Stein-Weiss type inequality  on the upper half space by weighted Hardy inequality and the layer cake representation technique.

For $R>0,$ denote
\begin{eqnarray*}
&&B_R(x,t)=\{(y,z)\in\mathbb{R}^{n}\times\mathbb{R}\,: \,|(y,z)-(x,t)|<R, (x,t)\in\mathbb{R}^{n+1}\},\\
&&B_R^+(0)=\{(y,z)=(y_1,y_2,\cdots,y_n,z)\in B_R(0)\cap \mathbb{R}^{n+1}_+\}.
\end{eqnarray*}
For simplicity, write
$B_R=B_R(0), \, B_R^+=B_R^+(0)$ and $X=(x,t), Y=(y,z)$ in the following.

We firstly recall the weighted Hardy inequality on whole space established by Dr\`{a}bek, Heinig, and Kufner in \cite{DHK1997} (also see \cite{Maz2011}).
\begin{lemma}\label{Weighted Hardy ineq}(Theorem $2.1$ and $2.2$ in \cite{DHK1997}) Let $W$ and $U$ be weighted functions on $\mathbb{R}^{n+1},$  $1<p\le q<\infty$ and $f\ge0$ on $\mathbb{R}^{n+1}$. The inequality
\begin{equation}\label{WH-1}
\big(\int_{\mathbb{R}^{n+1}}W(X)\big(\int_{B_{|X|}}f(Y)dY\big)^qdX\big)^\frac1q
\le C_0\big(\int_{\mathbb{R}^{n+1}}f^p(Y)U(Y)dY\big)^\frac1p
\end{equation}
holds if and only if
\begin{equation*}
A_0=\sup_{R>0}\big\{\big(\int_{|X|\ge R}W(X)dX\big)^\frac 1q\big(\int_{|Y|\le R}U^{1-p'}(Y)dY\big)^\frac1{p'}\big\}<\infty.
\end{equation*}
Moreover, if $C_0$ is the smallest constant of inequality \eqref{WH-1}, then
\[A_0\le C_0\le A_0(p')^\frac1{p'}p^\frac1{q}.\]
On the other hand, the inequality
\begin{equation}\label{WH-cond-4}
\big(\int_{\mathbb{R}^{n+1}}W(X)\big(\int_{\mathbb{R}^{n+1}\setminus B_{|X|}}f(Y)dY\big)^qdX\big)^\frac1q
\le C_1\big(\int_{\mathbb{R}^{n+1}}f^p(Y)U(Y)dY\big)^\frac1p
\end{equation}
holds if and only if
\begin{equation*}
A_1=\sup_{R>0}\big\{\big(\int_{|X|\leq R}W(X)dX\big)^\frac 1q\big(\int_{|Y|\geq R}U^{1-p'}(Y)dY\big)^\frac1{p'}\big\}<\infty.
\end{equation*}
Moreover, if $C_1$ is the smallest constant of inequality \eqref{WH-cond-4}, then
\[A_1\le C_1\le A_1 (p')^\frac1{p'}p^\frac1{q}.
\]
\end{lemma}

We can establish a weighted Hardy inequality on $\mathbb{R}^{n+1}_+$  employing the same method of Lemma \ref{Weighted Hardy ineq}.  That is,

\begin{lemma}\label{Weighted Hardy ineq-1} Let $W$ and $U$ be weighted functions on $\mathbb{R}_+^{n+1},$ $1<p\le q<\infty$ and $f\ge0$ on ${\mathbb{R}^{n+1}_+}.$ The following inequality
\begin{equation}\label{WH-2}
\big(\int_{\mathbb{R}^{n+1}_+}W(X)\big(\int_{B_{R}^{+}}f(Y)dY\big)^qdX\big)^\frac1q
\le C_2\big(\int_{\mathbb{R}^{n+1}_+}f^p(Y)U(Y)dY\big)^\frac1p
\end{equation}
holds if and only if
\begin{equation*}\label{WH-cond-2}
A_2=\sup_{R>0}\big\{\big(\int_{\mathbb{R}^{n+1}_+\setminus B^+_R}W(X)dX\big)^\frac1q\big(\int_{B^+_R}U^{1-p'}(Y)dY\big)^\frac1{p'}\big\}
<\infty.
\end{equation*}
Moreover, if $C_2$ is the smallest constant of inequality \eqref{WH-2}, then
\[ A_2\le C_2\le A_2 (p')^\frac1{p'}p^\frac1{q}.
\]
On the other hand, the inequality
\begin{equation}\label{WH-1-rev}
\big(\int_{\mathbb{R}^{n+1}_+}W(X)\big(\int_{\mathbb{R}^{n+1}_+\setminus B_{R}^{+}}f(Y)dY\big)^qdX\big)^\frac1q
\le C_3\big(\int_{\mathbb{R}^{n+1}_+}f^p(Y)U(Y)dY\big)^\frac1p
\end{equation}holds if and only if
\begin{equation*}\label{WH-cond-3}
A_3=\sup_{R>0}\big\{\big(\int_{B^+_R}W(X)dX\big)^\frac1q\big(\int_{\mathbb{R}^{n+1}_+\setminus B^+_R}U^{1-p'}(Y)dY\big)^\frac1{p'}\big\}<\infty.
\end{equation*}Moreover, if $C_3$ is the smallest constant of inequality \eqref{WH-1-rev}, then
\[A_3\le C_3\le A_3(p')^\frac1{p'}p^\frac1{q}.\]
\end{lemma}
\textbf{Proof of Theorem \ref{WHLSD-theo}.} Without loss of generality, assume that $f$ is a nonnegative function. Define
\begin{eqnarray*}
E_1f(Y) &=&\int_{B^+_{|Y|/2}}\frac{f(X)}{z^\beta|X-Y|^{\lambda}}dX,\\
E_2f(Y)&=&\int_{B^+_{2|Y|}\setminus B^+_{|Y|/2}}\frac{f(X)}{z^\beta|X-Y|^{\lambda}}dX,\\
E_3f(Y)&=&\int_{\mathbb{R}^{n+1}_+\setminus B^+_{2|Y|}}\frac{f(X)}{z^\beta|X-Y|^{\lambda}}dX.
\end{eqnarray*}
Thus, $E_\lambda f (Y)   z^{-\beta}$ can be divided into  the following three parts
\[E_\lambda f (Y)   z^{-\beta}=\int_{\mathbb{R}^{n+1}_+}\frac{f(X)}{z^\beta|X-Y|^{\lambda}}dX
=E_1f(Y)+E_2f(Y)+E_3f(Y).\]
We only need to show
\[
\int_{\mathbb{R}^{n+1}_+}|E_if|^qdY\le N^q_{\alpha,\beta,\lambda,p}\|t^\alpha f\|^q_{L^p(\mathbb{R}^{n+1}_+)},~~~~~~~~~~~~~~i=1,2,3.
\]

To this aim, we split it into three cases.

$\textbf{(1)}$ If $|X|\le\frac{|Y|}2$, then $|X-Y|\ge\frac{|Y|}2$. For any $\lambda>0,\, |X-Y|^{-\lambda}\le2^{\lambda}|Y|^{-\lambda}$, thus
\[
\int_{\mathbb{R}^{n+1}_+}|E_1f|^qdY\le\int_{\mathbb{R}^{n+1}_+}\big(\int_{B^+_{|Y|/2}}f(X)dX\big)^q z^{-\beta q}|Y|^{-\lambda q}dY.
\]
Selecting $W(Y)=z^{-\beta q}|Y|^{-\lambda q}$ and $U(X)=t^{\alpha p}$ in \eqref{WH-2}, we claim
\begin{equation}\label{E1-es-1}
\int_{\mathbb{R}^{n+1}_+}|E_1f|^qdY\le N^q_{\alpha,\beta,\lambda,p}\|t^\alpha f\|^q_{L^p(\mathbb{R}^{n+1}_+)}.
\end{equation} %if condition \eqref{WH-cond-1} is satisfied.
Indeed, from \eqref{WH-exp-1} we know
\[\frac {n+1}q={\alpha+\beta+\lambda}-\frac{n+1}{p'}.\]
Since $\alpha<\frac{1}{p'},$ we get ${n+1}-(\beta+\lambda)q=q(\alpha-\frac{n+1}{p'})<-\frac{q n}{p'}<0.$ For any $R>0,$ we have
\begin{eqnarray}\label{WH-6}
\int_{\mathbb{R}^{n+1}_+\setminus B^+_R}W(Y)dY&=&\int_{\mathbb{R}^{n+1}_+\setminus B^+_R}z^{-\beta q}|Y|^{-\lambda q}dY\nonumber\\
&=& J_{\beta q}\int^\infty_{R} \rho^{-(\beta+\lambda)q}\rho^{n}d\rho\nonumber\\
&=&C_1(n,\beta,\lambda,q)R^{n+1-(\beta+\lambda)q},
\end{eqnarray}
where $C_1(n,\beta,\lambda,q)=\frac{J_{\beta q}}{n+1-(\beta+\lambda)q}=\frac{\pi^{\frac n2}}{n+1-(\beta+\lambda)q}\cdot\frac{\Gamma(\frac{1-\beta q}2)}{\Gamma(\frac{n+1-\beta q}2)},$ and
\begin{eqnarray*}
J_{\beta q}&=&\int^\pi_0(sin\theta_1)^{n-1-\beta q}d\theta_1\int^{\pi}_{0}(sin\theta_2)^{n-2-\beta q}d\theta_2\cdots\int_{0}^{\pi}(sin\theta_n)^{-\beta q}d\theta_n\\
&=&\pi^{\frac n2}\frac{\Gamma(\frac{1-\beta q}2)}{\Gamma(\frac{n+1-\beta q}2)}.
\end{eqnarray*}
The fact $\beta<\frac{1}q$ was used on the above.
Moreover, since $\alpha<\frac{1}{p'},$ we have $\alpha p(1-p')+n+1=n+1-\alpha p'>0.$ For any $R>0,$ one has
\begin{eqnarray}\label{WH-7}
\int_{B^+_R}U^{1-p'}(X)dX&=&\int_{B^+_R} t^{\alpha p(1-p')}dX\nonumber\\
&=&J_{\alpha p'}\int_0^R \rho^{\alpha p(1-p')}\rho^nd\rho\nonumber\\
&=&C_2(n,\alpha,\lambda,p) R^{n+1-\alpha p'},
\end{eqnarray} where $C_2(n,\alpha,\lambda,p)=\frac{J_{\alpha p'}}{n+1-\alpha p'}=\frac{\pi^{\frac n2}}{n+1-\alpha p'}\cdot\frac{\Gamma(\frac{1-\alpha p'}2)}{\Gamma(\frac{n+1-\alpha p'}2)}$.
The fact $\alpha<\frac {1} {p'}$ was used on the above.
Combining \eqref{WH-6} and \eqref{WH-7} yields
\begin{eqnarray*}
&&\big(\int_{\mathbb{R}^{n+1}_+\setminus B^+_R}W(Y)dY\big)^\frac1q\big(\int_{B^+_R}U^{1-p'}(X)dX\big)^\frac1{p'}\\
&=&C_1^{\frac1q}(n,\beta,\lambda,q)C_2^{\frac1{p'}}(n,\alpha,\lambda,p)R^{\frac{n+1} q+\frac{n+1}{p'}-(\alpha+\beta+\lambda)}\\
&=&C_1^{\frac1q}(n,\beta,\lambda,q)C_2^{\frac1{p'}}(n,\alpha,\lambda,p),
\end{eqnarray*}
where $\frac{n+1} q+\frac{n+1}{p'}-(\alpha+\beta+\lambda)=0$ due to \eqref{WH-exp-1}. Hence,  the condition  of $A_2$ is satisfied in Lemma \ref{Weighted Hardy ineq-1}. We have proved \eqref{E1-es-1}.

$\textbf{(2)}$ If $|Y|\le \frac{|X|}2$, it implies $|X-Y|\ge\frac{|X|}2$. Similar to the case of $E_1$, taking $W(Y)=z^{-\beta q}$ and $U(X)=t^{ \alpha p}|X|^{ \lambda p}$ in \eqref{WH-1-rev}, we have
\begin{equation}\label{E3-es-1}
\int_{\mathbb{R}^{n+1}_+}|E_3f|^qdY
\le
N^q_{\alpha,\beta,\lambda,p}\|t^\alpha f\|^q_{L^p(\mathbb{R}^{n+1}_+)}.
\end{equation}
Since $\beta<\frac{1}{q},$  we know from \eqref{WH-exp-1} that $n+1-(\alpha+\lambda)p'<-\frac{np'}{q}<0.$ Thus, for any $R>0$, we have
\begin{eqnarray}\label{WH-8}
\int_{\mathbb{R}^{n+1}_+\setminus B^+_R}U^{1-p'}(X)dX&=&\int_{\mathbb{R}^{n+1}_+\backslash B^+_R}\big(t^{\alpha}|X|^\lambda\big)^{p(1-p')}dX\nonumber\\
&=&J_{\alpha p'}\int_{R}^\infty\rho^{(\alpha+\lambda)p(1-p')}\rho^nd\rho\nonumber\\
&=&C_3(n,\alpha,\lambda,p)R^{n+1-(\alpha+\lambda)p'},
\end{eqnarray}
where $C_3(n,\alpha,\lambda,p)=\frac{J_{\alpha p'}}{n+1-(\alpha+\lambda)p'}=\frac{\pi^{\frac n2}}{n+1-(\alpha+\lambda)p'}\cdot\frac{\Gamma(\frac{1-\alpha p'}2)}{\Gamma(\frac{n+1-\alpha p'}2)}$.
And
\begin{eqnarray}\label{WH-9}
\int_{B^+_R}W(Y)dY&=&\int_{B^+_R}z^{-\beta q}dY\nonumber\\
&=&J_{\beta q}\int_0^R\rho^{-\beta q}\rho^nd\rho\nonumber\\
&=&C_4(n,\beta,\lambda,q)R^{n+1-\beta q},
\end{eqnarray}
where $C_4(n,\beta,\lambda,q)=\frac{J_{\beta q}}{n+1-\beta q}=\frac{\pi^{\frac n2}}{n+1-\beta q}\cdot\frac{\Gamma(\frac{1-\beta q}2)}{\Gamma(\frac{n+1-\beta q}2)}$.
Combining \eqref{WH-8} and \eqref{WH-9}, we show that the condition  of $A_3$ holds in  Lemma \ref{Weighted Hardy ineq-1}, hence \eqref{E3-es-1} is proved.

$\textbf{(3)}$ If $\frac{|Y|}2<|X|< 2|Y|,$ it is not easy to estimate $E_2f$ by the idea of Stein and Weiss in \cite{SW1958}. We employ the layer cake representation technique to split the domain into the hollow cylinders with two cases $p<q$ and $p=q,$ respectively. This method is basic and is also employed to prove HLS inequality on nonhomogeneous spaces, see e.g., \cite{DHK1997, KM2005}.

\textbf{Case $1$.} $p<q.$  In the sequel, we use the notation
\begin{eqnarray*}
F_k&=&\big\{2^{k}\le|y|<2^{k+1},\, 2^{k}<z<2^{k+1}\big\},\\
G_k&=&\big\{2^{k-1}\le|x|<2^{k+2},\, 2^{k-1}<t<2^{k+2}\big\}
\end{eqnarray*}for any positive integer $k\in \mathbb{N}.$ Obviously, $\sum\limits_{k\in \mathbb{N}}F_k=\sum\limits_{k\in \mathbb{N}}G_k=\mathbb{R}^{n+1}_+$.
Since $|Y|<2|X|< 4|Y|$, if $Y\in F_k$, then we can  choose $X\in G_k$. And  write $\chi_K$ is the  standard characteristic function on domain $K$.

 Now, we divide $\lambda$ into two cases $\lambda>\frac{n+1}{p'}$ and $\lambda\le\frac{n+1}{p'}.$

If $\lambda>\frac{n+1}{p'},$ noting that  any $\alpha+\beta\ge0,$ it is easy to see
\[\frac{n+1}q=\frac{n+1}p-(n+1-(\alpha+\beta+\lambda))\ge\frac{n+1}p-(n+1-\lambda)=:\frac {n+1}{q_1}>0,\]
Obviously,  $q_1\ge q.$

By H\"{o}lder's inequality and the double form of  inequality \eqref{class-HLS}, we have
\begin{eqnarray}\label{WH-10}
& &\int_{\mathbb{R}^{n+1}_+} |E_2f|^qdY\nonumber\\
&=&\int_{\mathbb{R}^{n+1}_+}z^{-\beta q}\big(\int_{B^+_{2|Y|}\setminus B^+_{|Y|/2}}\frac{f(X)}{|X-Y|^{\lambda}}dX\big)^qdY\nonumber\\
&\le&\sum\limits_{k\in \mathbb{N}}\int_{F_k}z^{-\beta q}\big(\int_{B^+_{2|Y|}\setminus B^+_{|Y|/2}}\frac{f(X)}{|X-Y|^{\lambda}}dX\big)^qdY\nonumber\\
&\le&\sum\limits_{k\in \mathbb{N}}\big(\int_{F_k}(z^{-\beta q})^{\frac{q_1}{q_1-q}}dY\big)^{\frac{q_1-q}{q_1}}\big(\sum\limits_{k\in \mathbb{N}}\int_{F_k}\big(\int_{B^+_{2|Y|}\setminus B^+_{|Y|/2}}\frac{f(X)} {|X-Y|^{\lambda}}dX\big)^{q_1}dY\big)^{\frac q{q_1}}\nonumber\\
&\le&\sum\limits_{k\in\mathbb{N}}|F_k|^{\frac{q_1-q}{q_1}}2^{-k\beta q}\big(\int_{\mathbb{R}^{n+1}_+}\big(\sum\limits_{k\in \mathbb{N}} E_\lambda
(f\chi_{G_k})\big)^{q_1}dY\big)^{\frac q{q_1}}\nonumber\\
&\le&N^q(n+1,p,\lambda)\frac{w_{n-1}}n(2^n-1)\sum\limits_{k\in \mathbb{N}}2^{k\big(\frac{(n+1)(q_1-q)}{q_1}-\beta q\big)}\big(\sum\limits_{k\in \mathbb{N}}\int_{ G_k}|f(X)|^pdX\big)^{\frac qp}\nonumber\\
&\le&N^q(n+1,p,\lambda)\frac{w_{n-1}}n(2^n-1)\big(\int_{\mathbb{R}^{n+1}_+}\sum\limits_{k\in \mathbb{N}}2^{\frac {k p}q\big(\frac{(n+1)(q_1-q)}{q_1}-\beta q\big)}|f(X)|^pdX\big)^{\frac qp}\nonumber\\
&\le&N^q(n+1,p,\lambda)\frac{w_{n-1}}n(2^n-1)\big(\int_{\mathbb{R}^{n+1}_+}{(t^{\alpha q})}^{\frac pq}|f(X)|^pdX\big)^{\frac qp}\nonumber\\
%&=&\frac{w_n}n\cdot\sum\limits_{k\in N}2^{k\alpha q}(\sum\limits_{l\in N}\int_{F_l}|f(X)|^pdX)^{\frac qp}\\
%&\le&N_{\alpha,\beta,\lambda,p,q}(\int_{\mathbb{R}^{n+1}_+}|f(X)|^pz^{\alpha p}dX)^{\frac qp}\\
&\le& N^q_{\alpha,\beta,\lambda,p}\|t^\alpha f\|^q_{L^p(\mathbb{R}^{n+1}_+)},
\end{eqnarray}
where $|F_{k}|=\int^{2^{k+1}}_{2^k}dz\int^{2^{k+1}}_{2^k}dy=\frac{w_{n-1}}n\cdot2^{k(n+1)}(2^n-1)$  and $\frac{(n+1)(q_1-q)}{q_1}-\beta q=\alpha q.$

If $\lambda\le\frac{n+1}{p'},$  we have  $n+1-\lambda p'\ge0$.   By H\"{o}lder's inequality, we have
\begin{eqnarray}\label{WH-10-1}
& &\int_{\mathbb{R}^{n+1}_+} |E_2f|^qdY\nonumber\\
%&=&\int_{\mathbb{R}^{n+1}_+}z^{-\beta q}\big(\int_{B^+_{2|Y|}\setminus B^+_{|Y|/2}}\frac{f(X)}{|X-Y|^{\lambda}}dX\big)^qdY\nonumber\\
&\le&\int_{\mathbb{R}^{n+1}_+}z^{-\beta q}\big(\int_{B^+_{2|Y|}\setminus B^+_{|Y|/2}}|f(X)|^p dX\big)^{\frac qp}\big(\int_{B^+_{2|Y|}\setminus B^+_{|Y|/2}}|X-Y|^{-\lambda p'}dX\big)^{\frac q{p'}}dY\nonumber\\
&\le&[\frac{w_n}{n+1}2^{\lambda p'}(2^{n}-\frac1{2^{n+2}})]^{\frac q {p'}}\sum\limits_{k\in \mathbb{N}}\int_{F_k}z^{-\beta q}|Y|^{(n+1-\lambda p')\frac q{p'}}dY\nonumber\\
\quad &&\times\big(\sum\limits_{k\in \mathbb{N}}\int_{G_k}|f(X)|^pdX\big)^{\frac qp}\nonumber\\
%&\le&\sum\limits_{k\in \mathbb{N}}|F_k|2^{k(-\beta q+(n+1-\lambda p')\frac q{p'})}\big(\int_{\mathbb{R}^{n+1}_+}f^p(X)\chi_{G_k}dX\big)^{\frac qp}\nonumber\\
&\le&2^{\lambda q}(2^n-1)[\frac{w_n}{n+1}(2^{n}-\frac1{2^{n+2}})]^{\frac q {p'}}\frac {w_{n-1}}n\sum\limits_{k\in \mathbb{N}}2^{k(n+1-\beta q+(n+1-\lambda p')\frac q{p'})}\nonumber\\
\quad &&\times\big(\sum\limits_{k\in \mathbb{N}}\int_{G_k}|f(X)|^pdX\big)^{\frac qp}\nonumber\\
%&=&\frac{w_n}n\sum\limits_{k\in \mathbb{N}}\big(\int_{\mathbb{R}^{n+1}_+}2^{k(n+1-\beta q+(n+1-\lambda p')\frac q{p'})\frac pq}f^p(X)\chi_{G_k}dX\big)^{\frac qp}\nonumber\\
&\le&2^{\lambda q}(2^n-1)[\frac{w_n}{n+1}(2^{n}-\frac1{2^{n+2}})]^{\frac q {p'}}\frac {w_{n-1}}n\big(\int_{\mathbb{R}^{n+1}_+}{(t^{\alpha q})}^{\frac pq}|f(X)|^pdX\big)^{\frac qp}\nonumber\\
&\le& N^q_{\alpha,\beta,\lambda,p}\|t^\alpha f\|^q_{L^p(\mathbb{R}^{n+1}_+)},
\end{eqnarray}
where $n+1-\beta q+(n+1-\lambda p')\frac q{p'}=\alpha q,$ and
\begin{eqnarray*}
\int_{B^+_{2|Y|}\setminus B^+_{|Y|/2}}|X-Y|^{-\lambda p'}dX
&\le&\int_{B^+_{2|Y|}\setminus B^+_{|Y|/2}}||X|-|Y||^{-\lambda p'}dX\\
&\le&(\frac{|Y|}2)^{-\lambda p'}\int_{B^+_{2|Y|}\setminus B^+_{|Y|/2}}dX\\
%&=&2^{\lambda p'-1}{|Y|}^{-\lambda p'}[w_n{(2|Y|)}^{n+1}-w_n{(\frac{|Y|}2)}^{n+1}]\\
&=&\frac{w_n}{n+1}|Y|^{n+1-\lambda p'}(2^{n+\lambda p'}-2^{-(n+2-\lambda p')})\\
&=&\frac{w_n}{n+1}|Y|^{n+1-\lambda p'}2^{\lambda p'}(2^{n}-\frac1{2^{n+2}}).
%&=& B(n,p,\lambda)|Y|^{n+1-\lambda p'}.
\end{eqnarray*}

\textbf{Case $2 $.} $p=q.$ It follows from \eqref{WH-exp-1} that $n+1-(\alpha+\beta+\lambda)=0.$ By Young's inequality, we have
\begin{eqnarray}\label{WH-11}
& &\int_{\mathbb{R}^{n+1}_+}|E_2f|^qdY\nonumber\\
&=&\sum\limits_{k\in \mathbb{N}}\int_{F_k}z^{-\beta p}\big(\sum\limits_{k\in \mathbb{N}}E_{\lambda}(f\chi_{G_k})\big)^pdY\nonumber\\
&\le&\sum\limits_{k\in \mathbb{N}}2^{-k\beta p }\big(\int_{\mathbb{R}^{n+1}_+}|Y|^{-\lambda}\chi_{F_k}dY\big)^p \sum\limits_{k\in \mathbb{N}}\int_{\mathbb{R}^{n+1}_+}|f(X)|^p\chi_{G_k}dX\nonumber\\
&\le&(2^n-1)\frac{w_{n-1}}n\sum\limits_{k\in N}2^{pk(n+1-\beta-\lambda)}\sum\limits_{k\in N}\int_{G_k}|f(X)|^pdX\nonumber\\
&\le&(2^n-1)\frac{w_{n-1}}n\sum\limits_{k\in N}2^{pk(n+1-\alpha-\beta-\lambda)}\sum\limits_{k\in N}\int_{G_k}t^{\alpha p}|f(X)|^pdX\nonumber\\
&\le&N^q_{\alpha,\beta,\lambda,p}\sum\limits_{k\in \mathbb{N}}\int_{ G_k}t^{\alpha p}|f(X)|^pdX.
\end{eqnarray}
Thus, we have
\[\|E_2f\|_{L^q(\mathbb{R}^{n+1}_{+})}\le N_{\alpha,\beta,\lambda,p}\|t^\alpha f\|_{L^p(\mathbb{R}^{n+1}_+)}.\]
Finally, we will estimate the constant $N_{\alpha,\beta,\lambda,p}.$
%\[N_{\alpha,\beta,\lambda,p}=\min\{C_1(n,\beta,\lambda,q)C_2(n,\alpha,\lambda,p)\}.\]
For the sake of simplicity, write
\begin{eqnarray*}
D_1&=&C_1^{\frac1q}(n,\beta,\lambda,q)C_2^{\frac1{p'}}(n,\alpha,\lambda,p),\quad
D_2=C_3^{\frac1{p'}}(n,\alpha,\lambda,p)C_4^{\frac1q}(n,\beta,\lambda,p),\\
D_3&=&\max\{\big[(2^n-1)\frac{w_{n-1}}n\big]^{\frac1q},\big[(2^n-1)\frac{w_{n-1}}n\big]^{\frac1q}N(n+1,p,\lambda),\\
& &2^{\lambda }\big[(2^{n}-\frac1{2^{n+2}})\frac{w_n}{n+1}\big]^{\frac 1{p'}}\big[(2^n-1)\frac{w_{n-1}}n\big]^{\frac1q} \}
\end{eqnarray*}
From Lemma \ref{Weighted Hardy ineq-1}, \eqref{WH-10}, \eqref{WH-10-1} and \eqref{WH-11},  the upper and low bound of constant $N_{\alpha,\beta,\lambda,p}$ are obtained as follows
\begin{eqnarray*}
\max\big\{ D_1, D_2, D_3\}
&\le& N_{\alpha,\beta,\lambda,p}
\le\min\big\{(p')^{\frac1{p'}}p^{\frac1p}D_1, (p')^{\frac1{p'}}p^{\frac1p}D_2\big\},
\end{eqnarray*} and  it is easy see
\begin{eqnarray*}
%D_0&=&{\color{blue}\big(\frac{w_n}n\big)^{\frac1q}},\quad N=N(n,p,\lambda), \quad B=B(n,p,\lambda),\\
D_1&=&\big[\frac{\pi^{\frac n2}}{n+1-(\beta+\lambda)q}\cdot\frac{\Gamma(\frac{1-\beta q}2)}{\Gamma(\frac{n+1-\beta q}2)}\big]^{\frac1q}[\frac{\pi^{\frac n2}}{n+1-\alpha p'}\cdot\frac{\Gamma(\frac{1-\alpha p'}2)}{\Gamma(\frac{n+1-\alpha p'}2)}]^{\frac1{p'}},\\
D_2&=&\big[\frac{\pi^{\frac n2}}{n+1-(\alpha+\lambda)p'}\cdot\frac{\Gamma(\frac{1-\alpha p'}2)}{\Gamma(\frac{n+1-\alpha p'}2)}\big]^{\frac1{p'}}\big[\frac{\pi^{\frac n2}}{n+1-\beta q}\cdot\frac{\Gamma(\frac{1-\beta q}2)}{\Gamma(\frac{n+1-\beta q}2)}\big]^{\frac1q}.
\end{eqnarray*}
The proof is complete.
\qed

%~~~~~~~~~~~~~~~~~~~~~~~~~~~~~~~~~~~~~~~~~~~~~~~~~~~~~~~~~~~~~~~~~~~~~~~~~~~~

\section{\textbf{Existence of extremal functions} \label{Section 3}}
In this section, we present that the best constant of inequality \eqref{WHLSD-O} is attained  by concentration compactness principle introduced by P. Lions in \cite{Lions 1985}. Because that inequality \eqref{WHLSD-O} is equivalent to inequality \eqref{WHLSD-1} or \eqref{WHLSD-2}. For convenience, we will discuss that the best constant of inequality \eqref{WHLSD-1} can be attained.

The best constant of inequality \eqref{WHLSD-1} is classified by
\begin{equation}\label{cont}
N_{\alpha,\beta,\lambda,p}=\sup\big\{\|z^{-\beta} E_\lambda f \|_{L^q(\mathbb{R}^{n+1}_+)}\,:\, f\in
L^p{({ \mathbb{R}^{n+1}_+})},\,\|t^\alpha f\|_{L^p(\mathbb{R}^{n+1}_+)}=1\big\},~~~~~~~~~~~~~~~~~~~~~~~~~~~~~~~~~~~~~~~~~~~~~~~~~~~~~~~~~~~~~~~~
\end{equation}obviously, $N_{\alpha,\beta,\lambda,p}>0$.
%Let $f$ be a function defined on $\overline{\mathbb{R}^{n+1}_+}.$

For $\tau>0,$ define $f^{\tau}(X)=\tau^{-\frac{ n+1+\alpha p}p}f\big(\frac{X}{\tau}\big),$ it is easy to verify that
\[
\|t^\alpha f^{\tau}\|_{L^p{(\mathbb{R}^{n+1}_+)}}=\|t^\alpha f\|_{L^{p}{(\mathbb{R}^{n+1}_+)}}\quad\text{and} \quad\|z^{-\beta}E_\lambda f^{\tau}\|_{L^q{(\mathbb{R}^{n+1}_+)}}=\|z^{-\beta}E_\lambda f\|_{L^q{(\mathbb{R}^{n+1}_+)}}.
\]
Hence, the variational problem \eqref{cont} has both translation and dilation invariance with respect to $x\in\mathbb{R}^{n}$.

Define
\begin{eqnarray*}
L^p_\alpha(\mathbb{R}^{n+1}_+)&=&\big\{X\in \mathbb{R}^{n+1}_+\,\big|\, \big(\int_{\mathbb{R}^{n+1}_+}t^{\alpha p}|f|^pdX\big)^{\frac1p}<\infty\big\},\\
L^p_{\alpha,loc}(\mathbb{R}^{n+1}_+)&=&\big\{X\in K\,\big| \big(\int_{K}t^{\alpha p}|f|^pdX\big)^{\frac1p}<\infty,   \forall K\subset\mathbb{R}^{n+1}_+\big\}.
\end{eqnarray*}

Our main result states as
\begin{theorem}\label{max} Let $\alpha,\beta,\lambda,p,q$ satisfy \eqref{WH-exp-1} and $1<p<q<\infty.$ Suppose $\{f_i\}_{i=1}^\infty$ is a maximizing sequence of functions for \eqref{cont}, then after passing to a subsequence, there exists $\tau_i>0,$ such that $f^{\tau_i}_i\rightarrow f$ in $L_\alpha^p(\mathbb{R}^{n+1}_+).$ In particular, there exists at least one maximizing function for the variational problem \eqref{cont}.
\end{theorem}
To our aim, we need to establish a compact lemma and the concentration compactness principle.
\begin{lemma}\label{Compact-lm} Let $\alpha,\beta,\lambda,p,q$ satisfy \eqref{WH-exp-1}  and $1<p\le s<q<\infty$. Then, $$E_\lambda f: L^p_\alpha(\mathbb{R}^{n+1}_+)\hookrightarrow  L^s_{\beta, loc}(\mathbb{R}^{n+1}_+)$$
is compact.
\end{lemma}
\begin{proof} Since $\alpha+\beta\ge0,$ we have $\beta\ge0$ when $\alpha\le0,$ and  $\alpha\ge0$ when $\beta\le0.$ Using the dual form of inequality \eqref{WHLSD-O}, without loss of generality, we only need to consider the case $\beta\ge0.$ Let $\{f_i(X)\}\subset L^p_\alpha(\mathbb{R}^{n+1}_+)$ be a bounded sequence, then there exists a subsequence (still denoted by $\{f_i(X)\}$) and some function $f\in L^p_\alpha(\mathbb{R}^{n+1}_+)$ such that
\[f_i\rightharpoonup f ~\text{weakly in}~ L^p_\alpha(\mathbb{R}^{n+1}_+), ~\text{as}~ i\to \infty.\]

Without loss of generality, let any  $\Omega\subset B^+_R\subset\mathbb{R}^{n+1}_+$  for some $R>0$. Choosing some small $\varepsilon>0,$ we have
\begin{eqnarray}\label{bdes-1}
&&\big(\int_{\Omega}\big|E_\lambda(f_i-f)(Y)z^{-\beta}\big|^sdY\big)^{\frac1s}\nonumber\\
&=&\big(\int_{\Omega\cap\{z\ge \varepsilon\}}\big|E_\lambda(f_i-f)(Y)z^{-\beta}\big|^sdY\big)^{\frac1s}
+\big(\int_{\Omega\cap\{z< \varepsilon\}}\big|E_\lambda(f_i-f)(Y)z^{-\beta}\big|^sdY\big)^{\frac1s}.\quad\quad~~~
\end{eqnarray}
By H\"{o}lder's inequality and inequality \eqref{WHLSD-1},  we have
\begin{eqnarray}\label{bdes-2}
&&\big(\int_{\Omega\cap\{z< \varepsilon\}}\big|E_\lambda(f_i-f)(Y)z^{-\beta}\big|^sdY\big)^{\frac1s}\nonumber\\
&\le&C(R^{n}\varepsilon)^{\frac{(q-s)}{qs}}\big(\int_{\mathbb{R}^{n+1}_+}\big|E_\lambda (f_i-f)(Y)\chi_{\Omega\cap\{z< \varepsilon\}}|^qz^{-\beta q}dY\big)^{\frac1q}\nonumber\\
&\le&C\varepsilon^{\frac{(q-s)}{qs}}\|f_i-f\|^p_{L^p_\alpha(\mathbb{R}^{n+1}_+)}\to0 \quad\text{as}~ \varepsilon\rightarrow0.
\end{eqnarray}
Let $\varrho>0$ sufficiently small. We divide the first term of \eqref{bdes-1} into the following two parts
\begin{eqnarray}\label{bdes-3}
&&\big(\int_{\Omega\cap\{z\ge \varepsilon\}}z^{-\beta s}\big|\int_{
\mathbb{R}^{n+1}_+}\frac{f_i(X)-f(X)}{|X-Y|^\lambda}dX\big|^sdY\big)^{\frac1s}\nonumber\\
&=&\big(\int_{\Omega\cap\{z\ge\varepsilon\}}z^{-\beta s}\big|\int_{\{|X-Y|<\varrho\}\cap\mathbb{R}^{n+1}_+}\frac{f_i(X)-f(X)}{|X-Y|^\lambda}dX\big|^sdY\big)^{\frac1s}\nonumber\\
\quad&&+\big(\int_{\Omega\cap\{z\ge\varepsilon\}}z^{-\beta s}\big|\int_{\mathbb{R}^{n+1}_+\setminus \{|X-Y|\ge\varrho\}}\frac{f_i(X)-f(X)}{|X-Y|^\lambda}dX\big|^sdY\big)^{\frac1s}\nonumber\\
&=:&I_1+I_2.
\end{eqnarray}

We first estimate $I_1$.
Since $\alpha<\frac1{p'},$ by H\"{o}lder's inequality, it holds
\begin{eqnarray}\label{bdes-4}
\int_{\{|X-Y|<\varrho\}\cap\mathbb{R}^{n+1}_+}|f_i(X)-f(X)|dX
&\le&\big(\int_{\{|X-Y|<\varrho\}\cap\mathbb{R}^{n+1}_+}t^{\alpha p}| f_i(X)-f(X)|^pdX\big)^{\frac1p}\nonumber\\
\quad&&\times\big(\int_{\{|X-Y|<\varrho\}\cap\mathbb{R}^{n+1}_+}t^{-\alpha p'}dX\big)^{\frac1{p'}}<\infty.
\end{eqnarray}
For the sake of simplicity,  write
\begin{eqnarray*}
 A_z^1&=&\{|Y|\le|X-Y|<\varrho\}\cap\mathbb{R}^{n+1}_+, \\ A_z^2&=&\{z<|X-Y|<|Y|\}\cap\{|X-Y|<\varrho\} \cap\mathbb{R}^{n+1}_+, \\
 A_z^3&=&\{|X-Y|\le z\}\cap\{|X-Y|<\varrho\} \cap\mathbb{R}^{n+1}_+.
\end{eqnarray*}
Noting that $\beta, \lambda \ge0$ and $z\ge\varepsilon$,  we have
\begin{eqnarray}\label{bdes-5}
&&\int_{\{|X-Y|<\varrho\}}\frac1{z^{\beta}|X-Y|^\lambda}dY
\nonumber\\
&\le&\int_{A_z^1}\frac1{z^{\beta}|Y|^\lambda}dY
+\int_{A_z^2}\frac1{z^{\beta}|X-Y|^\lambda}dY
+\int_{A_z^3}\frac1{|X-Y|^{\lambda+\beta}}dY\nonumber\\
&\le&\int_{A_z^1}\frac1{z^{\beta}|Y|^\lambda}dY
+\int_{\varepsilon\le |X-Y|<\varrho}\frac1{z^{\beta}|X-Y|^\lambda}dY
+\int_{A_z^3}\frac1{|X-Y|^{\lambda+\beta}}dY\nonumber\\
\nonumber\\
&\le&J_{\beta }\int_0^\varrho\rho^{n-\beta-\lambda}d\rho+ \varepsilon^{-\lambda}J_{\beta}\int_0^\varrho\rho^{n-\beta} d\rho +c\int_0^\varrho\rho^{n-\beta-\lambda}d\rho\nonumber\\
&\le& C\varrho^{n+1-\lambda-\beta}\big(1+ \big(\frac{\varrho}{\varepsilon}\big)^{\lambda}\big)\nonumber\\
&\le& C\varrho^{n+1-\lambda-\beta}.\quad\quad(\text{we~ can ~choose}~\varrho=100\varepsilon)
\end{eqnarray}
Since $\alpha+\beta+\lambda< n+1,$ by Young's inequality, it follows from \eqref{bdes-5} and \eqref{bdes-4} that
\begin{eqnarray}\label{obf3}
& &\int_{\Omega\cap\{z\ge\varepsilon\}}\big|\int_{\{|X-Y|<\varrho\}\cap \mathbb{R}^{n+1}_+}\frac{f_i(X)-f(X)}{z^{\beta}|X-Y|^\lambda}dX\big|dY\nonumber\\
%&=&\int_{\mathbb{R}^{n+1}_+}\big|\int_{\{|X-Y|<\varrho\}\cap \mathbb{R}^{n+1}_+}\frac{f_i(X)-f(X)\chi_{\Omega\cap\{z\ge\varepsilon\}}}{z^{\beta}|X-Y|^\lambda}dX\big|dY\nonumber\\
&=&\int_{\mathbb{R}^{n+1}_+}\big|\int_{\{|X-Y|<\varrho\}\cap \mathbb{R}^{n+1}_+}\frac{f_i(X)-f(X)\chi_{\Omega\cap\{z\ge\varepsilon\}}}{z^{\beta}|X-Y|^\lambda}dY\big|dX\nonumber\\
&\le&\|f_i-f\|_{L^1(B^+_R)}\cdot\big\|\int_{\{|X-Y|<\varrho\}\cap \mathbb{R}^{n+1}_+}\frac1{z^{\beta}|X-Y|^\lambda}dY\big\|_{L^\infty(\Omega)}\nonumber\\
&\le&C\varrho^{n+1-\lambda-\beta}\|f_i-f\|^p_{L^p_\alpha(\mathbb{R}^{n+1}_+)}\rightarrow0,\quad\quad\text{as}~\varrho\rightarrow0.
\end{eqnarray}
Combining  \eqref{obf3} yields
\begin{equation}\label{EX-c-1}
I_1=\big(\int_{\Omega\cap\{z\ge\varepsilon\}}z^{-\beta s}\big|\int_{\{|X-Y|<\varrho\}\cap\mathbb{R}^{n+1}_+}\frac{f_i(X)-f(X)}
{|X-Y|^\lambda}dX\big|^sdY\big)^{\frac1s}\to0,\quad\text{as}~ i\rightarrow\infty.
\end{equation}

Next we estimate $I_2.$ Since $\beta<\frac1q,$ we know $n+1-(\alpha+\lambda)p'<0.$ One has
\begin{eqnarray*}
\int_{\{|X-Y|\ge\varrho\}\cap\mathbb{R}^{n+1}_+}\big(\frac1{t^\alpha|X-Y|^\lambda}\big)^{p'}dX
&\le&\int_{\{|X|\ge\varrho\}\cap\mathbb{R}^{n+1}_+}\big(\frac1{(t-z)^\alpha|X|^\lambda}\big)^{p'}dX\\
&=&J_{\alpha p'}\int_\varrho^\infty(\frac1{\rho^{\alpha+\lambda}})^{p'}\rho^nd\rho\\
&\le& C \varrho^{n+1-(\alpha+\lambda){p'}}
<\infty.
\end{eqnarray*}
Since $f_i\rightharpoonup f  $ weakly,  we have
\begin{eqnarray}\label{E-x2-c-1}
\int_{\{|X-Y|\ge\varrho\}\cap\mathbb{R}^{n+1}_+}\frac{f_i(X)-f(X)}{|X-Y|^\lambda}dX\to0,\quad\text{as}~ i\rightarrow\infty.
\end{eqnarray}
Therefore, by dominated convergence theorem,
\begin{equation}\label{bdes-8}
I_2=\big(\int_{\Omega\cap\{z\ge\varepsilon\}}z^{-\beta s}\big|\int_{\{|X-Y|\ge\varrho\}\cap\mathbb{R}^{n+1}_+}
\frac{f_i(X)-f(X)}{|X-Y|^\lambda}dX\big|^sdY\big)^{\frac1s}\to0, \quad\text{as}~ i\rightarrow\infty.
\end{equation}
Choosing $\varepsilon,\varrho$ sufficiently small and substituting \eqref{EX-c-1} and \eqref{bdes-8} into \eqref{bdes-3}, we arrive at
\begin{equation}\label{bdes-9}
\big(\int_{\Omega\cap\{z\ge \varepsilon\}}z^{-\beta s}\big|\int_{
\mathbb{R}^{n+1}_+}\frac{f_i(X)-f(X)}{|X-Y|^\lambda}dX\big|^sdY\big)^{\frac1s}\to0,\quad\text{as}~ i\rightarrow\infty.
\end{equation}
Furthermore, combining \eqref{bdes-2} and \eqref{bdes-9} into \eqref{bdes-1}, we conclude that the embedding $E_\lambda f: L^p_{\alpha}(\mathbb{R}^{n+1}_+)\hookrightarrow L^s_{\beta, loc}(\mathbb{R}^{n+1}_+)$ is compact. The proof is finished.
\end{proof}

\begin{rem}\label{rem converge} For any bounded sequence $\{f_i\}\subset L_\alpha^p(\mathbb{R}^{n+1}_+),$ there exists a subsequence (still denoted by $\{f_i\}$) and some function $f\in L^p_\alpha(\mathbb{R}^{n+1}_+)$ such that
\begin{gather*}
f_i\rightharpoonup f\quad\text{weakly in}\quad L^p_\alpha(\mathbb{R}^{n+1}_+),\\
E_\lambda {(f_i)}\rightharpoonup E_\lambda f \quad\text{weakly in}\quad L^q_\beta(\mathbb{R}^{n+1}_+),\\
E_\lambda{(f_i)}\rightarrow E_\lambda f \quad\text{strongly in} \quad L^s_{\beta,loc}(\mathbb{R}^{n+1}_+)
\end{gather*}for all $s\in [p,q).$ Furthermore, $E_\lambda(f_i)\rightarrow E_\lambda f$ pointwisely a.e. in $\mathbb{R}^{n+1}_+.$
\end{rem}
\begin{lemma}\label{c-0-c} Assume that $\{f_i\}\subset L_\alpha^p(\mathbb{R}^{n+1}_+)$ is a bounded nonnegative sequence and there exists some function $f\in L^p_\alpha(\mathbb{R}^{n+1}_+)$ such that
$$f_i\rightharpoonup f \quad\text{weakly in}\quad L^p_\alpha(\mathbb{R}^{n+1}_+).$$
After passing to a subsequence, assume that
\[t^{\alpha p}|f_i|^p{dX}\rightharpoonup\mu ~\text{weakly~in}~ \mathcal{M}(\mathbb{R}^{n+1}_+) ~\text{and}~ z^{-\beta q}|E_\lambda{ (f_i)}|^q dY \rightharpoonup\nu ~\text{weakly~in} ~\mathcal{M}(\mathbb{R}^{n+1}_+),\]
where $\mathcal{M}(\mathbb{R}^{n+1}_+)$ denotes the space of all Radon measures on $\mathbb{R}^{n+1}_+.$ Then

$i)$ There{\color{blue}} exists some countable set $J;$ a family $\{P_j: j\in J\}$ of distinct points in ${\mathbb{R}^{n+1}_+},$ and a family $\{\nu_j: j\in J\}$ of nonnegative numbers, such that
\[
\nu=|E_\lambda f|^qz^{-\beta q}dY+\sum_{j\in J}\nu_j\delta_{P_j},
\]where $\delta_{P_j}$ is the Dirac-mass of mass 1 concentrated at $P_j\in\mathbb{R}^{n+1}_+;$

$ii)$ In addition,
\[
\mu\ge|f|^p t^{\alpha p}dX+\sum_{j\in J}\mu_j\delta_{P_j}
\]for some family $\{\mu_j: j\in J\},$ where $\mu_j$ satisfies
\[
\nu_j^{\frac1q}\le N_{\alpha,\beta,\lambda,p}\mu_j^{\frac1p} ~\text{for all}~ j\in J.
\]
%Moreover, if $f=0$ and $\nu^{\frac1q}\le N_{\alpha,\beta,\lambda,p}\mu^{\frac1p}$,  then $\mu$ and $\nu$ are concentrated at a single point.
\end{lemma}
\begin{proof} We firstly prove $i)$. Assume that $\{f_i\}\subset L_\alpha^p(\mathbb{R}^{n+1}_+)$ is a bounded nonnegative sequence.  From Remark \ref{rem converge}, we know  that
\begin{gather*}
E_\lambda{ (f_i)}\rightharpoonup E_\lambda f \quad \text{weakly in}\quad  L_\beta^q(\mathbb{R}^{n+1}_+),\\
E_\lambda{ (f_i)}\rightarrow E_\lambda f \quad \text{strongly in}\quad L_{\beta,loc}^s(\mathbb{R}^{n+1}_+),\\
E_\lambda { (f_i)}\rightarrow E_\lambda f \quad \text{pointwisely a.e. in}~ \mathbb{R}^{n+1}_+
\end{gather*}
for all $s\in[p,q)$. Using Br\'{e}zis-Lieb lemma, we have
\begin{eqnarray*}
0&=&\lim_{i\rightarrow\infty}\int_{\mathbb{R}^{n+1}_+}z^{-\beta q}\big(|E_\lambda(f_i)|^q-|E_\lambda(f_i-f)|^q-|E_\lambda f|^q\big)dY\\
&=&\int_{\mathbb{R}^{n+1}_+}d\nu-\int_{\mathbb{R}^{n+1}_+}z^{-\beta q}|E_\lambda f|^qdY-\lim_{i\rightarrow\infty}\int_{\mathbb{R}^{n+1}_+}z^{-\beta q}|E_\lambda(f_i-f)|^qdY.
\end{eqnarray*} Hence, it is sufficient to discuss the case $f\equiv0.$

Suppose that $\varphi(Y)\in C^\infty_0(\overline{\mathbb{R}^{n+1}_+}),$  similar to the idea of P. Lions (see \cite{Lions1985(1), Lions 1985}), we only need to show that there exists some positive constant $C,$ such that
\begin{equation}\label{ESV}
\big(\int_{\mathbb{R}^{n+1}_+}|\varphi|^qd\nu\big)^{\frac1q}\le C\big(\int_{\mathbb{R}^{n+1}_+}|\varphi|^pd\mu\big)^{\frac1p}.
%\quad \forall\varphi\in C^\infty_0(\overline{\mathbb{R}^{n+1}_+}).
\end{equation}
Indeed,
\begin{eqnarray}\label{claim}
&&\big(\int_{\mathbb{R}^{n+1}_+}z^{-\beta q}|\varphi E_\lambda {(f_i)}|^qdY\big)^{\frac1q}\nonumber\\
&\le&\big(\int_{\mathbb{R}^{n+1}_+}z^{-\beta q}|E_\lambda(\varphi f_i)|^qdY\big)^{\frac1q}+\big(\int_{\mathbb{R}^{n+1}_+}z^{-\beta q}|\varphi E_\lambda {(f_i)}-E_\lambda(\varphi f_i)|^qdY\big)^{\frac1q}\nonumber\\
&\le&C\big(\int_{\mathbb{R}^{n+1}_+}t^{\alpha p}|\varphi f_i|^pdY\big)^{\frac1p}+\big(\int_{\mathbb{R}^{n+1}_+}z^{-\beta q}|\varphi E_\lambda {(f_i)}-E_\lambda (\varphi f_i)|^qdY\big)^{\frac1q}.\nonumber\\
\end{eqnarray}
Note that
\begin{eqnarray*}
z^{-\beta}|\varphi E_\lambda{(f_i)}-E_\lambda(\varphi f_i)|
&=&\big|\int_{\mathbb{R}^{n+1}_+}\frac{\big(\varphi(Y)-\varphi(X)\big)f_i(X)}{z^{\beta}|X-Y|^{\lambda}}
dX\big|\\
&\le&\big|\int_{|X-Y|\le R}\frac{\big(\varphi(Y)-\varphi(X)\big)f_i(X)}{z^{\beta}|X-Y|^{\lambda}}dX\big|\\
\quad&&+\big|\int_{|X-Y|> R}\frac{\big(\varphi(Y)-\varphi(X)\big)f_i(X)}{z^{\beta}|X-Y|^{\lambda}}dX\big|\\
&=:&J_1+J_2.
\end{eqnarray*}
Write $R(X,Y):=(\varphi(Y)-\varphi(X))|X-Y|^{-\lambda}.$ It is easy to verify that $R(X,Y)\in L^{s}(\mathbb{R}^{n+1}_+)$ for $s\le+\infty$ if $\lambda\le1$ and $s<\frac{n+1}{\lambda-1}$ if $\lambda>1.$ Note that
\[
J_1\le C\big|\int_{\{|X-Y|\le R\}}\frac{f_i(X)}{z^{\beta}|X-Y|^{\lambda-1}}dX\big|.
\]
For $\lambda\le1,$ using dominated convergence theorem, it holds
\begin{equation}\label{ESJ}
\lim_{i\rightarrow\infty}\int_{\mathbb{R}^{n+1}_+}J_1^qdY=0.
\end{equation}
For $\lambda>1,$ by inequality \eqref{WHLSD-1}, we have
\[J_1\in L^{l}(\mathbb{R}^{n+1}_+)\] for $\frac1{l}<\frac1q.$
Furthermore, $L^p_\alpha(\mathbb{R}^{n+1}_+)\hookrightarrow L^s_{\beta,loc}(\mathbb{R}^{n+1}_+)$ is a compact embedding for $s\in [p,q),$  we can obtain \eqref{ESJ} again. Since $R(X,Y)$ is uniformly bounded for $|X-Y|>R,$ by dominated convergence theorem, we obtain
\begin{eqnarray*}
\lim_{i\to\infty}\int_{\mathbb{R}^{n+1}_+}J_2^qdY
&=&\lim_{i\to\infty}\int_{\mathbb{R}^{n+1}_+}\int_{\{|X-Y|> R\}}z^{-\beta q}|R(X,Y)f_i(X)|^qdXdY\\
&=&\int_{\mathbb{R}^{n+1}_+}\int_{\{|X-Y|> R\}}z^{-\beta q}|R(X,Y)f(X)|^qdXdY\\
&=&0.
\end{eqnarray*}
Hence,
\begin{equation}\label{ESN}
\lim_{i\rightarrow\infty}\big(\int_{\mathbb{R}^{n+1}_+}z^{-\beta q}|\varphi E_\lambda {(f_i)}-E_\lambda(\varphi f_i)|^qdY\big)^{\frac1q}=0.
\end{equation}
Letting $i\rightarrow\infty$ in \eqref{claim}, and we derive that \eqref{ESV} from \eqref{ESN}. Furthermore, by P. Loin's Lemma (see Lemma 1.2 in \cite{Lions 1985}),  there exists some countable set $J$, a family $\{P_j: j\in J\}$ of distinct points in $\mathbb{R}^{n+1}_+$ such that
\[\lim_{i\rightarrow\infty}|E_\lambda(f_i-f)|^qz^{-\beta q}dY=\sum_{j\in J} \nu_j\delta_{P_j},\]
and
\[\nu=|E_\lambda f|^q z^{-\beta q}dY+\sum_{j\in J}\nu_j\delta_{P_j},\]
where $\nu_j= \nu(\{P_j\}).$

Next, we are ready to show $ii).$ From the assumption
\[f_i\rightharpoonup f ~\text{weakly~in}~ L_\alpha^p(\mathbb{R}^{n+1}_+),\]
 we have $\mu\ge t^{\alpha p}|f|^pdX.$ We just need to show that for each fixed $j\in J,$
\[\nu_j^{\frac1q}=\nu(\{P_j\})^{\frac1q}\le N_{\alpha,\beta,\lambda,p}\mu(\{P_j\})^{\frac1p}=N_{\alpha,\beta,\lambda,p}\mu_j^{\frac1p}.\]
Define $\varphi_{\tau}(Y)=\varphi\big(\frac{Y}{\tau}\big),$ where $\varphi(Y)\in C^\infty_0(\overline{\mathbb{R}^{n+1}_+})$ satisfies $0\le\varphi(Y)\le1,$ $\varphi(0)=1,$ supp$\varphi\subset B_1^{+}(0)$ and $\tau>0.$ Then
\[\big(\int_{\mathbb{R}^{n+1}_+}z^{-\beta q}|E_\lambda(\varphi_{\tau}f_i)|^qdY\big)^{\frac1q}\le C\big(\int_{\mathbb{R}^{n+1}_+}t^{\alpha p}|\varphi_{\tau}f_i|^pdX\big)^{\frac1p}.\]
Combining the above, we can obtain
\begin{eqnarray*}
&&\big(\int_{\mathbb{R}^{n+1}_+}z^{-\beta q}|\varphi_{\tau}\cdot E_\lambda (f_i)|^qdY\big)^{\frac1q}\nonumber\\
&\le&\big(\int_{\mathbb{R}^{n+1}_+}z^{-\beta q}|E_\lambda(\varphi_{\tau}f_i)|^qdY\big)^{\frac1q}
+\big(\int_{\mathbb{R}^{n+1}_+}
z^{-\beta q}|\varphi_{\tau}\cdot E_\lambda (f_i)-E_\lambda(\varphi_{\tau}f_i)|^qdY\big)^{\frac1q}\nonumber\\
&\le&N_{\alpha,\beta,\lambda,p}\big(\int_{\mathbb{R}^{n+1}_+}t^{\alpha p}|\varphi_{\tau}f_i|^pdX\big)^{\frac1p}+I,
\end{eqnarray*}
where
\[I:=\big(\int_{\mathbb{R}^{n+1}_+}z^{-\beta q}|\varphi_{\tau}\cdot E_\lambda(f_i)-E_\lambda(\varphi_{\tau}f_i)|^qdY\big)^{\frac1q}.\]
Arguing as the second term of \eqref{claim},  we can show
\[
I\rightarrow\big(\int_{\mathbb{R}^{n+1}_+}z^{-\beta q}|\varphi_{\tau}\cdot E_\lambda f-E_\lambda (\varphi_{\tau}f)|^qdY\big)^{\frac1q}\quad \text{as}~i\rightarrow\infty.
\]
Hence, letting $i\rightarrow\infty,$ we obtain
\begin{eqnarray*}
&&\big(\int_{\mathbb{R}^{n+1}_+}|\varphi_{\tau}|^qd\nu\big)^{\frac1q}\\
&\le&N_{\alpha,\beta,\lambda,p}\big(\int_{\mathbb{R}^{n+1}_+}|\varphi_{\tau}|^pd\mu\big)^{\frac1p}
+\big(\int_{\mathbb{R}^{n+1}_+}z^{-\beta q}|\varphi_{\tau}\cdot E_\lambda f-E_\lambda (\varphi_{\tau}f)|^qdY\big)^{\frac1q}.
\end{eqnarray*}
Since
\[\int_{\mathbb{R}^{n+1}_+}z^{-\beta q}|\varphi_{\tau}\cdot E_\lambda f|^qdY\rightarrow 0,% \quad as
\]
and
\[\big(\int_{\mathbb{R}^{n+1}_+}z^{-\beta q}|E_\lambda (\varphi_{\tau}f)|^qdY\big)^{\frac1q}
\le N_{\alpha,\beta,\lambda,p}\big(\int_{\mathbb{R}^{n+1}_+}t^{\alpha p}|\varphi_{\tau}f|^pdX\big)^{\frac1p}\rightarrow0
\]as $\tau\rightarrow0^+,$ we arrive at
\[\nu_j^{1/q}\le N_{\alpha,\beta,\lambda,p}\mu_j^{1/p} ~\text{for all}~  j\in J.\]
The proof is  completed.
\end{proof}

{\bf Proof of Theorem \ref{max}.} Without loss of generality, we assume  $\beta\ge0$.  For $\rho>0,$ define
\[Q_i(\rho)=\sup_{\xi\in{\mathbb{R}^{n}}}\int_{B^+_\rho(\xi,0)}t^{\alpha p}|f_i|^pdX.\]
Since for every $i,$
\[\lim_{\rho\rightarrow0^+}Q_i(\rho)=0, \quad \lim_{\rho\rightarrow\infty}Q_i(\rho)=1.\]
By introducing translation $\xi_i$ and dilation $\rho_i$ invariance for the maximizing sequence (replace $f_i$ by $f^{\rho_i}_i$ and still denote it as $f_i$ ), we may assume
\begin{equation}\label{LQ}
\frac12=\int_{B^+_{1}(0)}t^{\alpha p}|f_i|^pdX=\sup_{\xi_i\in{\mathbb{R}^n}}\int_{B^+_1(\xi_i,0)}t^{\alpha p}|f_i|^pdX.
\end{equation}
After passing to a subsequence, we may find $f\in L_\alpha^p(\mathbb{R}^{n+1}_+)$ such that
\begin{gather*}
f_i\rightharpoonup f ~\text{in}~  L_\alpha^p(\mathbb{R}^{n+1}_+),\\
t^{\alpha p}|f_i|^pdX\rightharpoonup\mu ~\text{in}~ \mathcal{M}(\mathbb{R}^{n+1}_+),\\
z^{-\beta q}|E_\lambda(f_i)|^qdY\rightharpoonup\nu ~\text{in}~ \mathcal{M}(\mathbb{R}^{n+1}_+).
\end{gather*}
From \eqref{LQ}, we have $\int_{ B^+_{1}(\xi_i,0)}t^{\alpha p}|f_i|^pdX\le\frac12.$

Firstly, we show that $\mu(\mathbb{R}^{n+1}_+)=1.$ We prove it by contradiction. Let $\mu(\mathbb{R}^{n+1}_+)=\sigma$ for $\sigma\in(0,1)$ and any small $\varepsilon$,  then we can find $\rho_1>0,$ such that
\begin{equation}\label{ro-1}
\int_{B^+_{\rho_1}}t^{\alpha p}|f_i|^pdX>\sigma-\varepsilon.
\end{equation}
Also we can find $\rho_i>\rho_1 $ large enough, such that
\begin{equation}\label{ro-i}
\int_{B^+_{\rho_i}}t^{\alpha p}|f_i|^pdX\le\sigma<\sigma+\varepsilon.
\end{equation}
Combining \eqref{ro-1} and \eqref{ro-i}, we have
\[
\sigma-\varepsilon<\int_{B^+_{\rho_1}}t^{\alpha p}|f_i|^pdX\le\int_{B^+_{\rho_i}}t^{\alpha p}|f_i|^pdX<\sigma+\varepsilon.
\]
Write
\[g_i=f_i\chi_{B^+_{\rho_1}}\quad\text{and} \quad h_i=f_i\chi_{\mathbb{R}^{n+1}_+\setminus B^+_{\rho_i}}. \]
Since
\[\int_{B^+_{\rho_i}\setminus B^+_{\rho_1}}t^{\alpha p}|f_i|^pdX\le2\varepsilon,\]
one has
\[\big(\int_{\mathbb{R}^{n+1}_+}t^{p\alpha}(f_i-g_i-h_i)^pdX\big)^\frac1p\le (2\cdot\varepsilon)^{\frac1p}.\]

On the other hand, by inequality \eqref{WHLSD-1} we have
\begin{eqnarray*}
&&\int_{\mathbb{R}^{n+1}_+}|E_\lambda (f_i)|^q z^{-\beta q}dY\\
&=&\int_{\mathbb{R}^{n+1}_+\setminus\{B^+_{\rho_i}\setminus B^+_{\rho_1}\}}|E_\lambda(f_i)|^q z^{-\beta q}dY
+\int_{B^+_{\rho_i}\setminus B^+_{\rho_1}}|E_\lambda(f_i)|^q z^{-\beta q}dY\\
&\le&\int_{\mathbb{R}^{n+1}_+\setminus\{B^+_{\rho_i}\setminus B^+_{\rho_1}\}}|E_\lambda(g_i+h_i)|^q z^{-\beta q}dY+N^q_{\alpha,\beta,\lambda,p}(\int_{B^+_{\rho_i}\setminus B^+_{\rho_1}}t^{\alpha p}|f_i|^pdX)^{\frac qp}\\
&\le&\int_{\mathbb{R}^{n+1}_+}|E_\lambda (g_i+h_i)|^q z^{-\beta q}dY+o(\varepsilon).
\end{eqnarray*}
We claim that
\begin{eqnarray}\label{E-claim}
& &\int_{\mathbb{R}^{n+1}_+}|E_\lambda(g_i+h_i)|^q z^{-\beta q}dY\nonumber\\
&=&\int_{\mathbb{R}^{n+1}_+}|E_\lambda (g_i)|^q z^{-\beta q}dY
 +\int_{\mathbb{R}^{n+1}_+}|E_\lambda(h_i)|^q z^{-\beta q}dY+o(\varepsilon).
\end{eqnarray}
In fact, using the fact ${n+1}-(\alpha+\lambda)p'<0$  and H\"{o}lder's inequality, it holds
\begin{eqnarray*}
&&|E_\lambda(g_i+h_i)(Y)-E_\lambda(g_i)(Y)|\\
&=&\int_{\mathbb{R}^{n+1}_+}\frac{f_i\chi_{\mathbb{R}^{n+1}_+\setminus B^+_{\rho_i}}}{|X-Y|^\lambda}dX\\
&\le&\big(\int_{\mathbb{R}^{n+1}_+\setminus B^+_{\rho_i}}|f_i(X)|^pt^{\alpha p}dX\big)^{\frac1p}\big(\int_{\mathbb{R}^{n+1}_+\setminus B^+_{\rho_i}}\big(\frac1{t^\alpha |X-Y|^\lambda}\big)^{p'}dX\big)^{\frac1{p'}}\\
&\to&0
\end{eqnarray*}as $\rho_i\to\infty$.
From the above inequality, we have
\begin{equation}\label{gie}
\lim_{\rho_i\rightarrow\infty}E_\lambda(g_i)(Y)=E_\lambda(g_i+h_i)(Y).
\end{equation}
By Brezis-Lieb lemma  and \eqref{gie}, the claim \eqref{E-claim} holds.

Since
$$\int_{\mathbb{R}^{n+1}_+}t^{\alpha p}|g_i|^pdX\le\sigma+\varepsilon,\quad \int_{\mathbb{R}^{n+1}_+}t^{\alpha p}|h_i|^pdX\le1-\sigma+\varepsilon,$$
we know that
\begin{eqnarray*}
N^q_{\alpha,\beta,\lambda,p}+o(1)&=&\int_{\mathbb{R}^{n+1}_+}z^{-\beta q}|E_\lambda (f_i)|^qdY\\
&=&\int_{\mathbb{R}^{n+1}_+}z^{-\beta q}|E_\lambda (g_i)|^qdY+\int_{\mathbb{R}^{n+1}_+}z^{-\beta q}|E_\lambda (h_i)|^qdY+o(\varepsilon)\\
&\le&N^q_{\alpha,\beta,\lambda,p}\big[\big(\int_{\mathbb{R}^{n+1}_+}t^{\alpha p}|g_i|^pdX\big)^{\frac qp}+ \big(\int_{\mathbb{R}^{n+1}_+}t^{\alpha p}|h_i|^pdX\big)^{\frac qp}\big]+o(\varepsilon).
\end{eqnarray*}
Letting $\varepsilon\rightarrow{0^+},$ we obtain
\[1\le\sigma^{\frac qp}+(1-\sigma)^{\frac qp}.\]
It contradicts $\frac qp>1.$ Hence $\mu(\mathbb{R}^{n+1}_+)=1.$

Secondly, we show that $\nu(\mathbb{R}^{n+1}_+)=N^q_{\alpha,\beta,\lambda,p}.$ For any $\varepsilon>0$ small, we can find $\rho_2>0$ such that
\[\int_{B^+_{\rho_2}}t^{\alpha p}|f_i|^pdX>1-\varepsilon\]
 for $i$  large enough. And then,
\begin{equation}\label{r-b2}
\int_{\mathbb{R}^{n+1}_+\setminus B^+_{\rho_2}}t^{\alpha p}|f_i|^pdX\le\varepsilon.
\end{equation}
From \eqref{r-b2} and \eqref{WHLSD-1}, we know that
\begin{eqnarray*}
\big(\int_{\mathbb{R}^{n+1}_+}\big|\int_{\mathbb{R}^{n+1}_+}\frac{f_i\chi_{\mathbb{R}^{n+1}_+\setminus B^+_{\rho_2}}}{z^\beta|X-Y|^\lambda}dX\big|^qdY\big)^{\frac1q}
&\le&N_{\alpha,\beta,\lambda,p}\big(\int_{\mathbb{R}^{n+1}_+\setminus B^+_{\rho_2}}t^{\alpha p} |f_i(X)|^pdX\big)^{\frac1p}\\
&\le&N_{\alpha,\beta,\lambda,p}\cdot\varepsilon^{\frac1p}.
\end{eqnarray*}
Letting $i\rightarrow\infty$, it yields
\begin{eqnarray*}
&&\lim_{i\rightarrow\infty}\int_{\mathbb{R}^{n+1}_+}|E_\lambda (f_i)|^q  z^{-\beta q}dY\\
&\ge&\lim_{i\rightarrow\infty}\int_{\mathbb{R}^{n+1}_+\cap\{z\le M\}}|E_\lambda (f_i)|^q z^{-\beta q}dY\\
&\ge&N^q_{\alpha,\beta,\lambda,p}-\lim_{i\rightarrow\infty}\int_{\mathbb{R}^{n+1}_+\cap\{z\ge M\}}|E_\lambda(f_i)|^q z^{-\beta q}dY\\
&=&N^q_{\alpha,\beta,\lambda,p}-\lim_{i\rightarrow\infty}\int_{\mathbb{R}^{n+1}_+\cap\{z\ge M\}}|E_\lambda (f_i\chi_{B^+_{\rho_2}})|^q z^{-\beta q}dY\\
& &-\lim_{i\rightarrow\infty}\int_{\mathbb{R}^{n+1}_+\cap\{z\ge M\}}|E_\lambda (f_i\chi_{\mathbb{R}^{n+1}_+\setminus B^+_{\rho_2}})|^q z^{-\beta q}dY\\
&\ge&N^q_{\alpha,\beta,\lambda,p}-C M^{n+1-(\lambda+\beta)q}\cdot\|t^{-\alpha}\|^q_{L^{p'}(B^+_{\rho_2})}\cdot\|t^\alpha f_i\|^q_{L^p(B^+_{\rho_2})}-N^q_{\alpha,\beta,\lambda,p}\cdot\varepsilon^{\frac qp}.
\end{eqnarray*}
Noting that ${n+1}-(\beta+\lambda)q<0,$ we obtain that $\nu(\mathbb{R}^{n+1}_+)=N^q_{\alpha,\beta,\lambda,p} $ as $M\rightarrow\infty$ and $\varepsilon\rightarrow0^+$ in the above.

Now, from Lemma \ref{c-0-c}, we have
\begin{eqnarray}\label{mume}
N^q_{\alpha,\beta,\lambda,p}&=&\|(E_\lambda f) z^{-\beta}\|^q_{L^q(\mathbb{R}^{n+1}_+)}+\sum_{j\in J}\nu_j\nonumber\\
&\le&N^q_{\alpha,\beta,\lambda,p}\|t^\alpha f\|^q_{L^p(\mathbb{R}^{n+1}_+)}+\sum_{j\in J}N^q_{\alpha,\beta,\lambda,p}(\mu_j)^{\frac q p}\nonumber\\
&=&N^q_{\alpha,\beta,\lambda,p}\big(\|t^\alpha f\|^q_{L^p(\mathbb{R}^{n+1}_+)}+\sum_{j\in J}(\mu_j)^{\frac q p}\big)\nonumber\\
&\le&N^q_{\alpha,\beta,\lambda,p}\big(\sigma^{\frac q p}+(1-\sigma)^{\frac q p}\big)\nonumber\\
&\le&N^q_{\alpha,\beta,\lambda,p},
\end{eqnarray} thus \eqref{mume} must be equality. Since $\frac pq<1,$ if $\|t^\alpha f\|^p_{L^p(\mathbb{R}^{n+1}_+)}=\sigma<1,$ then, $\mu=\delta_{y^*},$ $\nu=N^q_{\alpha,\beta,\lambda,p}\delta_{y^*}$ and $f=0.$ But we claim this is impossible to happen. We discuss this by two following cases.

{\bfseries Case 1.} If $y^*=(\xi^*,t^*)$ with $t^*=0,$ then $\int_{B^+_1(\xi^*,0)}d\mu=1,$ which contradicts  the initial assumption
\[\sup_{\xi^*\in{\mathbb{R}^n}}\int_{B^+_1(\xi^*,0)}t^{\alpha p}|f_i|^pdX=\frac12.\]

{\bfseries Case 2.} If $y^*=(\xi^*,t^*)$ with $t^*\neq0,$ we claim that there exists some $\delta>0$ such that
\[\lim_{i\rightarrow\infty}\int_{B_\delta(\xi^*,t^*)\cap\mathbb{R}^{n+1}_+}z^{-\beta q}|E_\lambda (f_i)|^qdY=0.\]
It contradicts $\int_{B_\delta(\xi^*,t^*)\cap\mathbb{R}^{n+1}_+}d\nu=N_{\alpha,\beta,\lambda,p}^q.$ In fact, since
\begin{equation*}
E_\lambda(f_i)(y^*)\rightarrow E_\lambda f(y^*)=0,\quad \forall y^*\in\mathbb{R}^{n+1}_+.
\end{equation*}
Let $l>1$ satisfy $\frac1l=\frac1p-\frac{n+1-(\alpha+\lambda)}{n+1}<\frac 1q,$ (note that $\beta\ge0$). By H\"{o}lder's inequality, we have
\begin{eqnarray}\label{rq}
&&\int_{B_\delta(\xi^*,t^*)\cap\mathbb{R}^{n+1}_+}|E_\lambda (f_i)|^q z^{-\beta q}dY\nonumber\\
&\le&\big(\int_{B_\delta(\xi^*,t^*)\cap\mathbb{R}^{n+1}_+}|E_\lambda (f_i)|^{l}dY\big)^{\frac ql}
\big(\int_{B_\delta(\xi^*,t^*)\cap\mathbb{R}^{n+1}_+}z^{-\frac{\beta ql}{l-q}}dY\big)^{\frac{l-q}l}\nonumber\\
&\le&C\cdot\delta^{\frac{(l-q)(n+1)}l-\beta q}\big(\int_{B_\delta(\xi^*,t^*)\cap\mathbb{R}^{n+1}_+}|E_\lambda (f_i)|^{l}dY\big)^{\frac ql}.
\end{eqnarray}
Using inequality \eqref{WHLSD-1}, one has
\begin{eqnarray}\label{Bd}
\int_{B_\delta(\xi^*,t^*)\cap\mathbb{R}^{n+1}_+}z^{-\beta q}|E_\lambda(f_i)|^qdY
&\le&C\cdot\delta^{\frac{(l-q)(n+1)}l-\beta q}\big(\int_{B_\delta(\xi^*,t^*)\cap\mathbb{R}^{n+1}_+}|E_\lambda (f_i)|^{l}dY\big)^{\frac qt}\nonumber\\
&\le&C\cdot\delta^{\frac{(l-q)(n+1)}l-\beta q}\big(\int_{\mathbb{R}^{n+1}_+}t^{\alpha p}|f_i|^pdX\big)^{\frac qp}.\nonumber\\
%&\le&C\cdot\delta^{\frac{(l-q)(n+1)}l-\beta q}.
\end{eqnarray}
Combining this  with $E_\lambda(f_i)(y^*)\rightarrow 0$ a.e. $y^*\in\mathbb{R}^{n+1}_+$, we  arrive at

\[\lim_{i\rightarrow\infty}\int_{B_\delta(\xi^*,t^*)\cap\mathbb{R}^{n+1}_+}z^{-\beta q}|E_\lambda (f_i)|^qdY=0.\]
We finish the proof of Theorem \ref{max}.
\qed

From Theorem \ref{max}, we know that the Euler-Lagrange equation of the variational problem \eqref{cont}, up to a positive constant multiplier, satisfies the following equation
\begin{equation*}\label{Euler-eq-1}
f^{p-1}(X)=\int_{\mathbb{R}^{n+1}_+} \frac{\big(Sf(Y)\big)^{q-1}}{t^\alpha |X-Y|^{\lambda}z^\beta }dY,~~~~~~~~~~\forall\,X\in\mathbb{R}^{n+1}_+,
\end{equation*}
which is equivalent to the system \eqref{Euler-syst-2} if $g=\big(Sf(Y)\big)^{q-1}$.

%~~~~~~~~~~~~~~~~~~~~~~~~~~~~~~~~~~~~~~~~~~~~~~~~~~~~~~~~~~~~~~~~~~~~~~~~~~~~

\section{\textbf{Cylindrical symmetry of the extremal functions} \label{Section 5}}
In this section, we will prove the cylindrical symmetry of  all nonnegative solutions to  the integral equation \eqref{Euler-syst-2} by the moving plane method in integral forms, which is introduced by Chen, Li and Ou \cite{CLO2006}.
%However, we will employ an improved process via Gluck \cite{G2020}.

For simplicity, write $u(X)=f^{p-1}(X),$ $v(Y)=Sf(Y)=g^{r-1}(Y)$ and let $\theta=\frac 1{p-1}, \,  \kappa=q-1.$ Then the system \eqref{Euler-syst-2} can be rewritten as
\begin{equation}\label{Euler-syst-1}
\begin{cases}
u(X)=\frac1{t^\alpha}\int_{\mathbb{R}^{n+1}_+} \frac{v^\kappa(Y)}{z^\beta|X-Y|^{\lambda}}dY, &\quad X\in\mathbb{R}^{n+1}_+,\\
v(Y)=\frac1{z^\beta}\int_{\mathbb{R}^{n+1}_+}\frac{u^ \theta(X)}{t^\alpha|X-Y|^{\lambda}}dX,&\quad Y\in\mathbb{R}^{n+1}_+
\end{cases}
\end{equation} with the conditions
\begin{equation}\label{WH-exp-2}
\begin{cases}
&0<\lambda<n+1,\, \kappa,\, \theta>0, \, \kappa\theta\ge1,\\
&\alpha<\frac1{\theta+1},\, \beta<\frac1{\kappa+1},\, \alpha+\beta\ge0, \\
&\frac1{\theta+1}+\frac1{\kappa+1}=\frac{\alpha+\beta+\lambda}{n+1}.
\end{cases}
\end{equation}
The result $(i)$ in Theorem \ref{classfy-1} is included in the following theorem.
\begin{theorem}\label{r-s-s} Let $\alpha,\beta,\lambda,\theta,\kappa$ satisfy \eqref{WH-exp-2}. If $(u,v)$ is a pair of positive solutions of \eqref{Euler-syst-1} with $u\in L^{\theta+1}(\mathbb{R}^{n+1}_+),$ then $u(X),v(X)$ are radially symmetric and monotone decreasing with respect to some $x_0\in\mathbb{R}^n$, that is $u(X)=u(|x-x_0|,t),\, v(X)=v(|x-x_0|,t).$
\end{theorem}For $\tau\in\mathbb{R}$, define
\[\Sigma_\tau=\{X=(x_1,x_2,\ldots x_n,t), x_1<\tau\},\]
and the  reflection  function
\[u_\tau(X)=u(X^\tau)\]
with $X^\tau=(2\tau-x_1,x_2,\ldots,x_n,t).$

For $X,Y\in\Sigma_\tau,$ if $\tau\le0$, then
\begin{equation}\label{notation a}
|X-Y^\tau|\ge|X-Y|,
\end{equation}and the inequality is strict when $\tau<0.$

We firstly show the representation formula with respect to system \eqref{Euler-syst-1} on $\Sigma_\tau.$
\begin{lemma}\label{notations} Let $(u,v)$ be a pair of solutions to \eqref{Euler-syst-1}. Then
\begin{eqnarray*}
&&u(X)-u_\tau(X)\nonumber\\
&=&\int_{\Sigma_\tau}\frac1{t^\alpha z^\beta}\big[\frac1{|X-Y|^\lambda}-\frac1{|X-Y^\tau|^\lambda}\big]
\big(v^\kappa(Y)-v^\kappa_\tau(Y)\big)dY,\nonumber\\
\end{eqnarray*}and
\begin{eqnarray}\label{v-v-0}
&&v(Y)-v_\tau(Y)\nonumber\\
&=&\int_{\Sigma_\tau}\frac1{t^\alpha z^\beta}\big[\frac{1}{|X-Y|^\lambda}-\frac{1}{|X-Y^\tau|^\lambda}\big]
\big(u^\theta(X)-u^\theta_\tau(X)\big)dX.\nonumber\\
\end{eqnarray}
\end{lemma}
\begin{proof} Using the change of variable and \eqref{Euler-syst-1}, we have
\begin{eqnarray}\label{rep-formula-1}
u(X)&=&\int_{\Sigma_\tau}\frac{v^\kappa(Y)}{t^\alpha |X-Y|^\lambda z^\beta }dY+\int_{\mathbb{R}^{n+1}_+\setminus\Sigma_\tau}\frac{v^\kappa(Y)}
{t^\alpha |X-Y|^\lambda z^\beta} dY\nonumber\\
&=&\int_{\Sigma_\tau}\frac{v^\kappa(Y)}{t^\alpha |X-Y|^\lambda z^\beta }dY+\int_{\Sigma{_\tau}}\frac{v^\kappa_\tau(Y)}
{t^\alpha|X-Y^\tau|^\lambda z^\beta}dY.\quad\quad
\end{eqnarray} Noting that $|X-Y^\tau|=|X^\tau-Y|$, it holds
\begin{eqnarray}\label{rep-formula-2}
u_\tau(X)&=&\int_{\Sigma_\tau}\frac {v^\kappa(Y)}{t^\alpha |X^\tau-Y|^\lambda z^\beta }dY+\int_{\Sigma{_\tau}}\frac {v^\kappa_\tau(Y)}{t^\alpha |X-Y|^\lambda z^\beta}dY\nonumber\\
&=&\int_{\Sigma_\tau}\frac {v^\kappa(Y)}{t^\alpha z^\beta |X-Y^\tau|^\lambda}dY+\int_{\Sigma_\tau}\frac{v^\kappa_\tau(Y)}{t^\alpha z^\beta |X-Y|^\lambda}dY.\quad\quad
\end{eqnarray}Combining \eqref{rep-formula-1} and \eqref{rep-formula-2} yields
\begin{eqnarray*}
&&u(X)-u_\tau(X)\\
&=&\int_{\Sigma_\tau}\frac {v^\kappa(Y)}{t^\alpha|X-Y|^\lambda z^\beta}dY+\int_{\Sigma{_\tau}}\frac{v^\kappa_\tau(Y)}
{t^\alpha|X-Y^\tau|^\lambda z^\beta}dY\\
\quad&&-\int_{\Sigma_\tau}\frac {v^\kappa(Y)}{t^\alpha|X- Y^\tau |^\lambda z^\beta}dY - \int_{\Sigma_\tau}\frac{v^\kappa_\tau(Y)}{t^\alpha |X-Y|^\lambda z^\beta}dY\\
&=&\int_{\Sigma_\tau}\frac1{t^\alpha z^\beta}\big[\frac{ 1}{|X-Y|^\lambda}-\frac{1}{|X-Y^\tau|^\lambda}\big]
(v^\kappa(Y)-v^\kappa_\tau(Y))dY.
\end{eqnarray*}Similarly, we can verify \eqref{v-v-0}. The proof is completed.
\end{proof}

Write
\begin{equation*}
\Sigma^u_\tau=\{X\in\Sigma_\tau|u(X)>u_\tau(X)\}~\text{and}~\Sigma^v_\tau=\{Y\in\Sigma_\tau|v(Y)>v_\tau(Y)\}.
\end{equation*}
The following lemma is the key to turn out the moving plane process.
\begin{lemma}\label{k-theta-u} Under the hypotheses of Theorem \ref{r-s-s}, there is some constant $C$ such that the following estimates hold.

$(a)$ If $\kappa\geq1$ and $\theta\geq1,$ then
\begin{equation}\label{u-v-4}
\|u-u_\tau\|_{L^{\theta+1}(\Sigma^u_\tau)}\le C\|u\|^{\theta-1}_{L^{\theta+1}(\Sigma^u_\tau)}
\|v\|^{\kappa-1}_{L^{\kappa+1}(\Sigma^v_\tau)}\|u-u_{\tau}\|_{L^{\theta+1}({\Sigma^u_\tau})},
\end{equation}and
\begin{equation}\label{u-v-5}
\|v-v_\tau\|_{L^{\kappa+1}(\Sigma^v_\tau)}\le C\|u\|^{\theta-1}_{L^{\theta+1}(\Sigma^u_\tau)}
\|v\|^{\kappa-1}_{L^{\kappa+1}(\Sigma^v_\tau)}\|v-v_{\tau}\|_{L^{\kappa+1}({\Sigma^v_\tau})}.
\end{equation}

$(b)$ If $\kappa\geq1$ and $\theta<1,$ then there exists $1<\frac1\theta\le r_1\le\kappa,$ such that
\begin{equation}\label{u-v-6}
\|v-v_\tau\|_{L^{\kappa+1}(\Sigma^v_\tau)}\le C\|u\|^{\theta-\frac1{r_1}}_{L^{\theta+1}(\Sigma^u_\tau)}\|v\|^{\frac\kappa{r_1}-1}_{L^{\kappa+1}(\Sigma^v_\tau)}
\|v-v_{\tau}\|_{L^{\kappa+1}(\Sigma^v_\tau)}.
\end{equation}

$(c)$ If $\kappa<1$ and $\theta>1,$ then there exists $1<\frac1\kappa\le r_2\le\theta,$ such that
\begin{equation}\label{u-v-7}
\|u-u_\tau\|_{L^{\theta+1}(\Sigma^u_\tau)}\le C\|u\|^{\frac\theta{r_2}-1}_{L^{\theta+1}(\Sigma^u_\tau)}\|v\|^{\kappa-\frac1{r_2}}_{L^{\kappa+1}(\Sigma^v_\tau)}
\|u-u_\tau\|_{L^{\theta+1}(\Sigma^u_\tau)}.
\end{equation}
\end{lemma}
 \begin{proof} The proof of this lemma is similar to that in \cite{G2020, HWY2008}, but we need more tedious computation.
 Let $\tau\le0$ for a.e. $X\in\Sigma_\tau.$ By Lemma \ref{notations}, we have
\begin{eqnarray*}
&&u(X)-u_\tau(X)\\
&=&\int_{\Sigma_\tau\setminus\Sigma^v_\tau}\frac{1}{t^\alpha z^\beta}\big[\frac 1{|X-Y|^\lambda}-\frac1{|X-Y^\tau|^\lambda}\big]\big(v^\kappa(Y)-v^\kappa_\tau(Y)\big)dY\\
\quad&&+\int_{\Sigma^v_\tau}\frac1{t^\alpha z^\beta}\big[\frac1{|X-Y|^\lambda}-\frac 1{|X-Y^\tau|^\lambda}\big]\big(v^\kappa(Y)-v^\kappa_\tau(Y)\big)dY.
\end{eqnarray*}
By the definition of $\Sigma^v_\tau$ and \eqref{notation a}, we obtain
\begin{eqnarray}\label{estimate-1}
&&u(X)-u_\tau(X)\nonumber\\
&\le&C\int_{\Sigma^v_\tau}\frac {1}{t^\alpha z^\beta}\big[\frac1{|X-Y|^\lambda}-\frac 1{|X-Y^\tau|^\lambda}\big]\big(v^\kappa(Y)-v^\kappa_\tau(Y)\big)dY\nonumber\\
&\le&C\int_{\Sigma^v_\tau}\frac {v^\kappa(Y)-v^\kappa_\tau(Y)}{t^{\alpha} |X-Y|^\lambda z^{\beta}}dY.
\end{eqnarray}Similarly,
\begin{eqnarray}\label{estimate-2}
&&v(Y)-v_\tau(Y)\nonumber\\
&\le&C\int_{\Sigma^u_\tau}\frac1{t^\alpha z^\beta}\big[\frac1{|X-Y|^\lambda}-\frac 1{|X-Y^\tau|^\lambda}\big]\big(u^\theta(X)-u^\theta_\tau(X)\big) dX \nonumber\\
&\le&C\int_{\Sigma^u_\tau}\frac {u^\theta(X)-u^{\theta}_\tau(X)}{t^{\alpha} |X-Y|^\lambda z^{\beta}}dX.
\end{eqnarray}

We split the remainder of the proof into three cases.

{\bf Case 1.} $\kappa\ge1$ and $\theta\ge1.$ For $X\in\Sigma^u_\tau,$ by the mean value theorem and \eqref{estimate-1}, we have
\[u(X)-u_\tau(X)\le \kappa\int_{\Sigma^v_\tau}\frac {v^{\kappa-1}(Y)\big(v(Y)-v_\tau(Y)\big)}{t^\alpha|X-Y|^\lambda z^\beta}dY.\]
Since $\frac1{\theta+1}+\frac1{\kappa+1}=\frac{\alpha+\beta+\lambda}{n+1}$, we know that $\frac{1}{\theta+1}=\frac\kappa{\kappa+1}-\frac{n+1-(\alpha+\beta+\lambda)}{n+1}$. Invoking the inequality \eqref{WHLSD-2}, we have
\[\|u-u_{\tau}\|_{L^{\theta+1}(\Sigma^u_\tau)}\le C\|v^{\kappa-1}(v-v_{\tau})\|_{L^{\frac{\kappa+1}\kappa}(\Sigma^v_\tau)}.\]
By H\"{o}lder's inequality, we obtain
\begin{equation}\label{u-v-3}
\|u-u_{\tau}\|_{L^{\theta+1}(\Sigma^u_\tau)}\le C\|v\|^{\kappa-1}_{L^{\kappa+1}(\Sigma^v_\tau)}
\|v-v_{\tau}\|_{L^{\kappa+1}({\Sigma^v_\tau})}.
\end{equation}% for some constant $C=C(n,\alpha,\beta,\lambda,\kappa,\theta)>0.$
Similarly, for any $Y\in\Sigma^v_\tau,$  the mean value theorem and  estimate  \eqref{estimate-2} show
\[v(Y)-v_\tau(Y)\le\theta\int_{\Sigma^u_\tau}\frac {u^{\theta-1}(X)\big(u(X)-u_\tau(X)\big)}{t^\alpha|X-Y|^\lambda z^\beta}dX.\]
Since $\frac1{\theta+1}+\frac1{\kappa+1}=\frac{\alpha+\beta+\lambda}{n+1}$, we  also have  $\frac1{\kappa+1}=\frac{\theta}{\theta+1}-\frac{n+1-(\alpha+\beta+\lambda)}{n+1}$.
Invoking the inequality \eqref{WHLSD-2} and H\"{o}lder's inequality, it yields
\begin{equation}\label{u-v-9}
\|v-v_\tau\|_{L^{\kappa+1}(\Sigma^v_\tau)}\le C
\|u\|^{\theta-1}_{L^{\theta+1}(\Sigma^u_\tau)}\|u-u_{\tau}\|_{L^{\theta+1}({\Sigma^u_\tau})}.
\end{equation}
%for some $C=C(n,\alpha,\beta,\lambda,\kappa,\theta)>0.$
Combining this into \eqref{u-v-3}, we arrive at \eqref{u-v-4}. And, substituting \eqref{u-v-3}
into \eqref{u-v-9}, we get \eqref{u-v-5}.

{\bf  Case 2.} $\kappa\ge1$ and $\theta<1.$ Let $r_1$ satisfy $1<\frac1\theta\le r_1\le\kappa.$  For $Y\in\Sigma^v_\tau,$ by the mean value theorem and \eqref{estimate-2}, we have
\[v(Y)-v_\tau(Y)\le\theta r_1\int_{\Sigma^u_\tau}\frac {u^{\theta-\frac1{r_1}}(X)\big(u^{\frac1{r_1}}(X)-u_\tau^{\frac1{r_1}}(X)\big)}{t^\alpha |X-Y|^\lambda z^\beta}dX.\]
Using \eqref{WHLSD-2} and H\"{o}lder's inequality, we have
\begin{eqnarray}\label{v-}
\|v-v_{\tau}\|_{L^{\kappa+1}(\Sigma^v_\tau)}&\le& C\|u^{\theta-\frac1{r_1}}(u^{\frac1{r_1}}-u^\frac1{r_1}_\tau)\|_{L^{\frac{\theta+1}\theta}(\Sigma^u_\tau)}\nonumber\\
&\le&C\|u\|^{\theta-\frac1{r_1}}_{L^{\theta+1}(\Sigma^u_\tau)}\|u^{\frac1{r_1}}-u^\frac1{r_1}_\tau\|_{L^{r_1(\theta+1)}({\Sigma^u_\tau})}.
\end{eqnarray}
 Define
\begin{eqnarray*}
I_1(X)&=&\int_{\Sigma^v_\tau}\frac{v^\kappa(Y)}{t^\alpha |X-Y|^\lambda z^\beta}dY,\quad
I_2(X)=\int_{\Sigma^v_\tau}\frac{v^\kappa_\tau(Y)}{t^\alpha |X-Y^\tau|^\lambda z^\beta}dY,\\
I_3(X)&=&\int_{\Sigma_\tau\setminus\Sigma^v_\tau}\frac{v^\kappa(Y)}{t^\alpha |X-Y|^\lambda z^\beta}dY,\quad
I_4(X)=\int_{\Sigma_\tau\setminus\Sigma^v_\tau}\frac{v^\kappa_\tau(Y)}{t^\alpha |X-Y^\tau|^\lambda z^\beta}dY.
\end{eqnarray*}We can write $u(X)=\Sigma^4_{j=1}I_j(X).$ Claim
\begin{equation}\label{I}
I_3(X^\tau)+I_4(X^\tau)\ge I_3(X)+I_4(X) ~ \text{for a.e.}~ X\in\Sigma_\tau.
\end{equation}To verify this claim, write
\[I_3(X)+I_4(X)-I_3(X^\tau)-I_4(X^\tau)=\int_{\Sigma_\tau\setminus\Sigma^v_\tau}Q(x,t,y,z)dY,\]
where the integrand $Q(X,Y)$ satisfies
\begin{eqnarray*}
&&t^\alpha z^\beta Q(X,Y)\\
&=&\frac{v^\kappa(Y)}{|X-Y|^\lambda}+\frac{v^\kappa_\tau(Y)}{|X-Y^\tau|^\lambda}
-\frac{v^\kappa(Y)}{|X-Y^\tau|^\lambda}-\frac{v^\kappa_\tau(Y)}{|X-Y|^\lambda}\\
&=&-\big(v^\kappa_\tau(Y)-v^\kappa(Y)\big)\big(\frac1{|X-Y|^\lambda}-\frac1{|X-Y^\tau|^\lambda}\big)\\
&\le&0
\end{eqnarray*}for a.e. $Y\in\Sigma_\tau\setminus\Sigma^v_\tau$ and $X\in\Sigma_\tau.$ It implies that \eqref{I} holds.

For $X\in\Sigma^u_\tau,$ define
\begin{eqnarray*}
&&a(X)=\int_{\Sigma^v_\tau}\frac{v^\kappa(Y)}{t^\alpha|X-Y|^\lambda z^\beta}dY+
\int_{\Sigma^v_\tau}\frac{v^\kappa(Y)}{t^\alpha|X-Y^\tau|^\lambda z^\beta}dY,\\
&&b(X)=\int_{\Sigma^v_\tau}\frac{v^\kappa_\tau(Y)}{t^\alpha |X-Y|^\lambda z^\beta}dY+\int_{\Sigma^v_\tau}\frac{ v^\kappa_\tau(Y)}{t^\alpha |X-Y^\tau|^\lambda z^\beta}dY,
\end{eqnarray*} we have both
\[
u(X)<a(X)+I_3(X)+I_4(X) \quad\text{and}\quad u_\tau(X)\ge b(X)+I_3(X^\tau)+I_4(X^\tau)
\]
for a.e. $X\in\Sigma^u_\tau.$ In addition, since $u_\tau(X)\ge b(X)$ in $\Sigma^u_\tau,$ there is a nonnegative function $c(X),$ such that
\[u(X)-a(X)\le c(X)\le u_\tau(X)-b(X)\quad\text{for a.e.}~X\in\Sigma^u_\tau.\]

For the sake of simplicity, write
\begin{eqnarray*}
\varphi_1&=&\big(\frac{v^\kappa(Y)}{t^\alpha |X-Y|^\lambda z^\beta}\big)^{\frac1{r_1}},\quad
\varphi_2=\big(\frac{v^\kappa(Y)}{t^\alpha |X-Y^\tau|^\lambda z^\beta}\big)^{\frac1{r_1}},\\
\psi_1&=&\big(\frac{v^\kappa_\tau(Y)}{t^\alpha |X-Y|^\lambda z^\beta}\big)^{\frac1{r_1}},\quad
\psi_2=\big(\frac{v^\kappa_\tau(Y)}{t^\alpha |X-Y^\tau|^\lambda z^\beta}\big)^{\frac1{r_1}}.
\end{eqnarray*}
We need the following basic inequality (see e.g., \cite{HWY2008}): Let $a\ge b\ge0,\ c\ge0,\ 0<m\le1$, it holds
\begin{equation}\label{nt-3}
(a+c)^m-(b+c)^m\le a^m-b^m.
\end{equation}For $r_1>1,$ by \eqref{nt-3}, it yields
\begin{eqnarray}\label{varphi-psi}
u^{\frac1{r_1}}(X)-u^{\frac1{r_1}}_\tau(X)&\le&\big(a(X)+c(X)\big)^{\frac1{r_1}}-\big(b(X)+c(X)\big)^{\frac1{r_1}}\nonumber\\
&\le&a^{\frac1{r_1}}(X)-b^{\frac1{r_1}}(X)\nonumber\\
&=&\big(\int_{\Sigma^v_\tau}\big|(\varphi_1,\varphi_2)\big|^{r_1}_{l^{r_1}}dY\big)^{\frac1{r_1}}-
\big(\int_{\Sigma^v_\tau}\big|(\psi_1,\psi_2)\big|^{r_1}_{l^{r_1}}dY)^{\frac1{r_1}}\nonumber\\
&\le&\big(\int_{\Sigma^v_\tau}\big|(\varphi_1-\psi_1,\varphi_2-\psi_2)\big|^{r_1}_{l^{r_1}}dY\big)^{\frac1{r_1}}\nonumber\\
&=&\big(\int_{\Sigma^v_\tau}|\varphi_1-\psi_1|^{r_1}+|\varphi_2-\psi_2|^{r_1}dY\big)^{\frac1{r_1}}
\end{eqnarray}for a.e. $X\in\Sigma^u_\tau.$ Moreover, for
$(X,Y)\in \Sigma^u_\tau\times\Sigma^v_\tau,$
\begin{eqnarray*}
&&|\varphi_1-\psi_1|^{r_1}=(\varphi_1-\psi_1)^{r_1}=\frac1{t^\alpha |X-Y|^\lambda z^\beta}\big(v^{\frac\kappa{r_1}}(Y)-
v_\tau^{\frac\kappa{r_1}}(Y)\big)^{r_1},\\
\text{and}
\\
&&|\varphi_2-\psi_2|^{r_1}=\frac1{t^\alpha |X-Y^\tau|^\lambda z^\beta}\big(v^{\frac\kappa{r_1}}
(Y)-v_\tau^{\frac\kappa{r_1}}(Y)\big)^{r_1}\le(\varphi_1-\psi_1)^{r_1}.
\end{eqnarray*}
Substituting the above estimates into \eqref{varphi-psi}, and then, using the mean value theorem, we obtain
\begin{eqnarray*}
0&<&u^{\frac1{r_1}}(X)-u^{\frac1{r_1}}_\tau(X)\\
&\le&2\big(\int_{\Sigma^v_\tau}\frac1{t^\alpha |X-Y|^\lambda z^\beta}\big(v^{\frac\kappa{r_1}}(Y)-v_\tau^{\frac\kappa{r_1}}(Y)\big)^{r_1}\big)^{\frac1{r_1}}\\
&\le&2\kappa\big(\int_{\Sigma^v_\tau}\frac1{t^\alpha |X-Y|^\lambda z^\beta}v^{\kappa-r_1}(Y)\big(v(Y)-v_\tau(Y)\big)^{r_1}\big)^{\frac1{r_1}},
\end{eqnarray*} for a.e. $X\in\Sigma^u_\tau.$ By \eqref{WHLSD-2} and H\"{o}lder's inequality, we have
\begin{eqnarray*}
\|u^{\frac1{r_1}}-u^{\frac1{r_1}}_\tau\|_{L^{r_1(\theta+1)}(\Sigma^u_\tau)}&\le& C\|v^{\kappa-r_1}(v-v_\tau)^{r_1}\|^{\frac1{r_1}}_{L^{\frac{\kappa+1}\kappa}(\Sigma^v_\tau)}\\
&\le&C\|v\|^{\frac\kappa{r_1}-1}_{L^{\kappa+1}(\Sigma^v_\tau)}\|v-v_\tau\|_{L^{\kappa+1}(\Sigma^v_\tau)}.
\end{eqnarray*}Combining the above into \eqref{v-}, we get \eqref{u-v-6}.

{\bf Case 3.} Assume $\kappa<1.$ Similar to the process proof of  \textbf{Case $2$},  we can obtain \eqref{u-v-7}. We omit it in here.
\end{proof}
\begin{lemma}\label{MP-1} Under the hypotheses of Theorem \ref{r-s-s}, there is $\tau_0$ sufficiently negative such that for all $\tau\le\tau_0,$ both
\begin{equation}\label{u-tau-0}
u(X)\le u_\tau(X) \quad\text{and}\quad v(Y)\le v_{\tau}(Y), \quad\text{in}~ \Sigma_\tau.
\end{equation}
\end{lemma}
\begin{proof} We firstly show $v\in L^{\kappa+1}(\mathbb{R}^{n+1}_+)$. In fact, from that condition $\frac1{\theta+1}+\frac1{\kappa+1}=\frac{\alpha+\beta+\lambda}{n+1}$, we know that $\frac1{\kappa+1}=\frac\theta{\theta+1}-\frac{n+1-(\alpha+\beta+\lambda)}{n+1}$. Using HLS inequality \eqref{WHLSD-2} and system \eqref{Euler-syst-1},  we have
\begin{eqnarray}\label{Re-v}
\big(\int_{\mathbb{R}^{n+1}_+}|v(Y)|^{\kappa+1}dY\big)^{\frac1{\kappa+1}}
&\le&\big(\int_{\mathbb{R}^{n+1}_+}(\int_{\mathbb{R}^{n+1}_+}\frac{u^\theta(X)}{ t^\alpha|X-Y|^\lambda z^\beta}dX)^{\kappa+1}dY\big)^{\frac1{\kappa+1}}\nonumber\\
&\le&\big(\int_{\mathbb{R}^{n+1}_+}u^{\theta+1}(X)dX\big)^{\frac\theta{\theta+1}}.
\end{eqnarray}Since $u\in L^{\theta+1}(\mathbb{R}^{n+1}_+),$ we deduce $v\in
L^{\kappa+1}(\mathbb{R}^{n+1}_+).$ Furthermore, we have both
$\|u\|_{L^{\theta+1}(\Sigma^u_\tau)}\rightarrow0$ and $\|v\|_{L^{\kappa+1}(\Sigma^v_\tau)}\rightarrow0$, as $\tau\rightarrow{-\infty}.$

If $\kappa\ge1$ and $\theta\ge1,$
we may choose $\tau_0$ sufficiently negative such that, for $\tau\le\tau_0$,
\[C\|v\|^{\kappa-1}_{L^{\kappa+1}(\Sigma^v_\tau)}\|u\|^{\theta-1}_{L^{\theta+1}(\Sigma^u_\tau)}\le\frac 12.\]
From \eqref{u-v-4} and \eqref{u-v-5}, if $ \tau\le\tau_0 $, we obtain both
\[\|u-u_\tau\|_{L^{\theta+1}(\Sigma^u_\tau)}\le\frac12\|u-u_\tau\|_{L^{\theta+1}(\Sigma^u_\tau)}
~\text{and}~\|v-v_\tau\|_{L^{\kappa+1}(\Sigma^v_\tau)}\le\frac12\|v-v_\tau\|_{L^{\kappa+1}(\Sigma^v_\tau)}.\]
It implies $|\Sigma^u_\tau|=|\Sigma^v_\tau|=0.$ Hence, both of $\Sigma^v_\tau$ and $\Sigma^u_\tau$ are measure zero sets.

If $\theta<1,$ then for any $\frac1\theta\le r_1\le\kappa,$ there is $\tau_0$ sufficiently
negative such that
\[C\|u\|^{\theta-\frac1{r_1}}_{L^{\theta+1}(\Sigma^u_\tau)}\|v\|^{\frac\kappa{r_1}-1}_{L^{\kappa+1}(\Sigma^v_\tau)}<\frac12\]
for $\tau\le\tau_0.$ For any such $\tau$ and $\tau_0,$ inequality \eqref{u-v-6} guarantees that
$|\Sigma^v_\tau|=0.$ Furthermore, it follows from \eqref{u-v-3} that $|\Sigma^u_\tau|=0$ for any $\tau\le\tau_0.$

Similarly, we find that $\Sigma^v_\tau$ and $\Sigma^u_\tau$ are measure zero sets for $\kappa<1.$

 Hence, combining the above argument, we arrive at
\[u(X)\le u_\tau(X) \quad\text{and}\quad v(Y)\le v_\tau(Y),\quad\text{for}~ X, Y\in\Sigma_\tau\]
for $\tau\le\tau_0$ negative enough.
It completes the proof.
\end{proof}

Define
\[\overline{\tau}=\sup\{\mu| \, u(X)\le u_\tau(X)~\text{and}~v(Y)\le v_\tau(Y),~\text{for}~ \forall\tau\le\mu, X,Y\in\Sigma_\tau\}.\]
Obviously, $\overline{\tau}$ is well defined according to Lemma \ref{MP-1}. Thus we can move the plane $x_1=\tau$ to right as long as \eqref{u-tau-0} holds. That is

\begin{lemma}\label{MP-2}
If $ \overline{\tau}\le0,$ then $u_{\overline{\tau}}\equiv u$ and $v_{\overline{\tau}}\equiv v,$  a.e. in $\Sigma_{\overline{\tau}}.$
\end{lemma}
\begin{proof} We prove it by contradiction argument. Otherwise, for $\overline{\tau}\le0,$  we show that the plane can be moved further to the right.  More precisely, there exists an $\epsilon,$ such that for all $\tau$ in
$[\overline{\tau},\overline{\tau}+\epsilon),$
\begin{equation}\label{u-v-13}
u(X)\le u_{\tau}(X)\quad\text{and}\quad v(Y)\le v_{\tau}(Y), \quad\text{in}~ \Sigma_{\overline{\tau}}.
\end{equation} In the case $v(Y)\not\equiv v_{\overline{\tau}}(Y)$ in $\Sigma_{\overline{\tau}}$, by Lemma \ref{notations} and the definition of $\overline{\tau}$, we have
$$ u(X)<u_{\overline{\tau}}(X)\quad\text{and}\quad v(Y)< v_{\overline{\tau}}(Y) ,\quad\text{in}~\Sigma_{\tau}.$$

Let
\begin{eqnarray*}
&&\Sigma^{u}_{\overline{\tau}}=\{X\in\Sigma_{\overline{\tau}}|u(X)\ge u_{\overline{\tau}}(X)\}  \quad\text{and}\quad
\Sigma^{v}_{\overline{\tau}}=\{Y\in\Sigma_{\overline{\tau}}|v(Y)\ge v_{\overline{\tau}}(Y)\}.
\end{eqnarray*}It is obvious that $\Sigma^{u}_{\overline{\tau}}$ has measure zero  set,  and $\lim_{\tau\rightarrow{\overline{\tau}}}\Sigma^{u}_{\tau}\subset\Sigma^{u}_{\overline{\tau}}$. It is same to that of $v$.

If $\kappa\geq1$ and $\theta\geq1,$ similar to \eqref{u-v-4} and \eqref{u-v-5}, we deduce
\begin{equation}\label{u-v-8}
\|u-u_\tau\|_{L^{\theta+1}(\Sigma^u_\tau)}\le C\|u\|^{\theta-1}_{L^{\theta+1}(\Sigma^u_\tau)}
\|v\|^{\kappa-1}_{L^{\kappa+1}(\Sigma^v_\tau)}\|u-u_\tau\|_{L^{\theta+1}({\Sigma^u_\tau})},
\end{equation} and
\begin{equation}\label{u-v-10}
\|v-v_\tau\|_{L^{\kappa+1}(\Sigma^v_\tau)}\le C\|u\|^{\theta-1}_{L^{\theta+1}(\Sigma^u_\tau)}
\|v\|^{\kappa-1}_{L^{\kappa+1}(\Sigma^v_\tau)}\|v-v_\tau\|_{L^{\kappa+1}({\Sigma^v_\tau})}
\end{equation}for $\tau\in[\overline{\tau},\overline{\tau}+\varepsilon)$.

 If $\kappa\ge1$ and $\theta<1,$ or $\kappa<1$ and $\theta>1,$ similar to \eqref{u-v-6} and
 \eqref{u-v-7}, we have
\begin{eqnarray}\label{u-v-11}
\|v-v_\tau\|_{L^{\kappa+1}(\Sigma^v_\tau)}&\le& C\|u\|^{\theta-\frac1{r_1}}_{L^{\theta+1}(\Sigma^u_\tau)}\|v\|^{\frac\kappa{r_1}-1}_{L^{\kappa+1}(\Sigma^v_\tau)}
\|v-v_\tau\|_{L^{\kappa+1}(\Sigma^v_\tau)},\\ \label{u-v-11a}
\|u-u_{\tau}\|_{L^{\theta+1}(\Sigma^u_\tau)}&\le& C\|v\|^{\kappa-1}_{L^{\kappa+1}(\Sigma^v_\tau)}
\|v-v_{\tau}\|_{L^{\kappa+1}({\Sigma^v_\tau})},
\end{eqnarray} and
\begin{eqnarray}\label{u-v-12}
\|u-u_\tau\|_{L^{\theta+1}(\Sigma^u_\tau)}&\le& C\|u\|^{\frac\theta{r_2}-1}_{L^{\theta+1}(\Sigma^u_\tau)}\|v\|^{\kappa-\frac1{r_2}}_{L^{\kappa+1}(\Sigma^v_\tau)}
\|u-u_\tau\|_{L^{\theta+1}(\Sigma^u_\tau)},\\ \label{u-v-12a}
\|v-v_{\tau}\|_{L^{\kappa+1}(\Sigma^v_\tau)}&\le& C\|v\|^{\theta-1}_{L^{\theta+1}(\Sigma^v_\tau)}
\|u-u_{\tau}\|_{L^{\theta+1}({\Sigma^u_\tau})}
\end{eqnarray}
for $\tau\in[\overline{\tau},\overline{\tau}+\varepsilon)$.

Since $u\in L^{\theta+1}(\mathbb{R}^{n+1}_+)$ and $v\in L^{\kappa+1}(\mathbb{R}^{n+1}_+)$
ensure that one can choose $\epsilon$ sufficiently small, such that for all $\tau$ in $[\overline{\tau},\overline{\tau}+\epsilon),$
\begin{eqnarray}\label{v-vs0}
&&C\|v\|^{\kappa-1}_{L^{\kappa+1}(\Sigma^v_\tau)}\|u\|^{\theta-1}_{L^{\theta+1}(\Sigma^u_\tau)}\le\frac 12,\\ \label{v-vs1}
& &C\|u\|^{\theta-\frac1{r_1}}_{L^{\theta+1}(\Sigma^u_\tau)}
\|v\|^{\frac\kappa{r_1}-1}_{L^{\kappa+1}(\Sigma^v_\tau)}\le\frac 12,\\ \label{v-vs2}
& &C\|u\|^{\frac\theta{r_2}-1}_{L^{\theta+1}(\Sigma^u_\tau)}
\|v\|^{\kappa-\frac1{r_2}}_{L^{\kappa+1}(\Sigma^v_\tau)}\le\frac 12.
\end{eqnarray}
Combining \eqref{u-v-8}, \eqref{u-v-10} and \eqref{v-vs0}, we have $\|u-u_{\tau}\|_{L^{\theta+1}(\Sigma^u_\tau)}=\|v-v_{\tau}\|_{L^{\kappa+1}(\Sigma^v_\tau)}=0$.
Substituting the \eqref{v-vs1} into \eqref{u-v-11}, we obtain $\|v-v_{\tau}\|_{L^{\kappa+1}(\Sigma^v_\tau)}=0.$ Furthermore, it follows from \eqref{u-v-11a} that $\|u-u_{\tau}\|_{L^{\theta+1}(\Sigma^u_\tau)}=0.$  Similarly, combining \eqref{u-v-12}, \eqref{u-v-12a} and \eqref{v-vs2}, we have $\|u-u_{\tau}\|_{L^{\theta+1}(\Sigma^u_\tau)}=\|v-v_{\tau}\|_{L^{\kappa+1}(\Sigma^v_\tau)}=0$. These imply that $|\Sigma^u_\tau|=|\Sigma^v_\tau|=0$ for $\kappa,\theta>0$.  This verifies \eqref{u-v-13},  which contradicts the definition of $\overline{\tau}$. Hence
$$u_{\overline{\tau}}\equiv u~\text{and}~ v_{\overline{\tau}}\equiv v,~\text{ a.e. ~in}~ \Sigma_{\overline{\tau}}$$
 for some $\overline{\tau}$.
Therefore, it completes the proof.
\end{proof}

\textbf{Proof of Theorem \ref{r-s-s}.} From  Lemma \ref{MP-2}, we show that if the plane stops at $x_1=\overline{\tau}$ for some $\overline{\tau}\le0$, then  $u$ and $v$ are radially  symmetric and monotone decreasing with respect to the plane $x_1=\overline{\tau}$.
%Otherwise, we can move the plane all the way $x_1=0$, so do it.
Moreover, since the $x_1$ direction can be replaced by the $x_i$ direction for $i=1,2,\cdots,n$, we deduce that $u$ and $v$ are radially symmetric and monotone decreasing with respect to $x\in\mathbb{R}^n$ about some $x_0\in\mathbb{R}^n$. This is, $u$ and $v$ are cylindrical symmetric about $t$-axis.  We complete the proof.
\qed

At the end of this section, we also give the proof of the cylindrical symmetric and monotone decreasing of positive solutions to \eqref{Euler-syst-3}.

Let $u(x,\hat{x})=f^{p-1}(x,\hat{x}),$  $v(y,\hat{y})=g^{r-1}(y,\hat{y})$ and take $\theta_1=\frac 1{p-1}, \kappa_1=q-1.$ Then, the system \eqref{Euler-syst-3} can be rewritten as
\begin{equation}\label{Euler-syst-5}
\begin{cases}
u(x,\hat{x})=\frac1{|\hat{x}|^\alpha}\int_{\mathbb{R}^{n+m}} \frac{v^{k_1}(y,\hat{y})}{|\hat{y}|^\beta|(x,\hat{x})-(y,\hat{y})|^{\lambda}}dyd\hat{y},&\quad (x,\hat{x})\in\mathbb{R}^{n+m},\\
v(y,\hat{y})=\frac1{|\hat{y}|^\beta}\int_{\mathbb{R}^{n+m}} \frac{u^{\theta_1}(x,\hat{x})}{|\hat{x}|^\alpha|(x,\hat{x})-(y,\hat{y})|^{\lambda}}dxd\hat{x},&\quad (y,\hat{y})\in\mathbb{R}^{n+m}
\end{cases}
\end{equation} with the conditions
\begin{equation}\label{WH-exp-5}
\begin{cases}
&0<\lambda<n+m,\,\kappa_1,\,\theta_1>0,\,\kappa_1\theta_1\ge1,\\
&\alpha+\beta\ge0,\, \alpha<\frac{m}{\theta_1+1},\, \beta<\frac {m}{\kappa_1+1},\\
&\frac1{\theta_1+1}+\frac1{\kappa_1+1}=\frac{ \alpha+\beta+\lambda}{n+m}.
\end{cases}
\end{equation}

The result $(i)$ of Theorem \ref{classfy-2} is included in the following theorem.

\begin{theorem}\label{Cylind-symmetry-nm} Let $\alpha,\beta,\lambda,\theta_1,\kappa_1$ satisfy \eqref{WH-exp-5}. If $(u,v)$ is a pair of positive solutions of \eqref{Euler-syst-5} with $u \in L^{\theta_1+1}(\mathbb{R}^{n+m}),$ then $u(x,\hat{x}),v(x,\hat{x})$ are radially symmetric and monotone decreasing with respect to $x$ about some $x_0$ in  $\mathbb{R}^n$; For $\alpha,\beta\ge0,$ $u(x,\hat{x}),v(x,\hat{x})$ are radially symmetric and monotone decreasing with respect to $\hat{x}$ about the origin in $\mathbb{R}^m$,
respectively. That is $u(x,\hat{x})=u(|x-x_0|,|\hat{x}|), \,  v(x,\hat{x})=v(|x-x_0|,|\hat{x}|)$.
\end{theorem}
\begin{proof} Similar to that one of Theorem \ref{r-s-s}, we can show that
$u(x,\hat{x}),v(x,\hat{x})$ are radially symmetric and monotone decreasing with respect to $x$ about some $x_0$ in  $\mathbb{R}^n$, thus we only need to show the rough sketch of $u(x,\hat{x}),v(x,\hat{x})$ are radially symmetric and monotone decreasing about the origin in $\mathbb{R}^m$ with respect to $\hat{x}$ and $\alpha,\beta\ge0$.

%For simplicity, we write $(x,x')=(x,\hat{x})\in\mathbb{R}^{n+m}$.
For any $\tau\in\mathbb{R}$, define
\[\Lambda_\tau=\{(x,\hat{x})=(x_1,x_2,\ldots x_n,\hat{x}_1,\hat{x}_2,\ldots \hat{x}_m),\hat{x}_1<\tau\}\]
and $(x,\hat{x}^\tau)=(x_1,x_2,\ldots x_n,2\tau-\hat{x}_1,\hat{x}_2,\ldots \hat{x}_m).$ We still write $u_{\tau}(x,\hat{x})=u(x,\hat{x}^\tau) $.

Let $(u,v)$ be a pair of positive solutions to \eqref{Euler-syst-5}. Similar to Lemma
\ref{notations}, it is easy to verify
\begin{eqnarray}\label{rep-nm-1}
&&u(x,\hat{x})-u_\tau(x,\hat{x})\nonumber\\
&=&\int_{\Lambda_\tau}\big[\frac1{|\hat{x}|^\alpha |(x,\hat{x})-(y,\hat{y})|^\lambda |\hat{y}|^\beta}-\frac1{|\hat{x}^\tau|^\alpha|(x,\hat{x})-(y,\hat{y}^\tau)|^\lambda|\hat{y}|^\beta}\big]
v^\kappa(y,\hat{y})dyd\hat{y}\nonumber\\
\quad&&-\int_{\Lambda_\tau}\big[\frac1{|\hat{x}^\tau|^\alpha |(x,\hat{x})-(y,\hat{y})|^\lambda|\hat{y}^\tau|^\beta}-\frac{1}{|\hat{x}|^\alpha |(x,\hat{x})-(y,\hat{y}^\tau)|^\lambda|\hat{y}^\tau|^\beta}\big]
v^\kappa_\tau(y,\hat{y}) dyd\hat{y},\nonumber\\
\end{eqnarray} and
\begin{eqnarray*}
&&v(y,\hat{y})-v_\tau(y,\hat{y})\nonumber\\
&=&\int_{\Lambda_\tau}\big[\frac1{|\hat{x}|^\alpha |(x,\hat{x})-(y,\hat{y})|^\lambda |\hat{y}|^\beta}-\frac1{|\hat{x}|^\alpha|(x,\hat{x})-(y,\hat{y}^\tau)|^\lambda|\hat{y}^\tau|^\beta}\big]
u^\theta(x,\hat{x})dxd\hat{x}\nonumber\\
\quad&&-\int_{\Lambda_\tau}\big[\frac1{|\hat{x}^\tau|^\alpha |(x,\hat{x})-(y,\hat{y})|^\lambda|\hat{y}^\tau|^\beta}-\frac{1}{|\hat{x}^\tau|^\alpha |(x,\hat{x})-(y,\hat{y}^\tau)|^\lambda|\hat{y}|^\beta}\big]
u^\theta_\tau(x,\hat{x})dxd\hat{x}.\nonumber\\
\end{eqnarray*}

Now, let $\tau\le0,$ for a.e. $(x,\hat{x})\in\Lambda_\tau.$ Using the fact $\alpha,\beta\ge0$, $|\hat{x}^\tau|\le |\hat{x}|,$ $|\hat{y}^\tau|\le|\hat{y}| $ and ${|(x,\hat{x})-(y,\hat{y})|\le
|(x,\hat{x})-(y,\hat{y}^\tau)|}$ on $\Lambda_\tau,$ it follows from \eqref{rep-nm-1} that
\begin{eqnarray*}
&&u(x,\hat{x})-u_\tau(x,\hat{x})\nonumber\\
&\le&\int_{\Lambda_\tau}\big[\frac1{|\hat{x}|^\alpha |(x,\hat{x})-(y,\hat{y})|^\lambda|\hat{y}|^\beta}-\frac1{|\hat{x}|^\alpha|(x,\hat{x})-(y,\hat{y}^\tau)|^\lambda|\hat{y}|^\beta}\big]
v^\kappa(y,\hat{y})dyd\hat{y}\nonumber\\
\quad&&-\int_{\Lambda_\tau}\big[\frac1{|\hat{x}|^\alpha |(x,\hat{x})-(y,\hat{y})|^\lambda|\hat{y}^\tau|^\beta}-\frac1{|\hat{x}|^\alpha |(x,\hat{x})-(y,\hat{y}^\tau)|^\lambda|\hat{y}^\tau|^\beta}\big]
v^\kappa_\tau(y,\hat{y})dyd\hat{y}\nonumber\\
&=&\int_{\Lambda_\tau}\frac{1}{|\hat{x}|^\alpha}\big[\frac{1}{ |(x,\hat{x})-(y,\hat{y})|^\lambda}-\frac{1}{|(x,\hat{x})-(y,\hat{y}^\tau)|^\lambda}\big]
\big(\frac{v^\kappa(y,\hat{y})}{|\hat{y}|^\beta}
-\frac{v_\tau^\kappa(y,\hat{y})}{|\hat{y}^\tau|^\beta}\big)dyd\hat{y}\nonumber\\
&\le&\int_{\Lambda_\tau}\frac{1}{|\hat{x}|^\alpha|\hat{y}|^\beta}\big[\frac{1}{ |(x,\hat{x})-(y,\hat{y})|^\lambda}-\frac{1}{|(x,\hat{x})-(y,\hat{y}^\tau)|^\lambda}\big]
\big(v^\kappa(y,\hat{y})
-v_\tau^\kappa(y,\hat{y})\big)dyd\hat{y}\nonumber\\
&\le&\int_{\Lambda_\tau}\frac{1}{|\hat{x}|^\alpha |(x,\hat{x})-(y,\hat{y})|^\lambda|\hat{y}|^\beta}
\big(v^\kappa(y,\hat{y})-v_\tau^\kappa(y,\hat{y})\big)dyd\hat{y}.
\end{eqnarray*}Similarly,
\begin{eqnarray*}
&&v(y,\hat{y})-v_\tau(y,\hat{y})\nonumber\\
&\le&\int_{\Lambda_\tau}\frac{1}{|\hat{x}|^\alpha |(x,\hat{x})-(y,\hat{y})|^\lambda|\hat{y}|^\beta}
\big(u^\theta(x,\hat{x})-u^\theta_\tau(x,\hat{x})\big)
dxd\hat{x}.\nonumber\\
\end{eqnarray*}
If $\alpha,\beta\ge0$, similar to Lemma \ref{MP-1} and Lemma \ref{MP-2}, we can employ the method of moving planes to prove the theorem. We omit the remainder of the proof in here.

\end{proof}

%~~~~~~~~~~~~~~~~~~~~~~~~~~~~~~~~~~~~~~~~~~~~~~~~~~~~~~~~~~~~~~~~~~~~~~~~~~~~
\section{\textbf{Classification results of positive solution to integral system} \label{Section 4}}
In this section, we present the explicit form on the boundary of  all nonnegative  solutions to integral equation \eqref{Euler-syst-1} via the method of moving spheres introduced by Li and Zhu in \cite{LZ1995}.

For simplicity, we write $\kappa^*=\frac{2(n+1)-\lambda-2\beta}{\lambda+2\beta}$ and
$\theta^*=\frac{2(n+1)-\lambda-2\alpha}{\lambda+2\alpha}.$

 The result $(ii)$ in Theorem \ref{classfy-1} is included in the following theorem.
\begin{theorem}\label{Solution-theo} Let $ 0<\lambda<n+1,\, \kappa,\, \theta>0, \, \kappa\theta\ge1,\alpha<\frac1{\theta+1},\, \beta<\frac1{\kappa+1}$ and $\alpha+\beta\ge0$. Assume that $(u,v)\in C^0(\overline{\mathbb{R}^{n+1}_+})\times
 C^0(\overline{\mathbb{R}^{n+1}_+})$ is a pair of positive solutions to \eqref{Euler-syst-1}, then the following two results hold.

$(i)$ If $\kappa=\kappa^*$ and $\theta=\theta^*,$
then $u,v$ must be the following forms on the boundary
\begin{eqnarray*}
u(\xi,0)=c_1\big(\frac1{|\xi-\xi_0|^2+d^2}\big)^{\frac{\lambda+2\alpha}2},\quad
v(\xi,0)=c_2\big(\frac1{|\xi-\xi_0|^2+d^2}\big)^{\frac{\lambda+2\beta}2},
\end{eqnarray*} where $c_1, c_2, d>0$ and $\xi, \xi_0\in\partial\mathbb{R}^{n+1}_+$.

$(ii)$ If $\kappa\le\frac{2(n+1)-\lambda-2\beta}{\lambda+2\beta}$ and
$\theta\le\frac{2(n+1)-\lambda-2\alpha}{\lambda+2\alpha}$ satisfying $\frac1{\theta+1}+\frac1{\kappa+1}>\frac{\alpha+\beta+\lambda}{n+1},$
%if $(u,v)$ are a pair of positive solutions to \eqref{Euler-syst-1} with $\color{blue}{u\in L^{\theta+1}(\mathbb{R}^{n+1}_+),v\in L^{\kappa+1}(\mathbb{R}^{n+1}_+)~\Red{?}}$
then $u,v=0.$
\end{theorem}

For $\tau>0$ and $(x,0)\in\partial\mathbb{R}^{n+1}_+, (\eta,z)\in \mathbb{R}^{n+1}_+$, define
\[(\eta,z)^{x,\tau}=(x,0)+\frac{\tau^2\big((\eta,z)-(x,0)\big)}{|(\eta,z)-(x,0)|^2}\]
as the Kelvin transformation of $(\eta,z)$ with respect to $B^+_\tau(x,0).$ For any such
$(\xi,t), (\eta,z)\in \mathbb{R}^{n+1}_+$, we have
\begin{equation}\label{notation-1}
|{(\xi,t)}^{x,\tau}-(\eta,z)^{x,\tau}|=\frac{\tau}{|(\xi,t)-(x,0)|}
\cdot\frac{\tau}{|(\eta,z)-(x,0)|}\cdot|(\xi,t)-(\eta,z)|,
\end{equation} and
\begin{equation*}
|(\eta,z)-(x,0)|\cdot|(\xi,t)-(\eta,z)^{x,\tau}|=|(\xi,t)-(x,0)|\cdot|(\xi,t)^{x,\tau}-(\eta,z)|,
\end{equation*} as well as
\begin{equation*}
t^{x,\tau}=\big(\frac{\tau}{|(\xi,t)-(x,0)|}\big)^2 t.
\end{equation*} We also need to define the transformation
\begin{eqnarray*}
&&\overline{u}_\tau(\xi,t)=\big(\frac\tau{|(\xi,t)-(x,0)|}\big)^{\lambda+2\alpha}u((\xi,t)^{x,\tau}), \ (\xi,t)\in\mathbb{R}^{n+1}_+\setminus\{(x,0)\},\\
&&\overline{v}_\tau(\eta,z)=\big(\frac\tau{|(\eta,z)-(x,0)|}\big)^{\lambda+2\beta}v((\eta,z)^{x,\tau}), \ (\eta,z)\in\mathbb{R}^{n+1}_+\setminus\{(x,0)\}.
\end{eqnarray*} For $R>0$, denote $\Sigma_{x,R}=\mathbb{R}^{n+1}_+\setminus\overline
{B^+_R(x,0)},$ and write $\overline{\omega}_\tau^\kappa(\eta,z) :=(\overline{\omega}_\tau(\eta,z))^\kappa.$
\begin{lemma}\label{one formula} Let $(u,v)$ be a pair of positive solutions to
\eqref{Euler-syst-1}. Then for any $(x,0)\in\partial\mathbb{R}^{n+1}_+,$ it holds
\begin{eqnarray}
\overline{u}_\tau(\xi,t)&=&\int_{\mathbb{R}^{n+1}_+}\frac {\overline{v}_\tau^\kappa(\eta,z)}{t^\alpha|(\xi,t)-(\eta,z)|^\lambda z^\beta}\big(\frac
\tau{|(\eta,z)-(x,0)|}\big)^{\mu_1}d\eta dz,\label{tow formula-1}\\
\overline{v}_\tau(\eta,z)&=&\int_{\mathbb{R}^{n+1}_+}\frac{\overline{u}_\tau^\theta(\xi,t)}{t^\alpha |(\xi,t)-(\eta,z)|^\lambda z^\beta}\big(\frac
\tau{|(\xi,t)-(x,0)|}\big)^{\mu_2}d\xi dt,\label{tow formula-2}
\end{eqnarray} where $(\xi,t),(\eta,z)\in\mathbb{R}^{n+1}_+, \, \mu_1=2n+2-(\kappa+1)(\lambda+2\beta),\ \mu_2=2n+2-(\theta+1)(\lambda+2\alpha).$ Moreover,
\begin{eqnarray}
\overline{u}_\tau(\xi,t)-u(\xi,t)&=&\int_{\Sigma_{x,\tau}}P_1\big[\big(\frac\tau{|(\eta,z)-(x,0)|}\big)^{\mu_1}\overline{v}_{\tau}^\kappa(\eta,z)-
v^\kappa(\eta,z)\big]d\eta dz,~~~~~~~~~~~~~~~~~~~~~~~~~\label{tow formula-3}\\
\overline{v}_\tau(\eta,z)-v(\eta,z)&=&\int_{\Sigma_{x,\tau}}P_2\big[\big(\frac\tau{|(\xi,t)-(x,0)|}\big)^{\mu_2}\overline{u}_{\tau}^\theta(\xi,t)-
u^\theta(\xi,t)\big]d\xi dt,~~~~~~~~~~~~~~~ \label{tow formula-4}
\end{eqnarray}
where
\begin{eqnarray*}
P_1&:=&P_1(x,\tau;(\xi,t),(\eta,z))\\
&=&\frac1{t^{\alpha}z^\beta}\big[\frac 1{|(\xi,t)-(\eta,z)|^\lambda}-\big(\frac \tau{|(\xi,t)-(x,0)|}\big)^{\lambda}\frac
1{|(\xi,t)^{x,\tau}-(\eta,z)|^\lambda}\big],\\
P_2&:=&P_2(x,\tau;(\xi,t),(\eta,z))\\
&=&\frac1{t^{\alpha}z^\beta}[\frac 1{|(\xi,t)-(\eta,z)|^\lambda}-(\frac \tau{|(\eta,z)-(x,0)|})^{\lambda}\frac
1{|(\eta,z)^{x,\tau}-(\xi,t)|^\lambda}],
\end{eqnarray*} and $P_1, P_2>0$  for $\forall(\xi,t), (\eta,z)\in\Sigma_{x,\tau}, \tau>0$.
\end{lemma}
\begin{proof} The proof of this lemma is similar to that in \cite{L2004}, we omit it in here.
\end{proof}

Define
\begin{eqnarray*}
&&\Sigma^u_{x,\tau}=\{(\xi,t)\in\Sigma_{x,\tau}\,| \, u(\xi,t)<\overline{u}_\tau(\xi,t)\},\\
&&\Sigma^v_{x,\tau}=\{(\eta,z)\in\Sigma_{x,\tau}\,|\, v(\eta,z)<\overline{v}_\tau(\eta,z)\}.
\end{eqnarray*} For simplicity, write $\overline{u}_\tau=\overline{u}_\tau(\xi,t),$ $\overline{v}_{\tau}=\overline{v}_{\tau}(\eta,z).$

 We give some estimates to carry out the method of moving spheres.}

\begin{lemma}\label{condion-estimate} Under the assumptions of Theorem \ref{Solution-theo}, let $1<q_1,q_2<\infty$, $(u,v)$ be a pair of positive solutions to \eqref{Euler-syst-1}, %with $u\in L^{\theta^*+1}(\mathbb{R}^{n+1}_+).$
 there exist some positive constants $C$ and $\tau$ such that the following estimates hold.

 $(a)$ If $\kappa\ge1$ and $\theta\ge1,$ then
\begin{eqnarray*}
\|\overline{u}_\tau-u\|_{L^{q_1}(\Sigma^u_{x,\tau})}
&\le&C\|\overline{v}_\tau\|^{\kappa-1}_{L^{\kappa+1}(\Sigma^v_{x,\tau})}
\|\overline{u}_\tau\|^{\theta-1}_{L^{\theta+1}(\Sigma^u_{x,\tau})}
\|\overline{u}_\tau-u\|_{L^{q_1}(\Sigma^u_{x,\tau})},\\
\|\overline{v}_\tau-v\|_{L^{q_2}(\Sigma^v_{x,\tau})}
&\le& C\|\overline{v}_\tau\|^{\kappa-1}_{L^{\kappa+1}(\Sigma^v_{x,\tau})}
\|\overline{u}_\tau\|^{\theta-1}_{L^{\theta+1}(\Sigma^u_{x,\tau})}
\|\overline{v}_\tau-v\|_{L^{q_2}(\Sigma^v_{x,\tau})},
\end{eqnarray*}
where $q_1,q_2$ satisfy $\frac1{q_1}-\frac1{q_2}=\frac{1}{\theta+1}-\frac{1}{\kappa+1}$.

 $(b)$ If $\kappa\ge1$ and $\theta<1,$  then there exists $s_1>1$ such that $1<\frac1{\theta}\le s_1\le\kappa$,
\begin{eqnarray*}
\|\overline{v}_\tau-v\|_{L^{q_2}(\Sigma^v_{x,\tau})}
&\le&C\|\overline{u}_\tau\|^{\theta-{\frac1{s_1}}}_{L^{\theta+1}(\Sigma^u_{x,\tau})}
\|\overline{v}_\tau\|^{\frac{\kappa}{s_1}-1}_{{L^{\kappa+1}}(\Sigma^v_{x,\tau})}\|
\overline{v}_\tau-v\|_{{L^{q_2}}(\Sigma^v_{x,\tau})},
\end{eqnarray*}where $q_1,q_2$ satisfy $\frac1{q_1}-\frac1{q_2}=\frac{1}{s_1(\theta+1)}-\frac{1}{\kappa+1}$.

 $(c)$ If $\kappa<1$ and $\theta\ge1,$ then there exists $s_2>1$ such that
$1<\frac1{\kappa}\le s_2\le\theta,$
\begin{equation*}
\|\overline{u}_\tau-u\|_{L^{q_1}(\Sigma^u_{x,\tau})}
\le C\|\overline{u}_\tau\|^{{\frac{\theta}{s_2}}-1}_{L^{\theta+1}(\Sigma^u_{x,\tau})}
\|\overline{v}_\tau\|^{\kappa-\frac 1{s_2}}_{{L^{\kappa+1}}(\Sigma^v_{x,\tau})}\|
\overline{u}_\tau-u\|_{{L^{q_1}}(\Sigma^u_{x,\tau})},~~~~~~~~~~~~~
\end{equation*}where $q_1,q_2$ satisfy $\frac1{q_1}-\frac1{q_2}=\frac{1}{\theta+1}-\frac{1}{s_2(\kappa+1)}$.
\end{lemma}
\begin{proof} Since $\kappa\le\frac{2n+2-\lambda-2\beta}{\lambda+2\beta},$ we have $\mu_1\ge0$. it follows from \eqref{tow formula-3} that
\begin{eqnarray}\label{uu-1}
0&\le&\overline{u}_\tau(\xi,t)-u(\xi,t)\nonumber\\
&=&\int_{\Sigma_{x,\tau}}P_1(x,\tau;(\xi,t),(\eta,z))\big[\big(\frac\tau{|(\eta,z)-(x,0)|}\big)^{\mu_1}\overline{v}_\tau^\kappa(\eta,z)-
v^\kappa(\eta,z)\big]d\eta dz\nonumber\\
&\le&\int_{\Sigma^v_{x,\tau}}P_1(x,\tau;(\xi,t),(\eta,z))\big(\overline{v}_\tau^\kappa(\eta,z)-v^\kappa(\eta,z)\big)d\eta dz\nonumber\\
&\le&\int_{\Sigma^v_{x,\tau}}\frac {\overline{v}_\tau^\kappa(\eta,z)-v^\kappa(\eta,z)}{t^\alpha z^\beta|(\xi,t)-(\eta,z)|^\lambda}d\eta dz.
\end{eqnarray}  Since $\mu_2\ge0$ for $\theta\le\frac{2n+2-\lambda-2\alpha}{\lambda+2\alpha},$ similar to the above,  we have for $(\eta,z)\in\Sigma^v_{x,\tau},$
\begin{eqnarray}\label{vv-1}
0&\le&\overline{v}_\tau(\eta,z)-v(\eta,z)\nonumber\\
&=&\int_{\Sigma_{x,\tau}}P_2(x,\tau,(\xi,t),(\eta,z))\big[\big(\frac\tau{|(\xi,t)-(x,0)|}\big)^{\mu_2}\overline{u}_\tau^\theta(\xi,t)-
u^\theta(\xi,t)\big]d\xi dt\nonumber\\
&\le&\int_{\Sigma^u_{x,\tau}}P_2(x,\tau,(\xi,t),(\eta,z))\big(\overline{u}_\tau^\theta(\xi,t)-u^\theta(\xi,t)\big)d\xi dt\nonumber\\
&\le&\int_{\Sigma^u_{x,\tau}}\frac {\overline{u}_\tau^\theta(\xi,t)-u^\theta(\xi,t)}{t^\alpha z^\beta|(\xi,t)-(\eta,z)|^\lambda}d\xi dt.
\end{eqnarray}
%\begin{eqnarray}\label{uu-1}
%0&\le&\overline{u}_\tau(\xi,t)-u(\xi,t)\nonumber\\
%&=&\int_{\Sigma_{x,\tau}}P_1(x,\tau;(\xi,t),(\eta,z))(\overline{v}_\tau^{\kappa^*}(\eta,z)-
%v^{\kappa^*}(\eta,z))d\eta dz\nonumber\\
%&\le&\int_{\Sigma^v_{x,\tau}}P_1(x,\tau;(\xi,t),(\eta,z))(\overline{v}_\tau^{\kappa^*}(\eta,z)-v^{\kappa^*}(\eta,z))d\eta dz\nonumber\\
%&\le&\int_{\Sigma^v_{x,\tau}}\frac {\overline{v}_\tau^{\kappa^*}(\eta,z)-v^{\kappa^*}(\eta,z)}{t^\alpha z^\beta|(\xi,t)-(\eta,z)|^\lambda}d\eta dz.
%\end{eqnarray} Note that $\theta=\theta^*=\frac{2n+2-\lambda-2\alpha}{\lambda+2\alpha},$ we have $\mu_2=0.$ A similar computation, for $(\eta,z)\in\Sigma^v_{x,\tau},$
%\begin{eqnarray}\label{vv-1}
%0&\le&\overline{v}_\tau(\eta,z)-v(\eta,z)\nonumber\\
%&=&\int_{\Sigma_{x,\tau}}P_2(x,\tau,(\xi,t),(\eta,z))(\overline{u}_\tau^{\theta^*}(\xi,t)-
%u^{\theta^*}(\xi,t))d\xi dt\nonumber\\
%&\le&\int_{\Sigma^u_{x,\tau}}P_2(x,\tau,(\xi,t),(\eta,z))(\overline{u}_\tau^{\theta^*}(\xi,t)-u^{\theta^*}(\xi,t))d\xi dt\nonumber\\
%&\le&\int_{\Sigma^u_{x,\tau}}\frac {\overline{u}_\tau^{\theta^*}(\xi,t)-u^{\theta^*}(\xi,t)}{t^\alpha z^\beta|(\xi,t)-(\eta,z)|^\lambda}d\xi dt.
%\end{eqnarray}
We split the remainder of the proof into two steps.

{\bfseries Step 1.} We discuss the case $\kappa\ge1.$ By the mean value theorem and \eqref{uu-1}, we have
\[0\le\overline{u}_\tau(\xi,t)-u(\xi,t)\le\kappa\int_{\Sigma^v_{x,\tau}}\frac{\overline{v}_\tau^{\kappa-1}(\eta,z)
(\overline{v}_\tau(\eta,z)-v(\eta,z))}{t^\alpha|(\xi,t)-(\eta,z)|^\lambda z^\beta}d\eta dz.\]
Let $q_1,p_1>1$ satisfy $\frac1{q_1}=\frac1{p_1}-\frac{n+1-(\alpha+\beta+\lambda)}{n+1}.$
Using inequality \eqref{WHLSD-2}, we have
\[\|(\overline{u}_{\tau}-u)_+\|_{L^{q_1}(\mathbb{R}^{n+1}_+)}\le C\|\overline{v}^{\kappa-1}_\tau(\overline{v}_\tau-v)_+\|_{L^{p_1}(\mathbb{R}^{n+1}_+)},\]
where $f_+(x,0)=\max\{f(x,0),0\}$. Assume that $q_2>1$ satisfies $\frac{\kappa-1}{\kappa+1}+\frac1{q_2}=\frac1{p_1}.$ By H\"{o}lder's inequality,  we obtain
\[\|\overline{v}^{\kappa-1}_\tau(\overline{v}_\tau-v)\|_{L^{p_1}(\Sigma^v_{x,\tau})}
\le\|\overline{v}_\tau\|^{\kappa-1}_{L^{\kappa+1}(\Sigma^v_{x,\tau})}\|\overline{v}_\tau-v\|_{L^{q_2}(\Sigma^v_{x,\tau})}.\]
We are also easy to check the fact $\frac1{q_1}-\frac1{q_2}=\frac{1}{\theta+1}-\frac{1}{\kappa+1}.$

Thus
\begin{equation}\label{u-formula}
\|\overline{u}_\tau-u\|_{L^{q_1}(\Sigma^u_{x,\tau})}\le C\|\overline{v}_\tau\|^{\kappa-1}_{L^{\kappa+1}
(\Sigma^v_{x,\tau})}\|\overline{v}_\tau-v\|_{L^{q_2}(\Sigma^v_{x,\tau})}.
\end{equation} We divide $\theta$ into two following cases to discuss $\overline{v}_\tau(\eta,z)-v(\eta,z)$.

{\bfseries Case 1.} $\kappa\ge1$ and $\theta\ge1.$ In this case, for any $(\eta,z)\in\Sigma^v_{x,\tau},$ by the mean value theorem and \eqref{vv-1} that
\[\overline{v}_\tau(\eta,z)-v(\eta,z)\le\theta\int_{\Sigma^u_{x,\tau}}\frac{\overline{u}^{\theta-1}_\tau(\xi,t)
(\overline{u}_\tau(\xi,t)-u(\xi,t))}{t^\alpha|(\eta,z)-(\xi,t)|^\lambda z^\beta}d\xi dt.\]
Now we choose $p_2,  q_2>1$ such that $\frac1{q_2}=\frac1{p_2}-\frac{n+1-(\alpha+\beta+\lambda)}{n+1}$ and $\frac{\theta-1}{\theta+1}+\frac1{q_1}=\frac1{p_2}.$  Similar to the above, using inequality \eqref{WHLSD-2}  and H\"{o}lder's inequality again, we have
\begin{eqnarray}\label{v-formula}
\|\overline{v}_\tau-v\|_{L^{q_2}(\Sigma^v_{x,\tau})}&\le&C \|\overline{u}^{\theta-1}_\tau(\overline{u}_\tau-u)\|_{L^{p_2}(\Sigma^u_{x,\tau})}\nonumber\\
&\le&C\|\overline{u}_\tau\|^{\theta-1}_{L^{\theta+1}(\Sigma^u_{x,\tau})}\| \overline{u}_\tau-u\|_{L^{q_1}(\Sigma^u_{x,\tau})}.
\end{eqnarray}
Combining \eqref{u-formula} and \eqref{v-formula}, we have
\begin{equation*}
\|\overline{u}_\tau-u\|_{L^{q_1}(\Sigma^u_{x,\tau})}
\le C\|\overline{v}_\tau\|^{\kappa-1}_{L^{\kappa+1}(\Sigma^v_{x,\tau})}
\|\overline{u}_\tau\|^{\theta-1}_{L^{\theta+1}(\Sigma^u_{x,\tau})}
\|\overline{u}_\tau-u\|_{L^{q_1}(\Sigma^u_{x,\tau})}.
\end{equation*} In the same way, we have
\begin{equation*}
\|\overline{v}_\tau-v\|_{L^{q_2}(\Sigma^v_{x,\tau})}
\le C\|v\|^{\kappa-1}_{L^{\kappa+1}(\Sigma^v_{x,\tau})}
\|u\|^{\theta-1}_{L^{\theta+1}(\Sigma^u_{x,\tau})}\|\overline{v}_\tau-v\|_{L^{q_2}(\Sigma^v_{x,\tau})}.
\end{equation*}
  Thus, we obtain the estimate $(a)$.

{\bfseries Case 2.} $\kappa\ge1$ and $\theta<1.$ Since $\frac1{\theta}\le \kappa,$ we can choose $s_1>1$ such that $1<\frac1{\theta}\le s_1\le\kappa$. For any $(\eta,z)\in\Sigma^v_{x,\tau},$ it follows from the mean value theorem and \eqref{vv-1} that
\begin{equation*}
0\le\overline{v}_\tau(\eta,z)-v(\eta,z)\le s_1\theta\int_{\Sigma^u_{x,\tau}}\frac{\overline{u}_\tau^{\theta-\frac1{s_1}}(\xi,t)
\big(\overline{u}^{\frac1{s_1}}_\tau(\xi,t)-u^{\frac1{s_1}}(\xi,t)\big)}{t^\alpha |(\eta,z)-(\xi,t)|^\lambda z^\beta}d\xi dt.
\end{equation*}
Assume that $q_1,q_2$ and $ p_2$ satisfy $\frac1{q_2}=\frac1{p_2}-\frac{n+1-(\alpha+\beta+\lambda)}{n+1}$ and $\frac{\theta-\frac1{s_1}}{\theta+1}+\frac1{q_1}=\frac1{p_2}$ in this case.
Invoking inequality \eqref{WHLSD-2} and H\"{o}lder's inequality, we have
\begin{eqnarray}\label{vformula-1}
\|\overline{v}_\tau-v\|_{L^{q_2}(\Sigma^v_{x_,\tau})}&\le&C\|\overline{u}^{\theta-\frac1{s_1}}_\tau \big(\overline{u}^{\frac1{s_1}}_\tau-u^{\frac1{s_1}}\big)\|_{L^{p_2}(\Sigma^u_{x,\tau})}\nonumber\\
&\le&C\|\overline{u}_\tau\|^{\theta-{\frac1{s_1}}}_{L^{\theta+1}(\Sigma^u_{x,\tau})} \|\overline{u}^{\frac1{s_1}}_\tau-u^{\frac1{s_1}}\|_{L^{q_1}(\Sigma^u_{x,\tau})}.
\end{eqnarray}
Assume that $p_1,q_1$ and $q_2$ satisfy $\frac1{q_1}=\frac1{p_1}-\frac{n+1-(\alpha+\beta+\lambda)}{n+1},$ and $\frac{\frac\kappa{s_1}-1}{\kappa+1}+\frac{1}{q_2}=\frac{1}{p_1}$.  Similar to that one of $(b)$ in Lemma \ref{k-theta-u}, we have
\begin{eqnarray}\label{u-fm-2}
\|\overline{u}_\tau^{\frac1{s_1}}-u^{\frac1{s_1}}\|_{L^{q_1}(\Sigma^v_{x,\tau})}&\le& C\|\overline{v}^{\kappa-s_1}(\overline{v}_\tau-v)^{s_1}\|^{\frac1{s_1}}_{{L^{p_1}}(\Sigma^v_{x,\tau})}\nonumber\\
&\le&C\|\overline{v}_\tau\|^{\frac{\kappa}{s_1}-1}_{{L^{\kappa+1}}(\Sigma^v_{x,\tau})}\| \overline{v}_\tau-v\|_{{L^{q_2}}(\Sigma^v_{x,\tau})}.
\end{eqnarray} Noting that $\frac1\kappa\le\frac1{s_1}\le\theta$ and combining \eqref{vformula-1} and \eqref{u-fm-2}, we arrive at
\begin{eqnarray*}
\|\overline{v}_\tau-v\|_{L^{q_2}(\Sigma^v_{x,\tau})}&\le& C\|\overline{u}^{\theta-\frac1{s_1}}_\tau(\overline{u}^{\frac1{s_1}}_\tau-u^{\frac1{s_1}})\| _{L^{p_2}\Sigma^v_{x,\tau})}\nonumber\\
&\le&C\|\overline{u}_\tau\|^{\theta-{\frac1{s_1}}}_{L^{\theta+1}(\Sigma^u_{x,\tau})}
\|\overline{v}_\tau\|^{\frac{\kappa}{s_1}-1}_{{L^{\kappa+1}}(\Sigma^v_{x,\tau})}\|
\overline{v}_\tau-v\|_{{L^{q_2}}(\Sigma^v_{x,\tau})}.
\end{eqnarray*}
 The estimate $(b)$ is showed.

{\bfseries Step 2.} We next discuss the case $\kappa<1.$  Since $\frac1{\theta}\le \kappa<1,$ we have $\theta\ge1,$ and then we can choose $s_2>1$ such that
$1<\frac1{\kappa}\le s_2\le\theta.$

Similar the process proof of \textbf{Case $2$} in step $1$, one has
\begin{equation*}
\|\overline{u}_\tau-u\|_{L^{q_1}(\Sigma^u_{x,\tau})}\le C\|\overline{u }_\tau\|^{{\frac{\theta}{s_2}}-1}_{L^{\theta+1}(\Sigma^u_{x,\tau})}\|\overline{v}_\tau\|^{\kappa-\frac 1{s_2}}_{{L^{\kappa+1}}(\Sigma^v_{x,\tau})}\|
\overline{u}_\tau-u\|_{{L^{q_1}}(\Sigma^u_{x,\tau})}.
\end{equation*} The details of this computation are omitted. The estimate $(c)$ is showed and the lemma is proved.
\end{proof}
\begin{lemma} \label{condion-1} Under the assumptions of Theorem \ref{Solution-theo},
suppose that $(u,v)$ is a  pair of positive solutions to \eqref{Euler-syst-1},
%with
%$\theta\le\frac{2n+2-\lambda-2\alpha}{\lambda+2\alpha}$ and $\frac1{\theta+1}+\frac1{\kappa+1}>\frac{\alpha+\beta+\lambda}{n+1}.$
then there exists $\tau_0(x,0)>0,$ such that for
$0<\tau<\tau_0(x,0)$ with $(x,0)\in\partial\mathbb{R}^{n+1}_+$,
\begin{equation}\label{MS-1}
\overline{u}_\tau(\xi,t)\le u(\xi,t),~\overline{v}_\tau(\eta,z)\le v(\eta,z), ~a.e.~\text{in}~\Sigma_{x,\tau}.
\end{equation}
\end{lemma}

\begin{proof}
We divide the proof into three cases.

$(i)$ Case $\kappa\ge1$ and $\theta\ge1.$ From $(a)$ of Lemma \ref{condion-estimate} and Kelvin transformation, we have
\begin{eqnarray}\label{u-f1-2}
\|\overline{u}_\tau-u\|_{L^{q_1}(\Sigma^u_{x,\tau})}
&\le&C\|\overline{v}_\tau\|^{\kappa-1}_{L^{\kappa+1}(\Sigma^v_{x,\tau})}
\|\overline{u}_\tau\|^{\theta-1}_{L^{\theta+1}(\Sigma^u_{x,\tau})}\|\overline{u}_\tau-u\|_{L^{q_1}(\Sigma^u_{x,\tau})}\nonumber\\
&\le&C\|v\|^{\kappa-1}_{L^{\kappa+1}(B^+_\tau{(x,0)})}
\|u\|^{\theta-1}_{L^{\theta+1}(B^+_\tau{(x,0)})}\|\overline{u}_\tau-u\|_{L^{q_1}(\Sigma^u_{x,\tau})},\nonumber\\
\end{eqnarray} and
\begin{equation}\label{v-fm-2}
\|\overline{v}_\tau-v\|_{L^{q_2}(\Sigma^v_{x,\tau})}
\le C\|v\|^{\kappa-1}_{L^{\kappa+1}(B^+_\tau(x,0))}
\|u\|^{\theta-1}_{L^{\theta+1}(B^+_\tau(x,0))}\|\overline{v}_\tau-v\|_{L^{q_2}(\Sigma^v_{x,\tau})}.
\end{equation}
%Since $u\in C^0(\overline{\mathbb{R}^{n+1}_+})$ and $v\in C^0(\overline{\mathbb{R}^{n+1}_+})$,  we follow that $ u$ and $v$ are bounded  on $\overline{B^+_\tau(x,0)}$.
Since $u\in L^{\theta+1}(\mathbb{R}^{n+1}_+)$ and $v\in L^{\kappa+1}(\mathbb{R}^{n+1}_+)$,
 we can choose $\tau<\tau_0(x,0)$ small enough such that
\[C\|v\|^{\kappa-1}_{L^{\kappa+1}(B^+_\tau(x,0))}\|u\|^{\theta-1}_{L^{\theta+1}(B^+_\tau(x,0))}
\le C \tau^{n(\frac{\kappa-1}{\kappa+1}+\frac{\theta-1}{\theta+1})} \le\frac12.\]
Combining the above with \eqref{u-f1-2} and \eqref{v-fm-2}, we get
\[\|\overline{u}_\tau-u\|_{L^{q_1}(\Sigma^u_{x,\tau})}\le\frac 12\|\overline{u}_\tau-u\|_{L^{q_1}(\Sigma^u_{x,\tau})},\]
and
\[\|\overline{v}_\tau-v\|_{L^{q_2}(\Sigma^v_{x,\tau})}\le\frac 12\|\overline{v}_\tau-v\|_{L^{q_2}(\Sigma^v_{x,\tau})}.\]
This implies $\|\overline{u}_\tau-u\|_{L^{q_1}(\Sigma^u_{x,\tau})}=\|\overline{v}_\tau-v\|_{L^{q_2}(\Sigma^v_{x,\tau})}=0,$ then $\Sigma^u_{x,\tau}$ and $\Sigma^v_{x,\tau}$ have measure zero sets. Hence, \eqref{MS-1} holds.

$(ii)$ Case $\kappa\ge1$ and $\theta<1.$ From $(b)$ of Lemma \ref{condion-estimate} and Kelvin transformation,  we have
\begin{eqnarray}\label{vl}
\|\overline{v}_\tau-v\|_{L^{q_2}(\Sigma^v_{x,\tau})}
\le C\| u\|^{\theta-{\frac1{s_1}}}_{L^{\theta+1}(B^+_\tau(x,0))}
\|v\|^{\frac{\kappa}{s_1}-1}_{{L^{\kappa+1}}(B^+_\tau(x,0))}\|
\overline{v}_\tau-v\|_{{L^{q_2}}(\Sigma^v_{x,\tau})}.~~~~~~~~~~~~~~~
\end{eqnarray}
Since $u\in L^{\theta+1}(\mathbb{R}^{n+1}_+)$ and $v\in L^{\kappa+1}(\mathbb{R}^{n+1}_+)$, there exists $\tau_0(x,0)$ small enough such that for $0<\tau<\tau_0(x,0),$
\begin{equation}\label{uv-s}
C\|u\|^{\theta-{\frac1{s_1}}}_{L^{\theta+1}(B^+_\tau(x,0))}
\|v\|^{\frac{\kappa}{s_1}-1}_{{L^{\kappa+1}}(B^+_\tau(x,0))}
\le C \tau^{n(\frac{\frac{\kappa}{s_1}-1}{\kappa+1}+\frac{\theta-{\frac1{s_1}}}{\theta+1})}
\le\frac12.
\end{equation}
Substituting \eqref{uv-s} into \eqref{vl}, we arrive at
\[\|\overline{v}_\tau-v\|_{L^{q_2}(\Sigma^v_{x,\tau})}\le\frac12\|\overline{v}_\tau-v\|_{L^{q_2}(\Sigma^v_{x,\tau})}.\]
That is
\begin{equation*}
\|\overline{v}_\tau-v\|_{L^{q_2}(\Sigma^v_{x,\tau})}=0.
\end{equation*}
This implies that $\Sigma^v_{x,\tau}$ must be measure $0.$
Furthermore, it follows from \eqref{u-formula} that
\begin{equation}\label{u-formula-2}
\|\overline{u}_\tau-u\|_{L^{q_1}(\Sigma^u_{x,\tau})}\le C\|v\|^{\kappa-1}_{L^{\kappa+1}
(B^+_\tau(x,0))}\|\overline{v}_\tau-v\|_{L^{q_2}(\Sigma^v_{x,\tau})}.
\end{equation}
Combining the above, we arrive at
\[\|\overline{u}_\tau-u\|_{L^{q_1}(\Sigma^u_{x,\tau})}=0.\]
This shows that both $\Sigma^u_{x,\tau}$ and $\Sigma^v_{x,\tau}$ must be measure $0.$ Hence, \eqref{MS-1} holds.

$(iii) $ Case $\kappa<1$. Similarly, using $(c)$ of Lemma \ref{condion-estimate}, we find that $\Sigma^u_{x,\tau}$ and $\Sigma^v_{x,\tau}$ are  measure  zero sets. And then we obtain \eqref{MS-1} holds in this case.
The proof is complete.
\end{proof}

Now, define
\[\overline{\tau}(x,0)=\sup\{\mu>0\,|\, \overline{u}_\tau\le u(\xi,t),
\overline{v}_\tau\le v(\eta,z),\forall\tau\in(0,\mu),
(\xi,t),(\eta,z)\in\Sigma_{x,\tau}\}.\]

From Lemma \ref{condion-1}, we know that $\overline{\tau}(x_0,0)$ is well defined. So the sphere can be moved to the outside.
\begin{lemma} \label{Equv-Lemma}
If $\overline{\tau}(x_0,0)<\infty$ for some $(x_0,0)\in\partial\mathbb{R}^{n+1}_+,$ then
\begin{equation}\label{K-E}
\overline{u}_{\overline{\tau}(x_0,0)}=u,\quad\text{and}\quad\overline{v}_{\overline{\tau}(x_0,0)}= v, ~ \text{on}~ \mathbb{R}^{n+1}_+.
\end{equation}
Furthermore, it must
\[\kappa=\kappa^*=\frac{2n+2-\lambda-2\beta}{\lambda+2\beta}\quad \text{and}\quad \theta=\theta^*=\frac{2n+2-\lambda-2\alpha}{\lambda+2\alpha}.\]
\end{lemma}
\begin{proof} The proof is similar to Lemma 3.4 in \cite{DZ2015a}. For simplicity, we always write $\overline{\tau}=\overline{\tau}(x_0,0)$ in the following section. We only need to show
\begin{eqnarray*}
&&\overline{u}_{\overline{\tau}}(\xi,t)=u(\xi,t),\quad \forall(\xi,t)\in\Sigma_{x_0,\overline{\tau}},\\
&&\overline{v}_{\overline{\tau}}(\eta,z)=v(\eta,z),\quad \forall(\eta,z)\in\Sigma_{x_0,\overline{\tau}}.
\end{eqnarray*} By the definition of $\overline{\tau},$ we have
\begin{eqnarray*}
&&\overline{u}_{\overline{\tau}}(\xi,t)\le u(\xi,t),\quad \forall(\xi,t)\in\Sigma_{x_0,\overline{\tau}},\\
&&\overline{v}_{\overline{\tau}}(\eta,z)\le v(\eta,z),\quad \forall(\eta,z)\in\Sigma_{x_0,\overline{\tau}}.
\end{eqnarray*} If $\overline{u}_{\overline{\tau}}(\xi,t)\not\equiv u(\xi,t)$ or $\overline{v}_{\overline{\tau}}(\eta,z)\not\equiv v(\eta,z)$, we know from \eqref{tow formula-3} and \eqref{tow formula-4} (using $(x_0,0),\overline{\tau}$ to replace $(x,0),\tau$), that
\begin{eqnarray*}
&&\overline{u}_{\overline{\tau}}(\xi,t)<u(\xi,t),\quad \forall(\xi,t)\in\Sigma_{x_0,\overline{\tau}},\\
&&\overline{v}_{\overline{\tau}}(\eta,z)<v(\eta,z),\quad \forall(\eta,z)\in\Sigma_{x_0,\overline{\tau}}.
\end{eqnarray*} For fixed $R,\varepsilon_0$ and any $\delta>0,$  there exist $c_1,c_2$ such that
\begin{eqnarray*}
&&u(\xi,t)-\overline{u}_{\overline{\tau}}(\xi,t)>c_1,\quad \forall(\xi,t)\in\Sigma_{x_0,\overline{\tau}+\delta}\cap B^+_R(x_0,0),\\
&&v(\eta,z)-\overline{v}_{\overline{\tau}}(\eta,z)>c_2,\quad \forall(\eta,z)\in\Sigma_{x_0,\overline{\tau}+\delta}\cap B^+_R(x_0,0).
\end{eqnarray*} Invoking \eqref{tow formula-3} and \eqref{tow formula-4}, we can choose $\varepsilon<\delta$ sufficiently small such that for $\tau\in[\overline{\tau},\overline{\tau}+\varepsilon),$
\begin{eqnarray*}
&&u(\xi,t)\ge\overline{u}_{\tau}(\xi,t),\quad \forall(\xi,t)\in\Sigma_{x_0,\tau+\delta}\cap B^+_R(x_0,0),\\
&&v(\eta,z)\ge \overline{v}_{\tau}(\eta,z),\quad \forall(\eta,z)\in\Sigma_{x_0,\tau+\delta}\cap B^+_R(x_0,0).
\end{eqnarray*} This implies that $\Sigma^u_{x_0,\tau}$ and $\Sigma^v_{x_0,\tau}$ have no intersection with $\Sigma_{x_0,\tau+\delta}\cap B^+_R(x_0,0).$ Thus, $\Sigma^u_{x_0,\tau}$ and $\Sigma^v_{x_0,\tau}$ are contained in the union of
\[\mathbb{R}^{n+1}_+\setminus B^+_{R}(x_0,0)\quad\text{and}\quad \Sigma_{x_0,\tau}\setminus\Sigma_{x_0,\tau+\delta}.\]
For simplicity, write
\[\Omega_{\tau,R}=(\mathbb{R}^{n+1}_+\setminus
B^+_{R}(x_0,0))\cup(\Sigma_{x_0,\tau}\setminus\Sigma_{x_0,\tau+\delta}),\]
and the reflection of $\Omega_{\tau,R}$ is defined as
\[(\Omega_{\tau,R})^*=B^+_{\varepsilon_1}(x_0,0)\cup(B_\tau^+(x_0,0)\setminus B^+_{\tau^2/(\tau+\delta)}(x_0,0)),\]
where $\varepsilon_1=\tau/R$ is small as $R\rightarrow\infty.$

If $\kappa\ge1$ and $ \theta\geq 1,$ similar to \eqref{u-f1-2} and \eqref{v-fm-2}, we have
\begin{equation}\label{u-formula-3}
\|\overline{u}_\tau-u\|_{L^{q_1}(\Sigma^u_{x_0,\tau})}\le
C\|v\|^{\kappa-1}_{L^{\kappa+1}((\Omega_{\tau,R})^*)}
\|u\|^{\theta-1}_{L^{\theta+1}((\Omega_{\tau,R})^*)}
\|\overline{u}_\tau-u\|_{L^{q_1}(\Sigma^u_{x_0,\tau})},
\end{equation} and
\begin{equation}\label{v-formula-3}
\|\overline{v}_\tau-v\|_{L^{q_2}(\Sigma^v_{x_0,\tau})}\le
C\|v\|^{\kappa-1}_{L^{\kappa+1}((\Omega_{\tau,R})^*)}
\|u\|^{\theta-1}_{L^{\theta+1}((\Omega_{\tau,R})^*)}\|\overline{v}_\tau-v\|_{L^{q_2}(\Sigma^v_{x_0,\tau})}
\end{equation}for $\tau\in[\bar{\tau},\bar{\tau}+\varepsilon)$.

If $\kappa\ge1$ and $\theta<1,$ or $\kappa<1$ and $\theta\ge1,$ similar to \eqref{vl} and \eqref{u-formula-2}, we obtain
\begin{eqnarray}
\|\overline{v}_\tau-v\|_{L^{q_2}(\Sigma^v_{x_0,\tau})}
&\le& C\|u\|^{\theta-{\frac1{s_1}}}_{L^{\theta+1}((\Omega_{\tau,R})^*)}
\|v\|^{\frac{\kappa}{s_1}-1}_{{L^{\kappa+1}}((\Omega_{\tau,R})^*)}
\|\overline{v}_\tau-v\|_{{L^{q_2}}(\Sigma^v_{x_0,\tau})},\nonumber\\
\label{u-formula-4}\\
\|\overline{u}_\tau-u\|_{L^{q_1}(\Sigma^u_{x_0,\tau})}&\le& C\|v\|^{\kappa-1}_{L^{\kappa+1}
((\Omega_{\tau,R})^*)}\|\overline{v}_\tau-v\|_{L^{q_2}(\Sigma^v_{x_0,\tau})},\label{u-formula-5}
\end{eqnarray} and
\begin{eqnarray}
\|\overline{u}_\tau-u\|_{L^{q_1}(\Sigma^u_{x_0,\tau})}
&\le& C\|u\|^{{\frac{\theta}{s_2}}-1}_{L^{\theta+1}((\Omega_{\tau,R})^*)}\|v\|^{\kappa-\frac 1{s_2}}_{{L^{\kappa+1}}((\Omega_{\tau,R})^*)}
\|\overline{u}_\tau-u\|_{{L^{q_1}}(\Sigma^u_{x_0,\tau})}, \nonumber\\
\label{u-formula-6}\\
\|\overline{v}_\tau-v\|_{L^{q_2}(\Sigma^v_{x_0,\tau})}&\le& C\|u\|^{\theta-1}_{L^{\theta+1}
((\Omega_{\tau,R})^*)}\|\overline{u}_\tau-u\|_{L^{q_1}(\Sigma^u_{x_0,\tau})}\label{u-formula-7}
\end{eqnarray} for $\tau\in[\overline{\tau},\overline{\tau}+\varepsilon).$
Since $u\in L^{\theta+1}(\mathbb{R}^{n+1}_+)$ and $v\in L^{\kappa+1}(\mathbb{R}^{n+1}_+)$, for $\delta$ small enough and $R$ large enough,  we have
\[\int_{(\Omega_{\tau,R})^*}u^{\theta+1}(\xi,t)d\xi dt<\varepsilon_0(\varepsilon,\delta),\quad \int_{(\Omega_{\tau,R})^*}v^{\kappa+1}(\eta,z)d\eta dz<\varepsilon_0(\varepsilon,\delta).\]
Choosing $\varepsilon_0$ small enough (for $\varepsilon,\delta$ sufficiently small) such that for $\tau\in[\overline{\tau},\overline{\tau}+\varepsilon),$
\begin{eqnarray}\label{uvE0}
C\|v\|^{\kappa-1}_{L^{\kappa+1}((\Omega_{\tau,R})^*)}
\|u\|^{\theta-1}_{L^{\theta+1}((\Omega_{\tau,R})^*)}&\le&\frac 12,\\ \label{uvE1}
C\|u\|^{\theta-{\frac1{s_1}}}_{L^{\theta+1}((\Omega_{\tau,R})^*)}
\|v\|^{\frac{\kappa}{s_1}-1}_{{L^{\kappa+1}}((\Omega_{\tau,R})^*)}&\le&\frac12,\\
\label{uvE2}
C\|u\|^{{\frac{\theta}{s_2}}-1}_{L^{\theta+1}((\Omega_{\tau,R})^*)}
\|v\|^{\kappa-\frac 1{s_2}}_{{L^{\kappa+1}}((\Omega_{\tau,R})^*)}&\le&\frac12.
\end{eqnarray}
Combining \eqref{u-formula-3}, \eqref{v-formula-3} and \eqref{uvE0}, we obtain
\begin{equation}\label{uve3}
\|\overline{u}_\tau-u\|_{L^{q_1}(\Sigma^u_{x_0,\tau})}=\|\overline{v}_\tau-v\|_{L^{q_2}(\Sigma^v_{x_0,\tau})}=0.
\end{equation}
Substituting \eqref{uvE1} into \eqref{u-formula-4},  we have
$$\|\overline{v}_\tau-v\|_{L^{q_2}(\Sigma^v_{x_0,\tau})}=0.$$
Furthermore, combining the above and  \eqref{u-formula-5}, it holds
$$\|\overline{u}_\tau-u\|_{L^{q_1}(\Sigma^u_{x_0,\tau})}=0.$$
Similarly, from \eqref{u-formula-6}, \eqref{u-formula-7} and \eqref{uvE2}, we also get \eqref{uve3} holds. This  shows that both $\Sigma^u_{x_0,\tau}$ and $\Sigma^v_{x_0,\tau}$ must be measure $0.$ Hence, we can conclude that
\begin{eqnarray*}
&&\overline{u}_{\tau}(\xi,t)\le u(\xi,t),~ \forall(\xi,t)\in\Sigma_{x_0,\tau},\\
&&\overline{v}_{\tau}(\eta,z)\le v(\eta,z),~\forall(\eta,z)\in\Sigma_{x_0,\tau}
\end{eqnarray*} for $\tau\in[\overline{\tau},\overline{\tau}+\varepsilon)$, which contradicts  the definition of $\overline{\tau}.$

Finally, combining \eqref{K-E} with \eqref{tow formula-1} and \eqref{tow formula-2}, $\kappa,\theta$ must be the forms $\kappa=\kappa^*$ and $\theta=\theta^*$.
 The lemma is proved.
\end{proof}
\begin{lemma}\label{Equv not-Lemma} If
$\kappa\le\frac{2n+2-\lambda-2\beta}{\lambda+2\beta}$ and $\theta\le\frac{2n+2-\lambda-2\alpha}{\lambda+2\alpha}$ satisfying $\frac1{\theta+1}+\frac1{\kappa+1}>\frac{\alpha+\beta+\lambda}{n+1},$ then
$\overline{\tau}(x_0,0)=\infty$ for all $(x_0,0)\in\partial\mathbb{R}^{n+1}_+$.
\end{lemma}
\begin{proof} We prove it by contradiction argument. Suppose the contrary, there exists some $(x_0,0)\in\partial\mathbb{R}^{n+1}_+$ such that $\overline{\tau}(x_0,0)<\infty.$ By the definition of $\overline{\tau},$
\begin{eqnarray*}
&&\overline{u}_{\overline{\tau}}(\xi,t)\le u(\xi,t), ~\text{for all}~ (\xi,t)\in\Sigma_{x_0,\overline{\tau}},\\
&&\overline{v}_{\overline{\tau}}(\eta,z)\le v(\eta,z), ~\text{for all}~ (\eta,z)\in\Sigma_{x_0,\overline{\tau}}.
\end{eqnarray*} Since $\frac1{\theta+1}+\frac1{\kappa+1}>\frac{\alpha+\beta+\lambda}{n+1},$ we have at least one of $\mu_1$ and $\mu_2$ is positive. From \eqref{tow formula-3} and \eqref{tow formula-4} with $(x,0)=(x_0,0),\tau=\overline{\tau},$ it holds
\begin{eqnarray*}
&&\overline{u}_{\overline{\tau}}(\xi,t)<u(\xi,t), ~\text{for all}~ (\xi,t)\in\Sigma_{x_0,\overline{\tau}},\\
&&\overline{v}_{\overline{\tau}}(\eta,z)<v(\eta,z),~\text{for all}~ (\eta,z)\in\Sigma_{x_0,\overline{\tau}}.
\end{eqnarray*}

Similar to proof process of Lemma \ref{Equv-Lemma}, for $\tau\in[\overline{\tau},\overline{\tau}+\varepsilon),$ we can conclude that $\Sigma^u_{x_0,\tau}$ and $\Sigma^v_{x_0,\tau}$ must have measure zero sets. Thus we obtain that
\begin{eqnarray*}
&&\overline{u}_\tau(\xi,t)\le u(\xi,t),~ \forall (\xi,t)\in\Sigma_{x_0,\tau},\\
&&\overline{v}_\tau(\eta,z)\le v(\eta,z),~ \forall  (\eta,z)\in\Sigma_{x_0,\tau}
\end{eqnarray*} for $\tau\in[\overline{\tau},\overline{\tau}+\varepsilon),$ which contradicts the definition of $\overline{\tau}.$
\end{proof}

We need the following two calculus key lemmas  to carry out the moving sphere process. The first lemma was proved in \cite{L2004}, the second lemma was introduced by Dou and Zhu \cite{DZ2015a}.
\begin{lemma}(Lemma 5.8 in \cite{L2004})\label{f-s} Let $n\ge1$,and $\mu\in\mathbb{R},$ and $f\in C^0(\mathbb{R}^{n}).$ Suppose that for every $x\in\mathbb{R}^{n},$ there exists $\tau>0$ such that
\begin{equation*}
\big(\frac\tau{|y-x|}\big)^\mu f\big(x+\frac{\tau^2(y-x)}{|y-x|^2}\big)= f(y), ~\forall y\in\mathbb{R}^{n}\setminus\{x\}.
\end{equation*}
Then there exist  $a\ge0,d>0,$ and ${x_0}\in\mathbb{R}^{n},$ such that
\begin{equation*}
f(x)\equiv\pm a\big(\frac 1{d+|x-{x_0}|^2}\big)^{\frac\mu2}.
\end{equation*}
\end{lemma}
\begin{lemma}(Lemma 3.7 in \cite{DZ2015a})\label{f-s-1} For $n\ge0,$ and $\mu\in\mathbb{R},$ if $f\in C^0(\overline{\mathbb{R}^{n+1}_+})$ satisfying
\[\big(\frac\tau{|Y-X|}\big)^\mu f\big(X+\frac{\tau^2(Y-X)}{|Y-X|^2}\big)\le f(Y)\] for $\forall\tau>0, \, X=(x,t),Y=(y,z)\in\mathbb{R}^{n+1}_+,$ and $|Y-X|\ge\tau,$ then
\[f(X)=f(0,t), ~\forall X\in\mathbb{R}^{n+1}_+.\]
\end{lemma}

{\bf Proof of Theorem \ref{Solution-theo}.}

{ \bf Case 1.} We firstly show that if there exists some
$(x_0,0)\in\partial\mathbb{R}^{n+1}_+$ such that $\overline{\tau}(x_0,0)<\infty,$ then
$\overline{\tau}(x,0)<\infty$ for all $(x,0)\in\partial\mathbb{R}^{n+1}_+.$

Indeed, for any $(x,0)\in\partial\mathbb{R}^{n+1}_+,$ from the definition of $\overline{\tau}(x,0)$, we know $\forall\tau\in(0,\bar{\tau}(x,0)),$
\[\overline{u}_\tau(\xi,t)\le u(\xi,t),\quad\forall(\xi,t)\in\Sigma_{x,\tau}.\]
It implies that for $\forall\tau\in(0,\overline{\tau}(x,0)),$
\begin{equation}\label{a:}
a:=\liminf_{|(\xi,t)|\rightarrow\infty}(|(\xi,t)|^{\lambda+2\alpha}u(\xi,t))\ge\liminf_{|(\xi,t)|\rightarrow\infty}
(|(\xi,t)|^{\lambda+2\alpha}\overline{u}_\tau(\xi,t))=\tau^{\lambda+2\alpha}u(x,0).
\end{equation} On the other hand, since $\bar{\tau}(x_0,0)<\infty$, using Lemma \ref{Equv-Lemma}, we have
\begin{eqnarray}\label{a}
a&=&\liminf_{|(\xi,t)|\rightarrow\infty}(|(\xi,t)|^{\lambda+2\alpha}u(\xi,t))
=\liminf_{|(\xi,t)|\rightarrow\infty}(|(\xi,t)|^{\lambda+2\alpha}\overline{u}_{\overline{\tau}}
(\xi,t))\nonumber\\
&=&\overline{\tau}^{\lambda+2\alpha}u(x_0,0)\nonumber\\
&<&\infty.
\end{eqnarray} Combining \eqref{a:} with \eqref{a}, we get $\overline{\tau}(x,0)<\infty$ for all $(x,0)\in\partial\mathbb{R}^{n+1}_+.$

Invoking Lemma \ref{Equv-Lemma} again, we know
\begin{gather*}
\overline{u}_{\overline{\tau}}(\xi,t)=u(\xi,t), ~(\xi,t)\in\mathbb{R}^{n+1}_+,\\
\overline{v}_{\overline{\tau}}(\eta,z)=v(\eta,z), ~(\eta,z)\in\mathbb{R}^{n+1}_+,
\end{gather*} and $$\kappa=\kappa^*=\frac{2n+2-\lambda-2\beta}{\lambda+2\beta}\quad \text{and}\quad \theta=\theta^*=\frac{2n+2-\lambda-2\alpha}{\lambda+2\alpha}.$$
Furthermore, according to Theorem \ref{Regularity-fg} (or Theorem  \ref{Regularity}), we  know $u,v\in C^{0, \gamma}_{loc}(\overline{\mathbb{R}^{n+1}_+})$. It follows from Lemma \ref{f-s} that
\[u(\xi,0)=c_1\big(\frac 1{|\xi-\xi_0|^2+d^2}\big)^{\frac{\lambda+2\alpha}2}\quad\text{and}\quad v(\xi,0)=c_2\big(\frac 1{|\xi-\xi_0|^2+d^2}\big)^{\frac{\lambda+2\beta}2},\]
for some $c_1,c_2,d>0$, and $\xi,\xi_0\in\mathbb{R}^{n}.$

{\bf Case 2.} $\overline{\tau}(x,0)=\infty$ for all $(x,0)\in\partial\mathbb{R}^{n+1}_+$. By the definition of $\overline{\tau}(x,0)$ and Lemma \ref{Equv not-Lemma} (for case
$\kappa\le\frac{2n+2-\lambda-2\beta}{\lambda+2\beta}$ and $\theta\le\frac{2n+2-\lambda-2\alpha}{\lambda+2\alpha}$ with $\frac1{\theta+1}+\frac1{\kappa+1}>\frac{\alpha+\beta+\lambda}{n+1}$), it holds
\begin{gather*}
\overline{u}_\tau(\xi,t)\le u(\xi,t), \quad\forall(\xi,t)\in\Sigma_{x,\tau},\\
\overline{v}_\tau(\eta,z)\le v(\eta,z), \quad\forall(\eta,z)\in\Sigma_{x,\tau}
\end{gather*} for all $\tau$ and $(x,0)\in\partial\mathbb{R}^{n+1}_+.$
%\Red{Combining the above  with \eqref{tow formula-1} and  \eqref{tow formula-2},  we have
% $\kappa\le\frac{2n+2-\lambda-2\beta}{\lambda+2\beta}$ and $\theta\le\frac{2n+2-\lambda-2\alpha}{\lambda+2\alpha}$ satisfying $\frac1{\theta+1}+\frac1{\kappa+1}>\frac{\alpha+\beta+\lambda}{n+1}$.}
From Lemma \ref{f-s-1}, we know that
$u(X)=u(0,t)$ and $v(Y)=v(0,z)$.
Thus, we have
\begin{eqnarray*}
v(Y)&=&v(0,z)\\
&=&\int_{\mathbb{R}^{n+1}_+}\frac{u^\theta(0,t)}{t^\alpha z^\beta(|x|^2+|t-z|^2)^{\frac\lambda2}}dX\\
&=&Cz^{-\beta}\int^\infty_0 u^\theta(0,t)\cdot t^{-\alpha}|t-z|^{n-\lambda}\int^\infty_0\frac{\rho^{n-1}}{(\rho^2+1)^{\frac\lambda2}}d\rho dt\\
&=&Cz^{-\beta}\int^\infty_0 u^\theta(0,zs)\cdot(z s)^{-\alpha}|z s-z|^{n-\lambda}\int^\infty_0\frac{\rho^{\frac{n-2}2}}{(\rho+1)
^{\frac\lambda2}}d\rho d(z s)\\
%&=&Cz^{n+1-\alpha-\beta-\lambda}\int_0^\infty u^\theta(0,zs)\cdot s^{-\alpha}|s-1|^{n-\lambda} ds\\
&=&C z^{n+1-\alpha-\beta-\lambda}\int_0^\infty u^\theta(0,zs)\cdot s^{-\alpha}|s-1|^{n-\lambda}ds,
\end{eqnarray*}
and
\begin{eqnarray*}
u(X)&=&u(0,t)=\int_{\mathbb{R}^{n+1}_+}\frac{v^\kappa(0,z)}{t^\alpha z^\beta(|y|^2+|t-z|^2)^{\frac \lambda2}}dY\\
%&=&Ct^{-\alpha}\int_0^\infty v^\kappa(0,z)\cdot z^{-\beta}|t-z|^{n-\lambda}\int^\infty_0\frac{\rho^{n-1}}{(\rho^2+1)^{\frac\lambda2}}d\rho dz\\
%&=&Ct^{n+1-\alpha-\beta-\lambda}\int^\infty_0 v^\kappa(0,ts)\cdot s^{-\beta}|s-1|^{n-\lambda}\int^\infty_0\frac{\rho^{\frac{n-2}2}}{(\rho+1)
%^{\frac\lambda2}}d\rho ds\\
&=&Ct^{n+1-\alpha-\beta-\lambda}\int^\infty_0 v^\kappa(0,ts)\cdot s^{-\beta}|s-1|^{n-\lambda}ds.
\end{eqnarray*}
Since  $n+1-\alpha-\beta-\lambda>0,$ we conclude that $u(0,t),v(0,z)$ are infinite. This is impossible.
The theorem is proved.
\qed

On the same way, we can  show partial classification results of the positive solutions to \eqref{Euler-syst-5} as follows. It includes the result $(ii)$ in Theorem \ref{classfy-2}.
\begin{theorem}\label{Solution-theo-nm} Let $n,m\ge0, $ $0<\lambda<n+m,\,\kappa_1,\,\theta_1>0,\,\kappa_1\theta_1\ge1,$ $ \alpha+\beta\ge0,\, \alpha<\frac{m}{\theta_1+1},\, \beta<\frac {m}{\kappa_1+1}$. Suppose that
 and $(u,v)\in C^0(\mathbb{R}^{n+m})\times C^0(\mathbb{R}^{n+m})$ is a pair of positive solutions to \eqref{Euler-syst-5},  then the following two results hold.

$(i)$ If $\kappa_1=\kappa_1^*=\frac{2(n+m)-\lambda-2\beta}{\lambda+2\beta}$ and $\theta_1=\theta_1^*=\frac{2(n+m)-\lambda-2\alpha}{\lambda+2\alpha},$  then  there exist  $c_3, c_4, a>0$ and $x, x_0 \in\mathbb{R}^{n}, \hat{x}\in\mathbb{R}^{m}$ such that $u,v$ must be the  following forms on $\mathbb{R}^n\times\{\hat{x}=0\}$
\begin{gather*}
u(x,0)=c_3\big(\frac 1{|x-x_0|^2+a^2}\big)^{\frac{\lambda+2\alpha}2},\quad
v(x,0)=c_4\big(\frac 1{|x-x_0|^2+a^2}\big)^{\frac {\lambda+2\beta}2}.
\end{gather*}

$(ii)$ If $\kappa_1\le\frac{2(n+m)-\lambda-2\beta}{\lambda+2\beta}$ and $\theta_1\le\frac{2(n+m)-\lambda-2\alpha}{\lambda+2\alpha}$ satisfying $\frac1{\theta_1+1}+\frac1{\kappa_1+1}>\frac{\alpha+\beta+\lambda}{n+m},$  then $u,v=0.$
\end{theorem}

%------------------------------------------------------------------
\section{\textbf{Regularity of solutions to integral system} \label{Section 6}}
In this section, we discuss the regularity of solutions to integral system \eqref{Euler-syst-1}.

Set $U(X)=t^\alpha u(X),\, V(Y)=z^\beta v(Y),$ we rewrite \eqref{Euler-syst-1} as
\begin{equation}\label{Euler-syst-1-1}
\begin{cases}
U(X)=\int_{\mathbb{R}^{n+1}_+} \frac{v^{\kappa}(Y)}{z^\beta|X-Y|^{\lambda}}dY,&\quad X\in\mathbb{R}^{n+1}_+,\\
V(Y)=\int_{\mathbb{R}^{n+1}_+} \frac{u^
{\theta}(X)}{t^\alpha|X-Y|^{\lambda}}dX,&\quad Y \in\mathbb{R}^{n+1}_+,
\end{cases}
\end{equation} where $\theta$ and $\kappa$ satisfy \eqref{WH-exp-2}.

Theorem \ref{Regularity-fg} can be stated as follows.
\begin{prop}\label{Regularity} Let $\alpha,\beta,\lambda,\theta,\kappa$ satisfy \eqref{WH-exp-2}. If $(u,v)$ is a pair of nonnegative solutions of \eqref{Euler-syst-1-1} with $u\in L^{\theta+1}(\mathbb{R}^{n+1}_+),$
then

$(i)$ $U(X),V(Y)\in L^{\infty}(\mathbb{R}^{n+1}_+).$

$(ii)$ $U, V\in C^\infty(\mathbb{R}^{n+1}_+\setminus\{t=0\})\cup C^{0, \gamma}_{loc}(\overline{\mathbb{R}^{n+1}_+})$ for any $0< \gamma<1$.
 In particular, $u,v\in C^\infty(\mathbb{R}^{n+1}_+\setminus\{t=0\})\cup C^{0, \gamma}_{loc}(\overline{\mathbb{R}^{n+1}_+})$.
\end{prop}%Let $\alpha,\beta\ge0$ and $\alpha,\beta,\lambda,\theta,\kappa$ satisfy \eqref{WH-exp-2}. If $(u(X),v(Y))$
%are the solution of \eqref{Euler-syst-1-1} with $u\in L^{\theta+1}(\mathbb{R}^{n+1}_+)$, then $U(X),V(Y)\in L^{\infty}(\mathbb{R}^{n+1}_+)$. Moreover, $u,v$ be $C^\infty$ in $\mathbb{R}^{n+1}_+\setminus\{0\}.$
We firstly state the following regularity lift lemma introduced  by Chen and Li in \cite{CL2010}.

Assume that $Z$ is a given vector space, $\|\cdot\|_H$ and $\|\cdot\|_K$ are two norms on $Z.$ Denote $\|\cdot\|_Z$ as
\[\|\cdot\|_Z=\sqrt[q]{\|\cdot\|^q_H+\|\cdot\|^q_K},\]
where $1\le p\le\infty.$ For simplicity, we assume that $Z$ is complete with respect to the norm $\|\cdot\|_Z.$ Let $H$ and $K$ be the completion of $Z$ under $\|\cdot\|_H$ and $\|\cdot\|_K$, respectively. Obviously, $Z=H\cap K.$

\begin{lemma}(\textbf{Regularity lift lemma}, in \cite{CL2010})\label{Reli} Let $T$ be a contraction map from $H$ into itself and from $K$ into itself. Assume that $f\in H$ and that there exists a function $g\in K$ such that $f=Tf+g$, then $f$ also belongs to Z.
\end{lemma}

By Lemma \ref{Reli}, we  have
\begin{lemma}\label{Reg-lm-s} Assume $\alpha,\beta,\lambda,\theta,\kappa$ satisfy \eqref{WH-exp-2}. If $(u,v)$ is a pair of nonnegative solutions of \eqref{Euler-syst-1-1} with $u\in L^{\theta+1}(\mathbb{R}^{n+1}_+),$ then $(u,v)\in L^l(\mathbb{R}^{n+1}_+)\times L^s(\mathbb{R}^{n+1}_+),$ for any $l,s>\frac{n+1}{\lambda+\alpha+\beta}.$
\end{lemma}
\begin{proof} For any positive real number $a$, we define
\begin{equation*}
u_a(X)=\begin{cases}
u(X),&\quad \text{if}~ |u(X)|>a\ or\ |X|>a,\\
0,&\quad otherwise.
\end{cases}
\end{equation*} and $u_b(X)=u(X)-u_a(X)$.
Choose parameter $r_1,r_2>0$ and satisfy  $0<\frac1\kappa\le r_1\le\theta$  and and  $0<\frac1\theta\le r_2\le\kappa$, in particular,  we can  choose $r_1=1$ if $\kappa\ge1$ and $r_2=1$ if $\theta\ge1$.
Assume that $\varphi\in L^s(\mathbb{R}^{n+1}_+)$ and $\psi\in L^l(\mathbb{R}^{n+1}_+)$, respectively.  Define the linear operators $F_1,F_1$ and $T_1,T_2$ as
\begin{eqnarray*}
(F_1\psi)(Y)&=&\int_{\mathbb{R}^{n+1}_+} \frac{u_a^{\theta-r_1}(X)\psi^{r_1}(X)}{t^{\alpha}|X-Y|^\lambda z^{\beta}}dX,\quad\
(F_2\varphi)(X)=\int_{\mathbb{R}^{n+1}_+} \frac{v_a^{\kappa-r_2}(Y)\varphi^{r_2}(Y)}{t^{\alpha}|X-Y|^\lambda z^{\beta}}dY,\\
(T_1\psi)(X)&=&\int_{\mathbb{R}^{n+1}_+} \frac{v_a^{\kappa-\frac1{r_1}}(Y)}{t^{\alpha}|X-Y|^\lambda z^{\beta}}\big[(F_1\psi)(Y)\big]^\frac1{r_1}dY,\\
(T_2 \varphi)(Y)&=&\int_{\mathbb{R}^{n+1}_+} \frac{u_a^{\theta-\frac1{r_2}}(X)}{t^{\alpha}|X-Y|^\lambda z^{\beta}}\big[(F_2\varphi)(X)\big]^\frac1{r_2}dX.
\end{eqnarray*}And write
\begin{eqnarray*}
G(X)&=&\int_{\mathbb{R}^{n+1}_+} \frac{v_b^{\kappa-\frac1{r_1}}v^\frac1{r_1}}{t^{\alpha}|X-Y|^\lambda z^{\beta}}dY \quad\text{and}\quad
H(Y)=\int_{\mathbb{R}^{n+1}_+} \frac{u_b^{\theta-\frac1{r_2}}u^\frac1{r_2}}{t^{\alpha}|X-Y|^\lambda z^{\beta}}dX.
\end{eqnarray*} Denote the norm in the cross product $L^l(\mathbb{R}^{n+1}_+)\times L^s(\mathbb{R}^{n+1}_+)$ by
\[\|(\psi,\varphi)\|_{L^l\times L^s}=\|\psi\|_{L^l}+\|\varphi\|_{L^s},\]
 and define the operator $T :L^l(\mathbb{R}^{n+1}_+)\times L^s(\mathbb{R}^{n+1}_+)\rightarrow L^l(\mathbb{R}^{n+1}_+)\times L^s(\mathbb{R}^{n+1}_+)$ by
 \[T(\psi,\varphi)=(T_1\psi,T_2\varphi).\]

Now consider
\[
(u,v)=T(\psi,\varphi)+(G,H).\]
It is easy to see that $(u,v)$ is a pair of   positive solutions of \eqref{Euler-syst-1-1}.

Let $l,s>\frac{n+1}{\lambda+\alpha+\beta}$ and $d_1,d>1$ satisfy $\frac1d=\frac1{d_1}-\frac{n+1-(\alpha+\beta+\lambda)}{n+1},$ and $ \frac{r_1}{l}+\frac{\theta-r_1}{\theta+1}=\frac1{d_1}.$
Using Stein-Weiss type inequality \eqref{WHLSD-2} and H\"{o}lder's inequality, we have
\begin{equation}\label{Fa1}
\|F_1\psi\|_{L^{d}(\mathbb{R}^{n+1}_+)}\le C\|u_a^{\theta-r_1}\psi^{r_1}\|_{L^{d_1}(\mathbb{R}^{n+1}_+)} \le C\|u_a\|^{\theta-r_1}_{L^{\theta+1}(\mathbb{R}^{n+1}_+)}\|\psi\|^{r_1}_{L^{l}(\mathbb{R}^{n+1}_+)}
\end{equation}
Furthermore, let $l_1>1$ satisfy $\frac1l=\frac1{l_1}-\frac{n+1-(\alpha+\beta+\lambda)}{n+1},\,  \frac1{d}+\frac{\kappa r_1-1}{\kappa+1}=\frac{r_1}{l_1}.$ Using Stein-Weiss type inequality \eqref{WHLSD-2} and H\"{o}lder's inequality again, we have
\begin{equation}\label{Fa2}
\|T_1\psi\|_{L^l(\mathbb{R}^{n+1}_+)}\le C\|v^{\kappa-\frac1{r_1}}_a(F_1\psi)^{\frac1{r_1}}\|_{L^{l_1}(\mathbb{R}^{n+1}_+)}\le C\|v_a\|^{\kappa-\frac1{r_1}}_{L^{\kappa+1}(\mathbb{R}^{n+1}_+)}\|F_1\psi\|^{\frac1{r_1}}_{L^{d}(\mathbb{R}^{n+1}_+)}.
\end{equation}
Submitting \eqref{Fa1} into \eqref{Fa2} yields
\begin{equation}\label{T1p1}
\|T_1\psi\|_{L^l(\mathbb{R}^{n+1}_+)}\le C\|v_a\|^{\kappa-\frac1{r_1}}_{L^{\kappa+1}(\mathbb{R}^{n+1}_+)}\|u_a\|^{\frac{\theta-r_1}{r_1}}_{L^{\theta+1}(\mathbb{R}^{n+1}_+)}\|\psi\|_{L^{l}(\mathbb{R}^{n+1}_+)}.
\end{equation}
It is easy to check that $d$ and $l$ satisfy $\frac1d-\frac{r_1}{l}=\frac1{\kappa+1}-\frac{r_1}{\theta+1}.$

On the same way, it also holds
\begin{eqnarray}\label{T2f}
\|T_2\varphi\|_{L^s(\mathbb{R}^{n+1}_+)}
%&\le& C\|u^{\theta-\frac1{r_2}}_a(F_1\varphi)^{\frac1{r_2}}\|_{L^{s_1}(\mathbb{R}^{n+1}_+)}\le C\|u_a\|^{\theta-\frac1{r_2}}_{L^{\theta+1}(\mathbb{R}^{n+1}_+)}\|F_1\varphi\|^{\frac1{r_2}}_{L^{t}(\mathbb{R}^{n+1}_+)}\nonumber\\
&\le&C\|v_a\|^{\frac{\kappa-r_2}{r_2}}_{L^{\kappa+1}(\mathbb{R}^{n+1}_+)}\|u_a\|^{\theta-\frac1{r_2}}_{L^{\theta+1}(\mathbb{R}^{n+1}_+)}\|\varphi\|_{L^{s}(\mathbb{R}^{n+1}_+)}.
\end{eqnarray}
From \eqref{Re-v}, we know $ v\in L^{\kappa+1}(\mathbb{R}^{n+1}_+).$ It implies that $u_a\in L^{\theta+1}(\mathbb{R}^{n+1}_+),\,  v_a\in L^{\kappa+1}(\mathbb{R}^{n+1}_+).$  We  may choose a sufficiently large $a$ such that
\begin{equation}\label{uva}
C\|v_a\|^{\kappa-\frac1{r_1}}_{L^{\kappa+1}(\mathbb{R}^{n+1}_+)}\|u_a\|^{\frac{\theta-r_1}{r_1}}
_{L^{\theta+1}(\mathbb{R}^{n+1}_+)}\le\frac12,\,C\|v_a\|^{\frac{\kappa-r_2}{r_2}}_{L^{\kappa+1}(\mathbb{R}^{n+1}_+)}\|u_a\|^{\theta-\frac1{r_2}}_{L^{\theta+1}(\mathbb{R}^{n+1}_+)}\le\frac12.
\end{equation}
Substituting \eqref{uva} into \eqref{T1p1} and \eqref{T2f}, we arrive at
\begin{equation}\label{uva1}
\|T_1\psi\|_{L^l(\mathbb{R}^{n+1}_+)}\le\frac12\|\psi\|_{L^l(\mathbb{R}^{n+1}_+)},
\quad \|T_2\varphi\|_{L^s(\mathbb{R}^{n+1}_+)}\le\frac12\|\varphi\|_{L^s(\mathbb{R}^{n+1}_+)}.
\end{equation}
According to the definition of the norm in the cross product $L^l(\mathbb{R}^{n+1}_+)\times L^s(\mathbb{R}^{n+1}_+)$, we have
\begin{equation*}
\|T(\psi,\varphi)\|_{L^l(\mathbb{R}^{n+1}_+)\times L^s(\mathbb{R}^{n+1}_+)}\le\frac12\|\psi\|_{L^l(\mathbb{R}^{n+1}_+)}+\frac12\|\varphi\|_{L^s(\mathbb{R}^{n+1}_+)}.
\end{equation*}Thus we can obtain that $T$ is a contraction map from $L^l(\mathbb{R}^{n+1}_+)\times L^s(\mathbb{R}^{n+1}_+)$ to itself.

Next, we will show that $(G,H)\in L^l(\mathbb{R}^{n+1}_+)\times L^s(\mathbb{R}^{n+1}_+)$
for any $l,s>\frac{n+1}{\lambda+\alpha+\beta}.$ Involving Stein-Weiss type inequality \eqref{WHLSD-2} and H\"{o}lder's inequality, we have
\[\|G(x,t)\|_{L^l(\mathbb{R}^{n+1}_+)}\le C\|v^{\kappa-\frac1{r_1}}_b v^\frac1{r_1}\|_{L^{l_1}(\mathbb{R}^{n+1}_+)}\le C\|v_b\|^{\kappa-\frac1{r_1}}_{{L^{k_1}(\mathbb{R}^{n+1}_+)}}\|v\|^\frac1{r_1}_{{L^{k_2}(\mathbb{R}^{n+1}_+)}},\]
where $1<k_1,k_2<\infty$ with $\frac{\kappa-\frac1{r_1}}{k_1}+\frac1{k_2r_1}=\frac1l+\frac{n+1-(\alpha+\beta+\lambda)}{n+1}.$

If $k_2>\frac{n+1}{(n+1-(\alpha+\beta+\lambda))r_1}$ for any $l>\frac{n+1}{\lambda+\alpha+\beta},$
 \[\|G(x,t)\|_{L^l(\mathbb{R}^{n+1}_+)}<\infty.\]
Similarly, we have
\[\|H(y,z)\|_{L^s(\mathbb{R}^{n+1}_+)}<\infty ~~ \text{for any} ~s>\frac{n+1}{\lambda+\alpha+\beta}.\]
By Lemma \ref{Reli}, we can derive that $(u,v)\in L^l(\mathbb{R}^{n+1}_+)\times L^s(\mathbb{R}^{n+1}_+)$ for any $l,s>\frac{n+1}{\lambda+\alpha+\beta}.$

If $k_2\le\frac{n+1}{(n+1-(\alpha+\beta+\lambda))r_1},$ for any  $\frac{n+1}{\lambda+\alpha+\beta}<l\le\frac{k_2r_1(n+1)}{(n+1)r_2-(n+1-(\alpha+\beta+\lambda))k_2r_1},$ we also have
\[\|G(x,t)\|_{L^l(\mathbb{R}^{n+1}_+)}<\infty, \quad\|H(y,z)\|_{L^s(\mathbb{R}^{n+1}_+)}<\infty.\]
By Lemma \ref{Reli} again, we arrive at $(u,v)\in L^l(\mathbb{R}^{n+1}_+)\times L^s(\mathbb{R}^{n+1}_+)$ for all $\frac{n+1}{\alpha+\beta+\lambda}<l,s\le\frac{k_2r_1(n+1)}{(n+1)r_1-(n+1-(\alpha+\beta+\lambda))\kappa_2r_1}.$

Now, we begin to iterate. That is, choosing $k_2=k^1_2$ with $k_2^1=\frac{\kappa_2r_1(n+1)}{(n+1)r_1-(n+1-(\alpha+\beta+\lambda))k_2r_1},$  and arguing as the above, we can obtain that $(u,v)\in L^l(\mathbb{R}^{n+1}_+)\times L^s(\mathbb{R}^{n+1}_+)$ for any $\frac{n+1}{\lambda+\alpha+\beta}<l,s\le\frac{k^1_2(n+1)}{n+1-(n+1-(\alpha+\beta+\lambda))k^1_2}.$
Repeating the above process after few steps, there exists some $k_2^i$ such that $k_2^i>\frac{n+1}{(n+1-(\alpha+\beta+\lambda))r_1}.$ It contains in the case of $k_2>\frac{n+1}{(n+1-(\alpha+\beta+\lambda))r_1}.$ Therefore, we conclude that $(u,v)\in L^l(\mathbb{R}^{n+1}_+)\times L^s(\mathbb{R}^{n+1}_+)$ for any $l,s>\frac{n+1}{\lambda+\alpha+\beta}.$
\end{proof}
\begin{lemma}\label{Reg-lm-infty}
  Under the assumptions in Proposition \ref{Regularity}, it holds $U, V\in L^\infty (\mathbb{R}^{n+1}_+).$
\end{lemma}
\begin{proof} Let $a$ be some fixed  positive  constant, and any $X\in\mathbb{R}^{n+1}_+$. It follows from  \eqref{Euler-syst-1-1} that
\begin{eqnarray}\label{UES}\small
|U(X)|&\le&\int_{ B_a(X)\cap\mathbb{R}^{n+1}_+}\frac{|v(Y)|^{\kappa}}{z^{\beta}|X-Y|^{\lambda}}dY +\int_{\mathbb{R}^{n+1}_+\setminus B_a(X)}\frac{|v(Y)|^{\kappa}}{z^{\beta}|X-Y|^{\lambda}}dY\nonumber \\
&=:&U_1(X)+U_2(X).
\end{eqnarray}

We first estimate $U_1(X).$
From \eqref{WH-exp-2} we know
\[
\frac{n+1}{\kappa+1}+\frac{n+1}{\theta+1}=\alpha+\beta+\lambda.
\]
This implies
\[ n+1-(\beta+\lambda)(\kappa+1)=(\kappa+1)\big(\alpha-\frac{n+1}{\theta+1}\big)<- \frac{n(\kappa+1) }{\theta+1}<0\]
 for $\alpha<\frac{1}{\theta+1}$.  That is, $\frac{n+1}{\lambda+\beta}<\kappa+1$.

Choose $1<s_1<\frac{n+1}{\lambda+\beta}$. Obviously, $n+1-(\lambda+\beta) s_1>0$ and $s_1<\kappa+1, s'_1>\frac{\kappa+1}\kappa$. For any $X\in\mathbb{R}^{n+1}_+$,  by H\"{o}lder's inequality and Lemma \ref{Reg-lm-s},  we have
\begin{eqnarray}\label{U2ES}
U_1(X)&=&\int_{B_a(X)\cap\mathbb{R}^{n+1}_+}\frac{|v(Y)|^{\kappa}}{z^{\beta}|X-Y|^{\lambda}}dY\nonumber\\
&\le&\|v\|^\kappa_{L^{s'_1\kappa}(B^+_a(0))}\big(\int_{B_a(X)\cap\mathbb{R}^{n+1}_+}\big(\frac1{z^\beta |X-Y|^\lambda}\big)^{s_1}dY\big)^{\frac1{s_1}}\nonumber\\
%&\le&\|v\|^\kappa_{L^{s'_2}(\mathbb{R}^{n+1}_+)}\big(\int_{B_a(X)\cap\{z\le b\}\cap\mathbb{R}^{n+1}_+}\big(\frac1{z^\beta |X-Y|^\lambda}\big)^{s_2}dY\big)^{\frac1{s_2}}\nonumber\\
&\le&C(\int_{B_a(X)\cap\mathbb{R}^{n+1}_+}\big(\frac1{z^\beta |X-Y|^\lambda}\big)^{ s_1}dY)^{\frac1{s_1}}.
\end{eqnarray}
For the sake of simplicity,  write
\begin{eqnarray*}
 A_1&=&\{|Y|\le|X-Y|<a\}\cap\mathbb{R}^{n+1}_+, \\
 A_2&=&\{z<|X-Y|<|Y|\}\cap\{|X-Y|<a\} \cap\mathbb{R}^{n+1}_+, \\
 A_3&=&\{|X-Y|\le z\}\cap\{|X-Y|<a\}\cap\mathbb{R}^{n+1}_+.
\end{eqnarray*}
Obviously, $B_a(X)\cap\mathbb{R}^{n+1}_+=A_1\cup A_2\cup A_3$.
Note that $|X|< 2|Y|$ for $|X|-|Y|\le|X-Y|<|Y|$.
For $\beta\ge0,$  we have
\begin{eqnarray}\label{AES1}
& &\int_{B_a(X)\cap\mathbb{R}^{n+1}_+}\big(\frac1{z^\beta |X-Y|^\lambda}\big)^{s_1}dY\nonumber\\
&=&\int_{A_1}\big(\frac1{z^\beta |X-Y|^\lambda}\big)^{s_1}dY
+\int_{A_2}\big(\frac1{z^\beta |X-Y|^\lambda}\big)^{s_1}dY
+\int_{A_3}\big(\frac1{z^\beta |X-Y|^\lambda}\big)^{s_1}dY\nonumber\\
&\le&\int_{A_1}\frac1{z^{\beta s_1}|Y|^{\lambda s_1}}dY
+\int_{A_2}\big(\frac1{z^\beta |X-Y|^\lambda}\big)^{s_1}dY
+\int_{A_3}\frac1{|X-Y|^{(\lambda+\beta)s_1}}dY\nonumber\\
&\le&\int_{B_{a}(0)}\frac1{z^{\beta s_1}|Y|^{\lambda s_1}}dY
+\int_{A_2}\big(\frac1{z^\beta |Y|^\lambda|\frac{|X|}{|Y|}-1|^\lambda}\big)^{s_1}dY
+\int_{B_a(X)}\frac1{|X-Y|^{(\lambda+\beta)s_1}}dY\nonumber\\
&\le&J_{\beta s_1} \int_{0}^{a}\rho^{n-(\lambda+\beta)s_1}d\rho
+cJ_{\beta s_1} \int_{0}^{a}\rho^{n-(\lambda+\beta)s_1}d\rho
+\omega_{n}\int_{0}^{a}\rho^{n-(\lambda+\beta)s_1}d\rho\nonumber\\
&\le&C a^{n+1-(\lambda+\beta)s_1}<\infty,
\end{eqnarray}
where $J_{\beta s_1}$ is finite due to $\beta s_1<1$ for $s_1<\kappa+1$. On the other hand,
since $|X-Y|\le z\le|Y|$,  $|X|< 2|Y|$  still holds in  the domain $A_3$. Thus,  for $-1<\beta<0$, we have
\begin{eqnarray*}
\int_{A_3}\big(\frac1{z^\beta |X-Y|^\lambda}\big)^{s_1}dY
&\le&\int_{A_3}\big(\frac1{|Y|^\beta |Y|^\lambda|\frac{|X|}{|Y|}-1|^\lambda}\big)^{s_1}dY\nonumber\\
&\le&C \omega_{n}\int_{0}^{a}\rho^{n-(\lambda+\beta)s_1}d\rho\nonumber\\
&\le&C a^{n+1-(\lambda+\beta)s_1}.
\end{eqnarray*}
For $-1<\beta<0$, similar to \eqref{AES1} and combining the above, it still yields
\begin{eqnarray}\label{AES2}
 \int_{B_a(X)\cap\mathbb{R}^{n+1}_+}\big(\frac1{z^\beta |X-Y|^\lambda}\big)^{s_1}dY
&\le&C a^{n+1-(\lambda+\beta)s_1}<\infty,
\end{eqnarray}
Combining \eqref{AES1} and \eqref{AES2} into \eqref{U2ES}, we deduce that
\begin{equation}\label{U2ES1}
U_1(X)=\int_{B_a(X)\cap\mathbb{R}^{n+1}_+}\frac{|v(Y)|^{\kappa}}{z^{\beta}|X-Y|^{\lambda}}dY<\infty.
\end{equation} for $\beta>-1$.

Now, we  estimate $U_2(X).$
Let $a>|X|$. Using the fact $\beta<\frac1{\kappa+1}$ and $n+1-(\beta+\lambda)(\kappa+1)<0$, and by H\"{o}lder's inequality,  we have
\begin{eqnarray}\label{U3ES}
U_2(X)&\le&\big(\int_{\mathbb{R}^{n+1}_+}v^{\kappa+1}(Y)dY\big)^{\frac\kappa{\kappa+1}}
\big(\int_{\mathbb{R}^{n+1}_+\setminus B_a(X)}\big(\frac1{z^\beta|X-Y|^{\lambda}}\big)^{\kappa+1}dY\big)^{\frac1{\kappa+1}}
\nonumber\\
&\le&\|v\|^\kappa_{L^{\kappa+1}(\mathbb{R}^{n+1}_+)}\big(J_{\beta(\kappa+1)}\int_a^\infty
(\frac1{\rho^{\beta+\lambda}|1-\frac{|X|}\rho|^\lambda} )^{\kappa+1}\rho^nd\rho\big)^{\frac1{\kappa+1}}\nonumber\\
&\le&C\|v\|^\kappa_{L^{\kappa+1}(\mathbb{R}^{n+1}_+)}\big(J_{\beta(\kappa+1)}
\int_a^\infty\rho^{n-(\beta+\lambda)(\kappa+1)}d\rho\big)^{\frac1{\kappa+1}}\nonumber\\
&\le&C\|v\|^\kappa_{L^{\kappa+1}(\mathbb{R}^{n+1}_+)}a^{n+1-(\beta+\lambda)(\kappa+1)}\nonumber\\
&<&\infty.
\end{eqnarray}
Thus, from \eqref{U2ES1} and \eqref{U3ES}, we can obtain that $U\in L^\infty(\mathbb{R}^{n+1}_+)$. Similarly, we can obtain that $V\in L^\infty(\mathbb{R}^{n+1}_+)$.
\end{proof}

\medskip

{\bf Proof of Proposition \ref{Regularity}.}
For any $X\in \mathbb{R}^{n+1}_+\setminus\{t=0\}$ and some suitable constant $a>0,$ we divide  $U(X)$ into the following form
\begin{eqnarray}\label{B3a}
U(X)&=&\int_{B_{3a}(X)}\frac {v^{\kappa}(Y)}{z^\beta|X-Y|^\lambda}dY
+\int_{\mathbb{R}^{n+1}_+\setminus B_{3a}(X)}\frac {v^{\kappa}(Y)}{z^\beta|X-Y|^\lambda}dY\nonumber\\
&=:&W_1(X)+W_2(X).
\end{eqnarray}
 By standard singular integral estimate (Chapter 10 in \cite{EM2001}), it is easy to see that
\begin{equation}\label{contiB2r}
W_1(X)=\int_{B_{3a}(X)}\frac{v^{\kappa}(Y)}{z^\beta|X-Y|^\lambda}dY\in C^{\gamma_1}(\mathbb{R}^{n+1}_+\setminus\{t=0\})
\end{equation} for any $\gamma_1<n+1-\lambda.$

 Next,  we show
\begin{eqnarray*}\label{G(X)}
W_2(X)\in C^\infty(\mathbb{R}^{n+1}_+\setminus\{t=0\}).
\end{eqnarray*}
Indeed, for fixed $X$ and $h$ small enough, by the mean value theorem, it holds
\begin{eqnarray*}
\big|\frac{W_2(X+he_i)-W_2(X)}{h}\big|&=&\big|\frac{\frac{v^{\kappa}(Y)}{z^\beta}\big(\frac{1}
{|X+he_i-Y|^\lambda}-\frac{1}{|X-Y|^\lambda}\big)}{h}\big|\nonumber\\
&\le& C\frac{|v(Y)|^{\kappa}}{z^\beta|X+\sigma he_i-Y|^{\lambda+1}},
\end{eqnarray*}
where $e_i=\{0,\ldots,1,\ldots,0\}$ is the $ith$ unit vector, $0<\sigma<1$. Let $\varepsilon>0$  small,  we have
\begin{eqnarray}\label{first derivative}
&&\int_{\mathbb{R}^{n+1}_+\setminus B_{3a}(X+\sigma he_i)}\frac {v^{\kappa}(Y)}{z^\beta|X+\sigma he_i-Y|^{\lambda+1}}dY\nonumber\\
&\le&\int_{\mathbb{R}^{n+1}_+\setminus B_{2a}(X)}\frac {v^{\kappa}(Y)}{z^\beta|X+\sigma he_i-Y|^{\lambda+1}}dY\nonumber\\
&\le&\int_{\mathbb{R}^{n+1}_+\setminus B_{2a}(X)}\frac {v^{\kappa}(Y)}{z^\beta|X-Y|^{\lambda+1}}dY.
\end{eqnarray}
Similar to the estimate $U_3$ in \eqref{U3ES}, using the fact $\beta<\frac1{\kappa+1}$ and $n+1-(\beta+\lambda)(\kappa+1)<0$, and by H\"{o}lder's inequality,  it holds
\begin{eqnarray}\label{CES}
& &\int_{\mathbb{R}^{n+1}_+\setminus B_{2a}(X)}\frac {v^{\kappa }(Y)}{z^\beta|X-Y|^{\lambda+1}}dY\nonumber\\
&\le&\big(\int_{\mathbb{R}^{n+1}_+}v^{\kappa+1}(Y)dY\big)^{\frac\kappa{\kappa+1}}\big(\int_{\mathbb{R}^{n+1}_+\setminus B_{2a}(X)}(\frac1{z^\beta|X-Y|^{\lambda+1 }})^{\kappa+1}dY\big)^{\frac1{\kappa+1}}\nonumber\\
\nonumber\\
&<&\infty.
\end{eqnarray}

Substituting \eqref{CES} into \eqref{first derivative}, and by the Lebesgue dominated convergence theorem, we have $W_2(X)\in C^1(\mathbb{R}^{n+1}_+\setminus\{t=0\}).$
Furthermore, by the standard process, we have $W_2(X)\in C^\infty(\mathbb{R}^{n+1}_+\setminus\{t=0\})$.

Combining the above and  \eqref{contiB2r}, we derive that $U(X)\in C^{\gamma_1}(\mathbb{R}^{n+1}_+\setminus\{t=0\})$  with $\gamma_1<n+1-\lambda$. By the bootstrap technique, we can prove that $U(X)\in C^\infty(\mathbb{R}^{n+1}_+\setminus\{t=0\})$.

Finally, we estimate the regularity near $t=0$. For some constant $b>0$, we can write
\begin{eqnarray}\label{OEM}
U(X)-U(x,0)
&=&\int_{\mathbb{R}^{n+1}_+\cap B_b(X)}\frac {v^{\kappa}(Y)}{z^\beta}\big(\frac1{|X-Y|^\lambda}-\frac 1{|(x,0)-Y|^\lambda}\big)dY\nonumber\\
\quad&&+\int_{\mathbb{R}^{n+1}_+\setminus B_b(X)}\frac {v^{\kappa}(Y)}{z^\beta}\big(\frac1{|X-Y|^\lambda}-\frac 1{|(x,0)-Y|^\lambda}\big)dY.\nonumber\\
\end{eqnarray} Let $D\subset \mathbb{R}^{n+1}_+\cap\{t<b\}$ be a bound domain with $\{t=0\}\cap D\neq \emptyset $.  We can show that $U\in C^{0,\gamma}(D)$ for any $0<\gamma<1$.

Assume $m\ge1$ and $0<l<1$, and recall the following basic inequalities
\begin{eqnarray}
& &||x-z|^{-m}-|y-z|^{-m}|\le C|x-y|^l(|x-z|^{-m-l}+|y-z|^{-m-l}),\label{m-basic-inq-1}\\
& &|\ln|x-z|-\ln|y-z||\le |x-y|^l(|x-z|^{-l}+|y-z|^{-l})/l\label{m-basic-inq-2}
\end{eqnarray}
for any $x,y,z\in \mathbb{R}^n$ (see Chapter 10 in \cite{EM2001}).

We divide \eqref{OEM} into two following cases.

{\bfseries Case 1.} $\lambda\ge1.$ By \eqref{m-basic-inq-1}  we obtain
\begin{eqnarray*}
\sup_{X\in D}\frac{|U(X)-U(x,0)|}{t^{\gamma}}&\le&C\sup_{X\in D}\int_{\mathbb{R}^{n+1}_+\cap B_b(X)}
\frac{|v(Y)|^{\kappa}}{z^\beta|X-Y|^{\lambda+\gamma}}dY\\
\quad&&+C\sup_{X\in D}\int_{\mathbb{R}^{n+1}_+\setminus B_b(X)}\frac{|v(Y)|^{\kappa}}{ z^\beta|X-Y|^{\lambda+\gamma}}dY.
\end{eqnarray*}
Arguing as \eqref{U2ES1} and \eqref{U3ES}, we have
\begin{eqnarray*}\label{OES1}
\sup_{X\in D}\int_{\mathbb{R}^{n+1}_+\cap B_b(X)}\frac{|v(Y)|^{\kappa}}{z^\beta|X-Y|^{\lambda+\gamma}}dY
&<&\infty,\label{OES1}\\
\sup_{X\in D}\int_{\mathbb{R}^{n+1}_+\setminus B_b(X)}\frac{|v(Y)|^{\kappa}}{ z^\beta|X-Y|^{\lambda+\gamma}}dY
&<&\infty.\label{OES2}
\end{eqnarray*}
Therefore,
\begin{eqnarray}\label{OES3}
\sup_{X\in D}\frac{|U(X)-U(x,0)|}{t^{\gamma}}
&<&\infty.
\end{eqnarray}

{\bfseries Case 2.} $0<\lambda<1.$ Using inequality \eqref{m-basic-inq-2}, and arguing as \textbf{Case $1$}, we have
 \begin{eqnarray}\label{OES4}
\sup_{X\in D}\frac{|U(X)-U(x,0)|}{t^{\gamma}}&\le&C\sup_{X\in D}\int_{\mathbb{R}^{n+1}_+\cap B_b(X)}
\frac{|v(Y)|^{\kappa}}{z^\beta|X-Y|^{\gamma}}dY\nonumber\\
\quad&&+C\sup_{X\in D}\int_{\mathbb{R}^{n+1}_+\setminus B_b(X)}\frac{|v(Y)|^{\kappa}}{ z^\beta|X-Y|^{\gamma}}dY\nonumber\\
&<&\infty
\end{eqnarray} for any $0< \gamma<1.$   Therefore,   combining \eqref{OES3} and \eqref{OES4} yields $U\in C^{0,\gamma}(D)$. Since $D$ is arbitrary, we obtain $U\in C^{0,\gamma}_{loc}(\overline{\mathbb{R}^{n+1}_+}).$

Similar to the above process, we can show $V\in C^\infty(\mathbb{R}^{n+1}_+\setminus\{t=0\})\cup C^{0,\gamma}_{loc}(\overline{\mathbb{R}^{n+1}_+})$ for any $0<\gamma<1.$
The proof is complete.
\qed

%--------------------------------------------------------------
\section{\textbf{Some weighted Sobolev inequalities}}\label{section 7 }
In this section, we will apply Stein-Weiss type inequalities with partial variable weight functions  to deduce some weighted Sobolev inequalities on $\mathbb{R}^{n+1}_+$ and $\mathbb{R}^{n+m}$, respectively.

To prove Theorem \ref{Weighted-Sobolev-1}, we  firstly need the following representation formula.
\begin{lemma} \label{two formulas} Let $x\in\overline{\mathbb{R}^{n}_+}, \,  n\ge2$ and  for any $f\in C^\infty_0(\overline{\mathbb{R}^{n}_+})$, it has
\begin{equation}\label{representation-nabla-1}
f(x)=\frac{1}{C(n)}\int_{\mathbb{R}^n_+}\frac{\langle\nabla f(y),x-y\rangle}{|x-y|^{n}}dy,
\end{equation}
where $C(n)=\frac{\omega_{n-1}}2 =\frac{\pi^\frac n2}{\Gamma(\frac n2)}$, $\omega_{n-1}$  is the surface area of the unit sphere $\mathbb{S}^{n-1}.$

Furthermore, let $m=2,3\cdots,n-1$ be an integral number, $k=0,1, 2,\cdots, [\frac m2]-1$, where $[\frac m2]$ is the greatest integer part of $\frac m2.$
If $$\lim_{y_n\to0}\frac {\partial f(y)}{\partial{y_n}}=\lim_{y_n\to0}\frac {\partial \Delta f(y)}{\partial{y_n}}=\cdots=\lim_{y_n\to0}\frac {\partial \Delta^{k} f(y)}{\partial{y_n}}=0,$$
 then
 \begin{equation}\label{representation-Delta-2}
f(x)=
\begin{cases}
(-1)^{\frac{m-1}2}\frac{\Gamma(\frac{n-m+1}2)}{\pi^{\frac n2}2^{m-1}\Gamma(\frac{m+1}2)}\int_{\mathbb{R}^n_+}
\frac{\langle\nabla\Delta^{\frac{m-1}2} f(y),(x-y) \rangle}{|x-y|^{n-m+1}}dy,&\quad m=odd,\\
(-1)^\frac{m}2\frac{\Gamma(\frac{n-m}2)}{\pi^{\frac n2}2^{m-1}\Gamma(\frac{m}2)}\int_{\mathbb{R}^n_+}\frac{\Delta^{\frac{m}2} f(y)}{|x-y|^{n-m}}dy,&\quad m=even.
\end{cases}
\end{equation}
\end{lemma}
\begin{proof} We only show the formulas for $m=1,2$, the remaining formulas is easily computed by exactly the same way.

Let $z\in \partial B^+_1(x)$ and write
\begin{equation*}
f(x)=\int^\infty_0\frac{d}{dt}f(x-tz)dt=\int^\infty_0\langle\nabla f(x-tz),z\rangle dt.
\end{equation*}
Integrating both sides on $\partial B^+_1(x)$ with respect to $z$ variable, we have
\begin{eqnarray}\label{rep-f-2}
f(x)&=&\frac1{C(n)}\int_{\partial B_1^+(x)}\int^\infty_0\langle\nabla f(x-tz),z\rangle  dt dS_{z}\nonumber\\
% & & -\int_0^{\frac\pi 2}\sin^{n-2}\theta d\theta \int_{S^{n-2}}\int^\infty_0\langle\nabla_{x+tz} f(x+tz),z\rangle t^{n-1+1-n}dt dS_{z'}\\
&=&\frac1{C(n)}\int_{\mathbb{R}^n_+}\frac{\langle\nabla f(y),x-y\rangle}{|x-y|^n}dy.%\quad(z=\frac{x-y}{|x-y|},t=|x-y|)
\end{eqnarray}
Thus, \eqref{representation-nabla-1} is proved.

Since $\nabla_{y}(|x-y|^{2-n})=-(2-n)\frac{x-y}{|x-y|^n}$ and $\lim_{y_n\to0}\frac {\partial f(y)}{\partial{y_n}}=0,$ integrating by parts, it follows from \eqref{rep-f-2} that
\begin{eqnarray}\label{rep-f-3}
f(x)&=&\frac1{(n- 2)C(n)}\int_{\mathbb{R}^n_+} \langle\nabla f(y),\nabla_{y}(|x-y|^{2-n})\rangle dy\nonumber\\
&=&-\frac1{(n- 2)C(n)}\big[\int_{\mathbb{R}^n_+} \frac{\Delta f(y)}{|x-y|^{n-2}}dy
+\int_{\partial\mathbb{R}^n_+} \lim_{y_n\to0}\frac {\partial f(y)}{\partial{y_n}}\frac{1}{|x-y'|^{n-2}}dy'\big]\nonumber\\
&=&-\frac1{(n-2)C(n)}\int_{\mathbb{R}^n_+} \frac{\Delta f(y)}{|x-y|^{n-2}}dy.
\end{eqnarray} We obtain \eqref{representation-Delta-2} with $m=2.$

%Noting that $ |x-y|^{2-n}=\frac{\Delta_y(|x-y|^{4-n})}{2(4-n)},$ integrating by parts again, from \eqref{rep-f-3}, we have
%\begin{eqnarray*}
%f(x)
%&=&\frac{1}{2(n-4)(n-2)C(n)}\int_{\mathbb{R}^n_+}{\Delta f(y)\Delta(|x-y|^{4-n})}dy\nonumber\\
%&=&-\frac{1}{2(n-4)(n-2)C(n)}\big[\int_{\mathbb{R}^n_+}\langle\nabla(\Delta f(y)),\nabla_y(|x-y|^{4-n})\rangle dy\nonumber\\
%\quad&&-\int_{\partial\mathbb{R}^n_+}\lim_{y_n\to0} \frac{\partial (|x-y'|^{4-n})}{\partial{y_n}}{\Delta f(y)}dy'\big]\nonumber\\
%&=&-\frac{1}{2(n-2)C(n)}\int_{\mathbb{R}^n_+}\frac{\langle\nabla(\Delta f(y)),x-y\rangle}{|x-y|^{n-2}}dy.
%\end{eqnarray*}
% We deduce \eqref{representation-Delta-2} with $m=3.$
%
%Furthermore, due to $\lim_{y_n\to0}\frac {\partial [\Delta_{y}f(y)]}{\partial{y_n}}=0,$ an integration by part twice yields
%\begin{eqnarray*}
%f(x)&=&\frac{1}{2(n-4)(n-2)C(n)}\int_{\mathbb{R}^n_+}{\Delta f(y)\Delta(|x-y|^{4-n})}dy\nonumber\\
%&=&-\frac{1}{2(n-4)(n-2)C(n)}\int_{\mathbb{R}^n_+}\langle\nabla(\Delta f(y)),\nabla_{y}(|x-y|^{4-n})\rangle dy\nonumber\\
%&=&\frac{1}{2(n-2)(n-4)C(n)}\nonumber\\
%\quad &&\times\big[\int_{\mathbb{R}^n_+}\frac{\Delta_{y}^2 f(y)}{|x-y|^{n-4}}dy
%-\int_{\partial\mathbb{R}^n_+} \lim_{y_n\to0}\frac {\partial [\Delta_{y}f(y)]}{\partial{y_n}}\frac{1}{|x-y'|^{n-4}}dy'\big]\nonumber\\
%&=&\frac{1}{2(n-2)(n-4)C(n)}\int_{\mathbb{R}^n_+}\frac{\Delta_{y}^2 f(y)}{|x-y|^{n-4}}dy.
%\end{eqnarray*}We have \eqref{representation-Delta-2} with $m=4.$
\end{proof}

\textbf{Proof of Theorem \ref{Weighted-Sobolev-1}.}
We only show the proof of inequality \eqref{WS-1}. On the same way, we can obtain  that inequalities \eqref{WS-mn-odd} and  \eqref{WS-mn-even} hold.

For any $f\in C^\infty_0(\overline{\mathbb{R}^{n+1}_+})$ and $X\in \mathbb{R}^{n+1}_+,$ invoking \eqref{representation-nabla-1} on $\mathbb{R}^{n+1}_+$ we have
\begin{equation}\label{Tr-1}
|f(X)|\le\frac{1}{C(n+1)}\int_{\mathbb{R}^{n+1}_+}\frac{|\nabla f(Y)|}{|X-Y|^{n}}dY.
\end{equation}
Substituting \eqref{Tr-1} into Stein-Weiss type inequality \eqref{WHLSD-1} with $\lambda=n$, we arrive at
\begin{equation}\label{WS-1-a}
\big(\int_{\mathbb{R}^{n+1}_+}t^{-\beta q}|u|^{q}dX\big)^{\frac p{q}}\le S(n+1,p,\alpha,\beta)\int_{\mathbb{R}^{n+1}_+}t^{\alpha p}|\nabla u|^pdX
\end{equation}  for some positive constant $S(n+1, p, \alpha,\beta)$. Here $\alpha<\frac{1}{p'},\, \beta<\frac{1}q,\, \alpha+\beta\ge0, \, 1< p\le q<\infty,$ and $\frac1q=\frac1p-\frac{n+1-(\alpha+\beta+\lambda)}{n+1}.$

Next, we only need to calculate the exponents showed in inequality \eqref{WS-1}. Let $\alpha_1=p \alpha,\, \beta_1=-q \beta$ in \eqref{WS-1-a},  it is just the inequality \eqref{WS-1} by density. We can deduce that $\frac{\alpha_1}{p}-\frac{\beta_1}q\ge0,$ then $\alpha_1\ge\frac{p\beta_1}q,$ $\alpha_1<p-1$ and $\beta_1=-q\beta>-1.$

Since $\lambda=n$ and $\frac1q=\frac1p+\frac{\alpha+\beta+n}{n+1}-1,$ %=\frac1p+\frac{\frac{\alpha_1}p-\frac{\beta_1}q-1}{n+1},\]
we have
\begin{equation}\label{q-0}
\frac{n+1+\beta_1}q=\frac{n+1+\alpha_1}p-1,
\end{equation} or
\begin{equation}\label{q-1}
q=\frac{p(n+1+\beta_1)}{n+1+\alpha_1-p}:=p^*(\alpha_1, \beta_1).
\end{equation}
%and
%\[\frac{n+1}q+\frac{\beta_1}q=\frac{n+1}p+\frac{\alpha_1}p.\]
From \eqref{q-0} and the fact $\beta<\frac1q,$ we conclude that
\[\beta_1=\frac{(n+1+\alpha_1)q}p-(n+1)-q>-1,\]
and then
\begin{equation}\label{alpha1 fomula}
\alpha_1>\frac{pn}q-(n+1-p).
\end{equation}
Thus
\begin{equation}\label{EX01}
\frac{pn}q-(n+1-p)<\alpha_1<p-1.
\end{equation} (From \eqref{EX01}, we know that $\alpha_1$ doesn't exist for $p=q$.)

On the other hand,
since $q>p$ , it follows from \eqref{q-1}  that
\begin{equation*}
n+1+\beta_1> n+1+\alpha_1-p,
\end{equation*} then $\alpha_1<p+\beta_1.$
% In particular,  if $p^*= p$,  then  $\alpha_1=\beta_1$ from  \eqref{q-0}.

Since $\alpha_1\ge\frac{\beta_1p}q=\frac{\beta_1(n+1+\alpha_1-p)}{n+1+\beta_1}$ and $\beta_1>-1,$ we obtain
\[\alpha_1(n+1+\beta_1)\ge \beta_1(n+1+\alpha_1-p).\]
It implies
\begin{equation*}
\alpha_1\ge\frac{n+1-p}{n+1}\beta_1.
\end{equation*} Thus, we have
\begin{equation}\label{EX02}
\frac{(n+1-p)}{n+1}\beta_1\le\alpha_1< p+\beta_1.
\end{equation} Combining  the right side of \eqref{EX01} and \eqref{EX02}, we have
 \[\alpha_1<\min\{p-1, p+\beta_1\}=p-1.\]
So, we  have verified the condition of exponents as
\[ \frac{pn}q-(n+1-p)<\alpha_1<p-1, ~ \beta_1\le\alpha_1\frac{n+1}{n+1-p}.\]
\qed
%Meanwhile, we can check that
%\[p^*<\frac{p(n+1+\beta-1)}{n+1-p+\frac{(n+1-p)\beta_1}{n+1}}=\frac{p(n+1+\beta_1)}{\frac{(n+1-p)(n+1)}{n+1}+\frac{(n+1-p)\beta_1}{n+1}}=\frac{p(n+1)}{n+1-p}.\]
\begin{rem}
We give the following remark about exponents in Theorem \ref{Weighted-Sobolev-1}.

\begin{enumerate}
\item In fact,  the  exponent conditions in Theorem \ref{Weighted-Sobolev-1} is equivalent to $n\ge1, \beta_1>-1, \alpha_1\ge\frac{p\beta_1}q, 1<p<n+1,  p<q=p^* $ with
  \[\frac{n+1+\beta_1}q=\frac{n+1+\alpha_1-p}p.
  \]
  \item  If $p=q=p^*,$ we have no inequality \eqref{Weighted-Sobolev-1}. Because  of the restriction $\alpha<\frac1{p'}$  in inequality \eqref{WS-1}, it implies that $\alpha_1<p-1$.  But it follows from  \eqref{alpha1 fomula} that $\alpha_1>p-1$  when $p=q$. So we have no inequality  in this case. I also checked that the inequality (2.1.35)(in Maz'ya book \cite{Maz2011}) assume condition $\alpha>\frac1{p'}$, which implies $\alpha_1>p-1$ for $p=q.$
 \end{enumerate}
 \end{rem} Define the Green function of fractional Laplacian $(-\Delta)^\frac\gamma2$ on $\mathbb{R}^{n+1}_+$ as
 \begin{eqnarray*}
G_\infty(X,Y)&=&\frac{A_{n,\gamma}}{|X-Y|^{n+1-\gamma}}\int_0^{\frac{4tz}{|X-Y|^2}}
\frac{b^{\frac\gamma2-1}}{(1+b)^\frac{n+1}2}db,
\end{eqnarray*}
%\begin{eqnarray*}
%G_\infty(X,Y)&=&\frac{A_{n,\gamma}}{|X-Y|^{n+1-\gamma}}\big[1
%-B\frac{1}{(4tz+|X-Y|^2)^\frac{n-1}2}\\
%\quad&&\times
%\int_0^{\frac{|X-Y|^2}{4tz}}\frac{(|X-Y|^2-4tzb)^\frac{n-1}2}{b^\frac\gamma2(1+b)}db\big],
%\end{eqnarray*}
where $A_{n,\gamma}$ is some positive constant (see e.g.,  \cite{CFY2015, Chen2021}).
Moreover, it can be verified that
\begin{equation}\label{Fractional eq-2}
G_\infty(X,Y)\le \frac{A_{n,\gamma}}{|X-Y|^{n+1-\gamma}}.
\end{equation}
Chen, Li and Ma  \cite{Chen2021} proved that if $u\in L_\gamma:=\big\{u\in L^1_{loc} \,|\, \int_{\mathbb{R}^{n+1}}\frac{|u(X)|}{1+|X|^{n+1+\gamma}}dX\big\}$ is a locally bounded positive solution of
\begin{equation}\label{Fractional eq-1}
\begin{cases}
(-\Delta)^\frac\gamma2u(X)=f(X), \quad& X\in\mathbb{R}^{n+1}_+,\\
u=0,\quad& X\not\in\mathbb{R}^{n+1}_+
\end{cases}
\end{equation} in the sense of distribution, then  it is also a solution of the integral equation
\begin{equation*}
u(X)=\int_{\mathbb{R}^{n+1}_+}G_\infty(X,Y)f(Y)dY.
\end{equation*}

\textbf{Proof of Theorem \ref{Fractional Sobolev-1}.} From \eqref{Fractional eq-2} we have
\begin{eqnarray*}
u(X)&=&\int_{\mathbb{R}^{n+1}_+}G_\infty(X,Y)f(Y)dY\\
&\le& C_{n,\gamma}\int_{\mathbb{R}^{n+1}_+}\frac {f(Y)}{|X-Y|^{n+1-\gamma}}dY\\
&=& C_{n,\gamma}E_{\lambda} f
\end{eqnarray*}  with $\lambda=n+1-\gamma$.
Furthermore, it follows from \eqref{Fractional eq-1} that
 \[u(X)\le C_{n,\gamma}E_\lambda((-\Delta)^\frac\gamma2u(X)).\]
Substituting the above into the  inequality \eqref{WHLSD-1}, there exists some  constant $S(n+1, p,\alpha,\beta, \gamma)$ such that
\[\|t^{-\beta}u\|_{L^q({\mathbb{R}^{n+1}_+})}\le S(n+1, p,\alpha,\beta, \gamma)\|t^\alpha(-\Delta)^\frac\gamma2u\|_{L^p(\mathbb{R}^{n+1}_+)},\] for $u\in W_\alpha^{\gamma, p}(\mathbb{R}^{n+1}_+)\cap L_\gamma$ with $0<\gamma<n+1.$
\qed

%\begin{theorem}\label{Fractional Sobolev-nm}

\textbf{Proof of Theorem \ref{Fractional Sobolev-nm}.}
Recall
\[I_\lambda f(y,\hat{y})=\int_{\mathbb{R}^{n+m}}\frac{f(x,\hat{x})}{|(x,\hat{x})-(y,\hat{y})|^\lambda}dxd\hat{x},  ~~~~\forall (y,\hat{y})\in\mathbb{R}^{n+m}, \,  \lambda \in (0,n+m).\]
Let $u=I_\lambda f$ with $\lambda=n+m-\gamma, \, 0<\gamma<n+m.$ Chen, Li and Ou \cite{CLO2006} (also see \cite{CL2010}) showed that $u$ satisfies
\begin{equation}\label{F-E}
(-\Delta)^\frac\gamma2u=f,\quad\text{in}~\mathbb{R}^{n+m}.
\end{equation}
%Invoking \eqref{EWHLSD-1}, we also have the following weighted Sobolev inequality for fractional Laplacian $(-\Delta)^\frac\gamma2.$

Using $u=I_\lambda f $ and the Stein-Weiss type inequality  \eqref{EWHLSD-1}, it yields
\[\||\hat{y}|^{-\beta}u\|_{L^q( \mathbb{R}^{n+m})}\le C_{\alpha,\beta,\lambda,p}\||\hat{y}|^\alpha f\|_{L^p(\mathbb{R}^{n+m})}.\]
Combining the above and \eqref{F-E}, there exists some constant $S(n,m, p,\alpha,\beta,\gamma)$ such that
\[\||\hat{y}|^{-\beta}u\|_{L^q({\mathbb{R}^{n+m}})}\le S(n,m, p,\alpha,\beta,\gamma)\||\hat{y}|^\alpha(-\Delta)^\frac\gamma2u\|_{L^p(\mathbb{R}^{n+m})}\]
for any $u\in W_\alpha^{\gamma, p}(\mathbb{R}^{n+m})$ with $0<\gamma<n+m.$

Similarly, choosing $\lambda=n+m-1,$ and employing \eqref{Tr-0} on $\mathbb{R}^{n+m}$,  we can show inequality \eqref{Sobolev-inq-mn} holds.
\qed

\medskip

 In the final, for the conformal case $p= p_\alpha=\frac{2(n+1)}{2(n+1)-\lambda-2\alpha}, \, r=r_\beta=\frac{2(n+1)}{2(n+1)-\lambda-2\beta}$, we show the Stein-Weiss type inequality \eqref{WHLSD-O} is equivalent  to  the  following inequality on the unit ball by Kelvin transformation. That is,
 \begin{crl}\label{cor-ball}
 Let $0<\lambda<n+1,$ $p=p_\alpha,r=r_\beta,$ $
\alpha<\frac{p_\alpha-1}{p_\alpha}, \, \beta<\frac {r_\beta-1}{r_\beta}, \,  \alpha+\beta\ge0$, and $x^1=(0,\cdots,0,-1)\in\mathbb{R}^{n+1}.$
There is a constant $N_{\alpha,\beta,\lambda,p_\alpha}:=C(n,\alpha,\beta,\lambda,p_\alpha)>0$, such that for all $F\in L^{p_\alpha}(B^{n+1}), \, G\in L^{r_\beta}(B^{n+1}),$
\begin{eqnarray}\label{WHLSBB}
\big|\int_{B^{n+1}}\int_{B^{n+1}}
\big(\frac 2{1-|\xi-x^1|^2}\big)^{\alpha}\big(\frac 2{1-|\eta-x^1|^2}\big)^{\beta}
\frac{F(\xi)G(\eta)}{|\xi-\eta|^\lambda}d\xi d\eta\big|\nonumber\\
\le N_{\alpha,\beta,\lambda,p_\alpha}\|F\|_{L^{p_\alpha}(B^{n+1})}\|G\|_{L^{r_\beta}(B^{n+1})}.
\end{eqnarray}
\end{crl}
\proof \,
Let $x^0=(0,\cdots,0,-2), \, x^1=(0,\cdots,0,-1)\in\mathbb{R}^{n+1},$ write $B^{n+1}=B_1(x^1)\subset\mathbb{R}^{n+1}.$ Define Kelvin transformation as
\[\eta^{x^0,2}=\frac{4(\eta-x^0)}{|\eta-x^0|^2}+x^0\in B^{n+1}\subset\mathbb{R}^{n+1}, \quad \eta\in\mathbb{R}^{n+1}_+,\]
and
\begin{equation*}
F(\xi)=f_{x^0,2}(\xi)=\big(\frac2{|\xi-x^0|}\big)^{\mu_1}f(\xi^{x^0,2}), \quad
G(\eta)=g_{x^0,2}(\eta)=\big(\frac2{|\eta-x^0|}\big)^{\mu_2}g(\eta^{x^0,2}),
\end{equation*} where $\mu_1=2(n+1)-\lambda-2\alpha, \, \mu_2=2(n+1)-\lambda-2\beta.$

It is easy to check the fact
\[\eta_{n+1}^{x^0,2}=\frac{2^2(\eta_{n+1}+2)}{|\eta-x^0|^2}-2=\frac4{|\eta-x^0|^2}\frac
{1-|\eta-x^1|^2}2.\]
If $p=p_\alpha, \, r=r_\beta,$ using Kelvin transform, rewrite  the left of the inequality \eqref{WHLSD-O} as
\begin{eqnarray}\label{EFB}
&&\int_{\mathbb{R}^{n+1}_+}\int_{\mathbb{R}^{n+1}_+}\frac{f(X) g(Y)}{t^\alpha|X-Y|^{\lambda} z^\beta} dXdY\nonumber\\
&=&\int_{B^{n+1}}\int_{B^{n+1}}\big(\frac2{|\xi-x_0|}
\cdot\frac2{|\eta-x_0|}\big)^{2(n+1)}\frac{f(\xi^{x^0,2})g(\eta^{x^0,2})}{(\xi_{n+1}^{x^0,2})^\alpha
(\eta_{n+1}^{x^0,2})^\beta|\xi^{x^0,2}-\eta^{x^0,2}|^\lambda}d\xi d\eta\nonumber\\
&=&\int_{B^{n+1}}\int_{B^{n+1}}\big(\frac{1-|\xi-x^1|}2\big)^{-\alpha}
\big(\frac{1-|\eta-x^1|}2\big)^{-\beta}\frac{F(\xi)G(\eta)}{|\xi-\eta|^\lambda}\nonumber\\
&\times&(\frac2{|\xi-x_0|})^{2(n+1)-\lambda-2\alpha-\mu_1}
(\frac2{|\eta-x_0|})^{2(n+1)-\lambda-2\beta-\mu_2}d\xi d\eta\nonumber\\
&=&\int_{B^{n+1}}\int_{B^{n+1}}\big(\frac 2{1-|\xi-x^1|}\big)^{\alpha}\big(\frac 2{1-|\eta-x^1|}\big)^{\beta}
\frac{F(\xi)G(\eta)}{|\xi-\eta|^\lambda}d\xi d\eta.
\end{eqnarray} A direct computation has
\begin{eqnarray}\label{FB}
\|f\|^{p_\alpha}_{L^{p_\alpha}(\mathbb{R}^{n+1}_+)}
&=&\int_{\mathbb{R}^{n+1}_+}|f(X)|^{p_\alpha}dX\nonumber\\
&=&\int_{B_1(x_1)}|f(\xi^{x^0,2})|^{p_\alpha}\big(\frac2{|\xi-x^0|}\big)^{2(n+1)}
d\xi\nonumber\\
&=&\int_{B_1(x_1)}\big|\big(\frac2{|\xi-x_0|}\big)^{-\mu_1}F(\xi)\big|^{p_\alpha}
\big(\frac2{|\xi-x_0|}\big)^{2(n+1)}d\xi\nonumber\\
%&=&\int_{B_1(x_1)}\big(\frac2{|\xi-x_0|}\big)^{-2(n+1)}\big(\frac2{|\xi-x_0|}\big)^{2(n+1)}|F(\xi)|^{p_\alpha}d\xi\nonumber\\
&=&\|F\|^{p_\alpha}_{L^{p_\alpha}(B^{n+1})}.
\end{eqnarray} Similarly,
\begin{equation}\label{GB}
\|g\|_{L^{r_\beta}(\mathbb{R}^{n+1}_+)}=\|G\|_{L^{r_\beta}(B^{n+1})}.
\end{equation}
Combining \eqref{EFB}, \eqref{FB} and \eqref{GB}, we  obtain inequality \eqref{WHLSBB}.
\qed

 \vskip 1cm
{\noindent {\bf Acknowledgements}\\
%The author would like to thank the referee for his/her careful reading of the manuscript and many good suggestions.
The author would like to thank Prof. Meijun Zhu for very helpful discussions and comments on this research. The project is supported by  the National Natural Science Foundation of China (Grant No. 12071269), Youth Innovation Team of Shaanxi Universities, Shaanxi Fundamental Science Research
Project for Mathematics and Physics (Grant No. 22JSZ012) and the Fundamental Research Funds for the Central Universities (Grant No. 2021CBLY001, GK202307001, GK202202007).

\medskip	

\noindent {\bf Declarations}\\
	
We declare that we have no conflict of interest.

\medskip	

\textbf{Data availability }The manuscript has no associated data.

	%% bibliography--------------------------------------------------------------------
%\begin{center}

\end{document}